\documentclass[10pt,oneside]{book}
\usepackage{amsmath,amssymb,color,graphicx,tikz,MnSymbol,stmaryrd}
\usetikzlibrary{positioning}
\usepackage[mathcal]{euscript}
\usetikzlibrary{matrix,arrows}
\DeclareMathAlphabet{\mathpzc}{OT1}{pzc}{m}{it}
\usepackage{dsfont,url}
\usepackage{bbm}
\usepackage[latin1]{inputenc}
\usepackage[russian,english]{babel}
\usepackage[textsize=footnotesize]{todonotes}
\usepackage{mathrsfs}
\usepackage[OT2,T1]{fontenc}
\usepackage{andsupp}
\usepackage{amscd}
\usepackage{boxedminipage}
\usepackage{cancel}
\usepackage[line,cmtip,all,matrix,dvips]{xy}

\makeatletter
\renewcommand{\tableofcontents}{%
  \@starttoc{toc}%
}
\makeatother

\newtheorem{theorem}{Theorem}[section]

\newtheorem{proposition}[theorem]{Proposition}
\newtheorem{lemma}[theorem]{Lemma}
\newtheorem{corollary}[theorem]{Corollary}

\newtheorem{definition}[theorem]{Definition}

\newtheorem{remark}[theorem]{Remark}

\def\beaq{\begin{eqnarray}}
\def\eeaq{\end{eqnarray}}

\newcommand{\beq}{\begin{equation}}
\newcommand{\eeq}{\end{equation}}
\newcommand{\bea}{\begin{eqnarray}}
\newcommand{\eea}{\end{eqnarray}}
\newcommand{\dd}{\mathrm{d}}

\newcommand{\Res}{\mathop{\,\rm Res\,}}

\definecolor{rouge}{rgb}{0.84,0.18,0.07}
\definecolor{bleu}{rgb}{0.22,0.41,0.74}
\definecolor{vertf}{rgb}{0.08,0.46,0.07}
\usepackage[pdfpagelabels]{hyperref}
\hypersetup{colorlinks,urlcolor=vertf,citecolor=rouge,linkcolor=bleu,filecolor=black}

\renewcommand{\epsilon}{\varepsilon}
\renewcommand{\rho}{\varrho}
\renewcommand{\phi}{\varphi}

\newcommand{\Surf}{{\mathcal{B}}}

\newcommand{\Bosonnormalord}{\stackrel{\textstyle \circ}{ \circ}}

\newcommand{\bt}{\begin{theorem}}
\newcommand{\et}{\end{theorem}}
\newcommand{\br}{\begin{remark}}
\newcommand{\er}{\end{remark}}

\renewcommand\thechapter{\Roman{section}}

\def\Bosonnormalord{\,\lower.8ex \hbox{$\circ$} \llap{\raise.8ex\hbox{$\circ$}} \,}
\def\normalord{\,\lower.8ex \hbox{$\cdot$} \llap{\raise.8ex\hbox{$\cdot$}} \,}

\def\ra{\rightarrow}

\def\Res{\mathop{\rm Res}\limits}

\addto\captionsenglish{%
  }

\begin{document}

\sloppy

\pagestyle{empty}
\addtolength{\baselineskip}{0.20\baselineskip}
\begin{center}

\vspace{26pt}

{\huge Geometric recursion}

\vspace{26pt}

\textsl{J\o{}rgen Ellegaard Andersen}\footnote{Centre for Quantum Mathematics, Danish Institute for Advanced Study, University of Southern Denmark, Campusvej 55, 5230 Odense M, Denmark. \\ \href{mailto:jea@sdu.dk}{\texttt{jea@sdu.dk}}}, \textsl{Ga\"etan Borot}\footnote{Max Planck Institut f\"ur Mathematik, Vivatsgasse 7, 53111 Bonn, Germany. \\ Humboldt-Universit\"at zu Berlin, Institut f\"ur Mathematik und Institut f\"ur Physik, Rudower Chaussee 25, 10247 Berlin, Germany, \href{mailto:gaetan.borot@hu-berlin.de}{\texttt{gaetan.borot@hu-berlin.de}}}, \textsl{Nicolas Orantin}\footnote{\'{E}cole Polytechnique F\'{e}d\'{e}rale de Lausanne, D\'epartement de Math\'ematiques, 1015 Lausanne, Switzerland. \\  Section de Math\'ematiques, Universit\'e de Gen\`eve, Uni Dufour, 24, rue du Général Dufour, Case postale 64, 1211 Genève 4, Switzerland. \\ \href{mailto:nicolas.orantin@unige.ch}{\texttt{nicolas.orantin@unige.ch}}}
\end{center}

\vspace{20pt}

%

\begin{center}
\textsc{Abstract}
\end{center}

We propose a general theory for constructing functorial assignments
\[
\Sigma \longmapsto \Omega_{\Sigma} \in E(\Sigma)
\]
for a large class of functors $E$ from a certain category of bordered surfaces to a suitable target category of topological vector spaces. The construction proceeds by successive excisions of homotopy classes of embedded pairs of pants, and thus by induction on the Euler characteristic. We provide sufficient conditions to guarantee the infinite sums appearing in this construction converge. In particular, we can generate mapping class group invariant vectors $\Omega_{\Sigma} \in E(\Sigma)$. The initial data for the recursion encode the cases when $\Sigma$ is a pair of pants or a torus with one boundary, as well as the ``recursion kernels'' used for glueing. We give this construction the name of Geometric Recursion.

As a first application, we demonstrate that our formalism produce a large class of measurable functions on the moduli space of bordered Riemann surfaces. Under certain conditions, the functions produced by the geometric recursion can be integrated with respect to the Weil--Petersson measure on moduli spaces with fixed boundary lengths, and we show that the integrals satisfy a topological recursion generalizing the one of Eynard and Orantin. We establish a generalization of Mirzakhani--McShane identities, namely that multiplicative statistics of hyperbolic lengths of multicurves can be computed by the geometric recursion, and thus their integrals satisfy the topological recursion. As a corollary, we show that the systole function can be obtained from the geometric recursion, and find an interpretation of the intersection indices of the Chern character of bundles of conformal blocks in terms of the aforementioned statistics.

The theory has however a wider scope than functions on Teichm\"uller space, which will be explored in subsequent papers;  one expects that many functorial objects in low-dimensional geometry could be constructed by variants of this geometric recursion.

\vfill

\newpage

\vspace{26pt}
\pagestyle{plain}
\vfill 

\vspace{0.5cm}

\newpage

\setcounter{tocdepth}{2}
\tableofcontents
\newpage

\section{Introduction}

This article grew out of the search for an intrinsic geometric meaning of the topological recursion of \cite{EORev}, inspired by the famous but isolated example of Mirzakhani--McShane identities and Mirzakhani's recursion for the Weil--Petersson volume of the moduli spaces of bordered Riemann surfaces \cite{Mirzakhani}. As an outcome, we devise a general formalism which constructs mapping class group invariant objects associated to surfaces of arbitrary topology. In order to speak of mapping class group invariant objects, we need to begin with spaces carrying actions of the mapping class groups. For the purposes of this article, this takes the form of a functor from a category of smooth surfaces to a category of topological vector spaces. Our construction also demands functorial maps that realise the disjoint union and glueing of surfaces. The functor together with these maps is, roughly speaking, what we will call a ``target theory'' (see Section \ref{PRETARGETDEF} for the precise definition), and geometry provides a wide range of target theories. Given a target theory, and given initial data which are mapping class group invariant vectors attached to pairs of pants and tori with one boundary, our construction yields mapping class group invariant vectors for arbitrary smooth surfaces by successive excisions of embedded pairs of pants, i.e. by induction on the Euler characteristic. In this text, we establish the foundations of this construction and explore its first applications. 

\subsection{Main construction}

In the interest of keeping the introduction short, yet covering our main results, we provide sketches of definitions of needed concepts here and refer the reader to the main text for the precise definitions. We work in the category $\Surf_{1}$ where
\begin{itemize}
\item[$\bullet$] objects are stable bordered, compact, oriented smooth surfaces $\Sigma$ with exactly one negatively oriented boundary per connected component (denoted $\partial_-\Sigma$), with non-empty boundary if the surface is non-empty. The stability condition means that the Euler characteristic of each connected component of $\Sigma$ must be negative (no disks or cylinders are allowed).
\item[$\bullet$] morphisms are isotopy classes of diffeomorphisms of the surface preserving the orientation of the surface and the given orientations of the boundaries.
\end{itemize}
The automorphism group of an object $\Sigma$ is denoted $\Gamma_{\Sigma}$ and the subgroup consisting of mapping classes which induce the identity permutation on the set of boundary components is denoted $\Gamma^{\partial}_{\Sigma}$. We shall denote $\partial_-\Sigma$ and $\partial_{+}\Sigma$ the set of boundaries whose orientation (dis)agrees with the one of $\Sigma$. We will consider functors $E\,:\,\Surf_1 \rightarrow \mathcal{C}$ to the category $\mathcal{C} = {\rm Pro}$-$\mathcal{V}$ of projective systems of objects in the category $\mathcal{V}$ of Hausdorff complete locally convex topological vector spaces over $\mathbb{K} = \mathbb{R}$ or $\mathbb{C}$ (Section~\ref{STarget}). There is a particular subspace $E'(\Sigma)$ of the projective limit in which our constructed mapping class group invariant vectors will be contained.

In $\Surf_1$ we can take disjoint unions and we have the operations of cutting along primitive multicurves $\gamma$ with ordered components (see Section~\ref{S1}). We introduce the notion of \emph{target theory} in Section~\ref{STargetT}. Roughly speaking, it is the data of a functor $E\,:\, \Surf_1 \rightarrow \mathcal{C}$ together with a collection of functorial multilinear morphisms $\sqcup$ and $\Theta$ representing the disjoint union and certain cutting operations, and a system of abstract ``length functions'' on the set of simple closed curves. 

For connected $\Sigma$, we will be especially interested in sets $\mathcal{P}_{\Sigma}^{\emptyset}$ (resp. $\mathcal{P}_{\Sigma}^{b}$) of homotopy classes of embedded pairs of pants $[P \hookrightarrow \Sigma]$ such that $\partial P \cap \partial \Sigma = \partial_-\Sigma$  (resp. $\partial P \cap \partial \Sigma = \partial_-\Sigma \cap b$ for the indicated boundary component $b$), see Section~\ref{glueinpant} for their definition. We shall denote $\Theta_{P}\,:\,E(P) \times E(\Sigma - P) \rightarrow E(\Sigma)$ the bilinear morphism that corresponds to reglueing along $\partial P \cap \Sigma$ in the given target theory.

\vspace{0.2cm}

Our main result is the existence of the following construction.

\begin{theorem}[Definition]
\label{thM} Let $E:\,\Surf_{1} \rightarrow \mathcal{C}$ be a target theory. Let $P$ (resp. $T$) be an object of $\Surf_{1}$ with the topology of a pair of pants (resp. a torus with one boundary), and $b$ be the choice of a positively oriented boundary of $P$. Assume we are given initial data
\[
A_{P} \in E'(P)^{\Gamma_{P}},\qquad B_{P}^{b} \in E(P)^{\Gamma^{\partial}_{P}},\qquad C_{P} \in E(P)^{\Gamma_{P}},\qquad D_{T} \in E'(T)^{\Gamma_{T}},
\] 
which are admissible according to Definition~\ref{DEFDECAY}. There exists a unique, well-defined, functorial assignment
\[
\Sigma \mapsto \Omega_{\Sigma} \in E'(\Sigma)
\]
specified by the formulas
\[
\Omega_{\emptyset} = 1_{E(\emptyset)},\qquad \Omega_{P} = A_{P},\qquad \Omega_{T} = D_{T},\qquad \Omega_{\Sigma_1 \cup \Sigma_{2}} = \Omega_{\Sigma_1} \sqcup \Omega_{\Sigma_{2}},
\]
and for any connected $\Sigma$ of Euler characteristic $\leq -2$
\begin{equation}
\label{firstseries}\Omega_{\Sigma} =  \sum_{b \in \pi_0(\partial_+\Sigma)} \sum_{[P] \in \mathcal{P}^{b}_{\Sigma}} \Theta_{P}(B_{P}^{b},\Omega_{\Sigma - P})+ \sum_{[P] \in \mathcal{P}^{\emptyset}_{\Sigma}} \tfrac{1}{2}\,\Theta_{P}(C_{P},\Omega_{\Sigma - P}).
\eeq 
\end{theorem}

The admissibility assumption of initial data includes decay in $B$ and $C$ with respect to the length functions provided by the target theory, which permits us to prove that the series \eqref{firstseries} are absolutely convergent with respect to semi-norms provided and finally is contained in $E'(\Sigma)$. This implies functoriality, as the action of isotopy classes of diffeomorphisms merely permutes the terms in \eqref{firstseries}, in the strong sense below.

\begin{theorem}[Naturality]  \label{natTh} Let $\eta:\,E \Rightarrow \tilde{E}$ be a natural transformation between two target theories $E$ and $\tilde{E}$. Let $(A,B,C,D)$ be an admissible initial data for $E$, leading by the above construction to the functorial assignment $\Omega$. Then $(\eta(A),\eta(B),\eta(C),\eta(D))$ is an admissible data for $\tilde{E}$, and by the above construction it leads to a functorial assignment equal to $\eta(\Omega)$.
\end{theorem}

Notice that an $E$-valued functorial assignment is nothing but a natural transformation $E \Rightarrow \mathbf{1}$ to the trivial functor $\mathbf{1}\,:\,\Sigma \mapsto \mathbb{K}$, so that Theorem~\ref{natTh} expresses the compatibility of $\Omega\,:\,E \Rightarrow \mathbf{1}$ with the pre-composition of natural transformations.

\vspace{0.2cm}

We give the name \emph{geometric recursion} (GR) to our construction, because many of the target theories we have in mind come from spaces of geometric structures on surfaces. We refer to the $\Omega_{\Sigma}$s as ``GR amplitudes''.

\subsection{Functions on Teichm\"uller spaces}

A fundamental example of a target theory, which is developed in Section~\ref{SfunSctin}, is the space $E(\Sigma) = {\rm Mes}(\mathcal{T}_{\Sigma})$ of measurable functions on the Teichm\"uller space $\mathcal{T}_{\Sigma}$ of the bordered surface $\Sigma$. We equip ${\rm Mes}(\mathcal{T}_{\Sigma})$ with the topology of convergence on every compact subset, and this space is seen as the limit of a projective system $({\rm Mes}(\mathcal{T}_{\Sigma}^{\epsilon}))_{\epsilon \in (0,1)}$ where $\mathcal{T}_{\Sigma}^{\epsilon} \subset \mathcal{T}_{\Sigma}$ is the $\varepsilon$-thick part. In this case $E'(\Sigma)$ is the subspace of $E(\Sigma) = {\rm Mes}(\mathcal{T}_{\Sigma})$ consisting of measurable functions which are uniformly bounded on each systole set $\mathcal{T}_{\Sigma}^{\epsilon}$. The cutting morphism relies on the cutting of hyperbolic structures along geodesics, and the length functions in the target theory are given by hyperbolic lengths $\ell_{\sigma}$ for $\sigma \in \mathcal{T}_{\Sigma}$.

In this context, mapping class group invariant measurable functions are simply just measurable functions defined on moduli spaces $\mathcal{M}_{\Sigma} = \mathcal{T}_{\Sigma}/\Gamma^{\partial}_{\Sigma}$ of bordered Riemann surfaces, and Theorem~\ref{thM} specializes to the following result.

\begin{corollary}
\label{cothbtM}Let $(A,B,C)$ be measurable functions on $\mathbb{R}_{+}^3$ and $D_T$ be a $\Gamma_{T}$-invariant measurable function on $\mathcal{T}_T$. We assume that $X(L_1,L_2,L_3) = X(L_1,L_3,L_2)$ for $X \in \{A,C\}$, and the existence of $t \geq 0$ and for any $s,\epsilon > 0$ the existence of $M_{s,\epsilon} > 0$ such that for any $L_1,L_2,\ell,\ell' > \epsilon$ and $\sigma \in \mathcal{T}_{T}^{\epsilon}$
\begin{equation}
\begin{split}
|A(L_1,L_2,L_3)| & \leq M_{\epsilon,1}(1 + L_1)^{t}(1 + L_2)^{t}(1 + L_3)^{t}, \\
|B(L_1,L_2,\ell)| & \leq \frac{M_{\epsilon,s}(1 + L_1)^{t}(1 + L_2)^{t}}{\big(1 + [\ell - L_1 - L_2]_{+}\big)^{s}}, \\
|C(L_1,\ell,\ell')| & \leq \frac{M_{\epsilon,s}(1 + L_1)^{t}}{\big(1 + [\ell + \ell' - L_1]_{+}\big)^{s}}, \\
|D_{T}(\sigma)| & \leq M_{\epsilon,1}\,(1 + \ell_{\sigma}(\partial T))^{t}.
\end{split}
\end{equation}
where $[x]_{+} = \max(x,0)$. We have a well-defined functorial assignment $\Sigma \mapsto \Omega_{\Sigma} \in {\rm Mes}(\mathcal{T}_{\Sigma})$ given by the formulas
\[
\Omega_{P}(\sigma) = A(\vec{\ell}_{\sigma}(\partial P)),\qquad \Omega_{T}(\sigma) = D_{T}(\sigma),  \qquad \Omega_{\Sigma_1 \cup \Sigma_2}(\sigma) = \Omega_{\Sigma_1}(\sigma|_{\Sigma_1}) \Omega_{\Sigma_2}(\sigma|_{\Sigma_2}),
\]
and for connected objects $\Sigma$ of Euler characteristic $\leq -2$
\begin{equation}
\label{fiunfgnGRseris}
\Omega_{\Sigma}(\sigma) = \sum_{b \in \pi_0(\partial_+\Sigma)} \sum_{[P] \in \mathcal{P}_{\Sigma}^{b}} B(\vec{\ell}_{\sigma}(\partial P))\,\Omega_{\Sigma - P}(\sigma|_{\Sigma - P}) + \frac{1}{2} \sum_{[P] \in \mathcal{P}_{\Sigma}^{\emptyset}} C(\vec{\ell}_{\sigma}(\partial P))\,\Omega_{\Sigma - P}(\sigma|_{\Sigma - P}).
\end{equation}
More precisely, the series \eqref{fiunfgnGRseris} converges absolutely, uniformly for $\sigma$ in any compact of $\mathcal{T}_{\Sigma}$.
\end{corollary}

When $\Omega$ is given by Corollary~\ref{cothbtM} and $L \in \mathbb{R}_{+}^n$, we denote by $\Omega_{g,n}$ the function it induces on $\mathcal{M}_{g,n}(L) \cong \mathcal{M}_{\Sigma}(L)$ where $\Sigma$ is a bordered surface of genus $g$ with $n$ ordered boundary components such that $\partial_-\Sigma$ is the first.  The following result -- stated with stronger assumptions than the result presented in Theorem \ref{intTRHTT} -- provides sufficient conditions for the integrability of GR amplitudes with respect to the Weil--Petersson measure $\mu_{{\rm WP}}$ and computes these integrals 
\[
V\Omega_{g,n}(L) := \int_{\mathcal{M}_{g,n}(L)} \Omega_{g,n}\dd\mu_{{\rm WP}}.
\]

\begin{theorem}
\label{ntinghg}In the situation of Corollary~\ref{cothbtM}, assume that one can choose $M_{s,\epsilon}$ independently of $\epsilon$. Then, for $2g - 2 + n > 0$, the functions $\Omega_{g,n}$ coming from Corollary~\ref{cothbtM} are integrable over $\mathcal{M}_{g,n}(L)$, for each $L \in \mathbb{R}_{+}$. They satisfy the topological recursion, for $2g - 2 + n \geq 2$
\begin{equation}
\label{TRintross} \begin{split}
& \quad V\Omega_{g,n}(L_1,\ldots,L_n) \\
& = \sum_{m = 2}^{n} \int_{\mathbb{R}_{+}}\,B(L_1,L_m,\ell) V\Omega_{g,n - 1}(\ell,L_2,\ldots,\widehat{L_m},\ldots,L_n)  \ell \,\dd \ell\\
& \quad +\frac{1}{2} \int_{\mathbb{R}_{+}^2} C(L_1,\ell,\ell') \bigg(V\Omega_{g-1,n + 1}(\ell,\ell',L_2,\ldots,L_n) + \!\!\!\!\!\!\!\!\! \sum_{\substack{J \cup J' = \{L_2,\ldots,L_n\} \\  h + h ' = g}}\!\!\!\!\!\!\!\!\!\!\!\!\! V\Omega_{h,1 + \#J}(\ell,J) V\Omega_{h',1 + \#J'}(\ell',J')\bigg) \ell\ell'  \, \dd \ell \dd \ell' ,
\end{split}
\end{equation}
with the base cases
\[
V\Omega_{0,3}(L_1,L_2,L_3) = A(L_1,L_2,L_3),\qquad V\Omega_{1,1} := V\!D(L_1) = \int_{\mathcal{M}_{T}(L_1)} D_T\,\dd\mu_{{\rm WP}}\,.
\]
\end{theorem}

The geometric recursion constructs a large class of functions on the moduli spaces $\mathcal{M}_{\Sigma}$, namely those satisfying the ``non-local glueing rule'' \eqref{fiunfgnGRseris} and they are designed such that their integrals against $\mu_{{\rm WP}}$ satisfies topological recursion. As we will now show, there are many geometrically meaningful functions on the moduli spaces which belong to this class. 

The first example, which in fact was an inspiration for the whole formalism, is the constant function $1$ on $\mathcal{T}_{\Sigma}$, as a result of Mirzakhani's generalisation of McShane identities.
\begin{theorem} \cite{Mirzakhani}\label{MMGR}
\label{Mirzathth} The geometric recursion for the initial data $A^{{\rm M}} = 1$, $D^{{\rm M}}_{T} = 1$ ,
\begin{equation*}
\begin{split}
B^{{\rm M}}(L_1,L_2,\ell) & = 1 - \frac{1}{L_1}\ln\bigg(\frac{{\rm cosh}(L_2/2) + {\rm cosh}[(L_1 + \ell)/2]}{{\rm cosh}(L_2/2) + {\rm cosh}[(L_1 - \ell)/2]}\bigg) \\
C^{{\rm M}}(L_1,\ell,\ell') & = \frac{1}{L_1}\ln\bigg(\frac{e^{L_1/2} + e^{(\ell + \ell')/2}}{e^{-L_1/2} + e^{(\ell + \ell')/2}}\bigg)
\end{split}
\end{equation*}
yields the function $\Omega_{\Sigma} = 1$ for all $\Sigma$.
\end{theorem}

We establish a generalisation of Mirzakhani's identities by showing that for a twist of the initial data, the GR amplitudes compute multiplicative statistics of the hyperbolic spectrum of multicurves. Let $M_{\Sigma}'$ be the set of primitive multicurves (including the empty one) in $\Sigma$. If $c \in \Sigma$, we denote $\Sigma^{c}$ the surface $\Sigma$ cut along $c$.

\begin{theorem}\label{Twistf}
Let $(A,B,C,D)$ satisfy the conditions of Corollary~\ref{cothbtM} and assume that the corresponding GR amplitudes $\Omega_{\Sigma}$ are invariant under braiding of all boundary components of $\Sigma$. Let $f$ be a measurable function on $\mathbb{R}_{+}$ such that $\sup_{\ell > 0} (1 + \ell)^{s}|f(\ell)| < +\infty$ for any $s > 0$. Then, the series
\[
\Omega_{\Sigma}[f](\sigma) := \sum_{c \in M'_{\Sigma}} \Omega_{\Sigma^c}(\sigma|_{\Sigma^{c}})\,\prod_{\gamma \in \pi_0(c)} f(\ell_{\sigma}(\gamma))
\]
converges absolutely, uniformly for $\sigma$ in any compact of $\mathcal{T}_{\Sigma}$, and it coincides with the GR amplitudes for the initial data
\begin{equation*}
\begin{split}
A[f](L_1,L_2,L_3) & = A(L_1,L_2,L_3) \\
B[f](L_1,L_2,\ell) & = B(L_1,L_2,\ell) + A(L_1,L_2,\ell)\,f(\ell) \\
C[f](L_1,\ell,\ell') & = C(L_1,\ell,\ell') + B(L_1,\ell,\ell')f(\ell) + B(L_1,\ell',\ell)f(\ell') + A(L_1,\ell,\ell')f(\ell)f(\ell') \\
D_{T}[f](\sigma) & =  D_{T}(\sigma) + \frac{1}{2} \sum_{\gamma \in S_{T}^{\circ}} A(\ell_{\sigma}(\partial T),\ell_{\sigma}(\gamma),\ell_{\sigma}(\gamma))\,f(\ell_{\sigma}(\gamma))
\end{split}
\end{equation*}
where $S_{T}^{\circ}$ is the set of simple closed curves in $T$.
\end{theorem}
The integrals of $\Omega_{g,n}[f]$ over $\mathcal{M}_{g,n}(L)$ against the Weil--Petersson measure can then be computed by the topological recursion from Theorem~\ref{ntinghg}, or by direct integration as a sum over stable graphs $\mathcal{G}_{g,n} \cong M'_{\Sigma}/\Gamma_{\Sigma}^{\partial}$ (see Lemma~\ref{Lemstab}). After integration, the twisting has some relation with the Givental group action which we discuss in Section~\ref{GUNIGUNFIOGSO}. 

\vspace{0.2cm}

\subsection{Applications and perspectives}

The formalism of this article has already found several applications.

\medskip

Theorem~\ref{Twistf} was used in \cite{MVpaper} to prove that the Masur--Veech volumes of the (top) stratum of the moduli space of quadratic differentials $\mathcal{Q}\mathfrak{M}_{g,n}$ with $4g - 4 + n$ zeroes and $n$ poles, can be computed by topological recursion.

\medskip

Applying Theorem~\ref{Twistf} to $f = t \cdot f_\epsilon$ is the Heaviside function shifted by some sufficiently small $\epsilon$ and multiplied by a variable $t$, we establish in Theorem \ref{sysL} in Section \ref{sysf} that
\[
\Omega^{{\rm M}}_\Sigma[f_\epsilon](\sigma) = (1+t)^{\#S_{\Sigma,\sigma}^{\circ,\epsilon}},
\]
where $S_{\Sigma,\sigma}^{\circ,\epsilon}$ is the set of simple closed curves on $\Sigma$ such that $\ell_\sigma(\gamma) \leq \epsilon$. In particular, if one makes the evaluation $t=-1$, one gets the indicator function of the $\epsilon$-thick part of the moduli space
$$ \tilde V_{g,n}^\epsilon = \{ [\sigma] \in {\mathcal M}_{g,n} \mid \text{\rm{ {\small all simple interior closed geodesics on} }} [\sigma]  \text{\rm{ {\small have length at least}  }} \epsilon \}.$$
In Section \ref{ASFT} we show that the sets
$$ V_{g,n}^\epsilon = \tilde V_{g,n}^\epsilon \cap {\mathcal M}_{g,n}(\epsilon, \ldots, \epsilon)$$
solves the quantum master equation (see Theorem \ref{sysQME}) as a consequence of $\Omega^{{\rm M}}_\Sigma[f_\epsilon]_{t=-1}$ satisfying GR and thus gives a construction of a topological vertex in string field theory, as introduced and defined by  Zwiebach in \cite{BZ}. See also the remarks right after Theorem \ref{sysQME} in Section \ref{ASFT}. \\

In \cite{WKarticle}, the space of functions on the combinatorial Teichm\"uller space is equipped with the structure of a target theory, and all results of Section~\ref{SfunSctin}--\ref{S7n7n7} have an analog in this context, including a combinatorial analog Mirzakhani's identity (GR for the constant function $1$), an analog of Theorem~\ref{Twistf} (GR for statistics of combinatorial lengths of multicurves) and a result that GR in this target theory implies topological recursion after integration against the Kontsevich volume form. These results were used to give new, uniform, and geometric proofs of Witten's conjecture/Kontsevich theorem, Norbury's enumeration of lattice point in the combinatorial moduli space, and to give yet another approach to Masur--Veech volumes of quadratic differentials which is independent of hyperbolic geometry. Besides, GR on the combinatorial Teichm\"uller space was related to GR on Teichm\"uller space by means of the Bowditch--Epstein flow. \\

We envision many future applications of Geometric Recursion, both in terms of showing that various well-known constructions satisfy GR, which could allow for explicit calculations not possible otherwise, and also constructions which are only possible via GR and therefore offering new construction possibilities. We will illustrate this with an array of examples of functorial assignment $\Omega_{\Sigma} \in  E(\Sigma)^{\Gamma_\Sigma}$ for every object $\Sigma$ of $\Surf_1$ for various pre-target theories $E$. 

\vskip.5cm 

\begin{center}
\begin{tabular}{| c | c |c| } 
\hline 
No.& Pre-target Theory & Functorial assignment   \\ [0.5ex] \hline
 {\bf 1}& $ E(\Sigma) = \mathcal{C}^0({\mathcal T}_\Sigma)$ & $\Omega_\Sigma = 1$  \\[0.5ex] \hline
{\bf 2}& $E(\Sigma) = \mathcal{C}^0({\mathcal T}_\Sigma)$ & $ \Omega_\Sigma(\sigma) = \sum_{\gamma\in S_\Sigma}  \prod_{c \in \pi_0(\gamma)} f(l_\sigma(\gamma_c))$  \\[0.5ex]  \hline
{\bf 3}& $ E(\Sigma) = \mathcal{C}^0({\mathcal T}_\Sigma)$ & $\Omega_\Sigma(\sigma) = \text{Tr}(f(-\Delta_\sigma)) $  \\[0.5ex]  \hline
{\bf 4}&  $E(\Sigma) = \Omega^2({\mathcal T}_\Sigma)$ & $\Omega_\Sigma = \omega_{\text{WP}}$ \\[0.5ex]  \hline
{\bf 5}&  $E(\Sigma) = C^\infty({\mathcal T}_\Sigma, \text{End}(T\mathcal{T}_\Sigma))$ & $\Omega_\Sigma = I_{\text{Bers}}$  \\[0.5ex]  \hline
{\bf 6}&  $E(\Sigma) = \Omega^*({\mathcal T}_\Sigma)$ & $\Omega_\Sigma \in \Omega^*({\mathcal M}_\Sigma)$, $\dd\Omega_\Sigma = 0$  \\[0.5ex]  \hline
{\bf 7}&  $E(\Sigma) = C^\infty({M}_G(\Sigma), \Lambda^2 TM_G(\Sigma))$ & $\Omega_\Sigma = P_{\text{FR}} \in E(\Sigma)^{\Gamma_\Sigma}$  \\[0.5ex]   \hline
{\bf 8} & $E(\Sigma) = C^\infty(\mathcal{T}_\Sigma, C^\infty({M}_G(\Sigma,c), \text{End}(TM_G(\Sigma,c)))$ & $\Omega_\Sigma = I_{NS}\in E(\Sigma)^{\Gamma_\Sigma}$  \\[0.5ex] \hline
{\bf  9}& $E(\Sigma) = C^\infty(\mathcal{T}_\Sigma, C^\infty({M}_G(\Sigma,c)))$ & $\Omega_\Sigma = F_{\text{Ricci}} \in E(\Sigma)^{\Gamma_\Sigma}.$ \\[0.5ex]  \hline
{\bf  10} & $E(\Sigma) = C^\infty(\mathcal{T}_\Sigma, C^\infty({M}_G(\Sigma,c),\text{End}(TM_G(\Sigma,c)))^{\times 3}$ & $\Omega_\Sigma = (I,J,K)_{\text{Hitchin}} \in E(\Sigma)^{\Gamma_\Sigma}.$\\[0.5ex]  \hline
 {\bf 11} & $E(\Sigma) = \Omega^1(\mathcal{T}, \mathcal{T}\times \text{End}(V))$ & $ \Omega_\Sigma = u_{\rho_\Sigma} \in E(\Sigma)^{\Gamma_\Sigma}$ \\[0.5ex]   \hline
 {\bf 12 }& $ E(\Sigma) = {\mathcal M}(Q{\mathcal T}_{\Sigma})$ & $ \mu_{{\rm MV}}\in {\mathcal M}(\mathcal{Q}{\mathcal T}_{\Sigma})$ \\[0.5ex]  \hline
 {\bf 13} & $ E(\Sigma) = \mathbb{C}[\text{\small Heegaard diagrams } (\alpha,\beta)\, \text{\small on }\Sigma]^*$ & $ \Omega_\Sigma =  I_3(\alpha,\beta) = I_3(X^3_{(\alpha,\beta)})$ \\[0.5ex]  \hline
{\bf  14} & $ E(\Sigma) = \mathbb{C}[\text{\small Tri-section diagrams } (\alpha,\beta,\gamma)\, \text{\small on }\Sigma]^*$ & $  \Omega_\Sigma =  I_4(\alpha,\beta,\gamma) = I_4(X^4_{(\alpha,\beta,\gamma)})$\\[0.5ex]  \hline
 {\bf 15} & $ E(\Sigma) = \Omega^*(\mathcal{T}_\Sigma)$ & $\Omega_\Sigma = \varphi_{GW} \in E(\Sigma)^{\Gamma_\Sigma}$\\[0.5ex]  \hline
 {\bf 16 }& $ E(\Sigma) = \Omega^{\rm{top}}({\mathcal T}_\Sigma)$ & $\Omega_\Sigma = A_\Sigma \in E(\Sigma)^{\Gamma_\Sigma}$ \\[0.5ex] \hline
 {\bf 17} & $ E(\Sigma) = \mathscr{M}(\left\{\right.$ collection of self-avoiding loops on $\Sigma$ $\left.\right\})$. & $\mu_{{\rm CLE}_{\kappa}} \in E(\Sigma)^{\Gamma_{\Sigma}}$  \\[0.5ex] \hline
\end{tabular}
\end{center}
\vskip.5cm
Examples \textbf{1} and \textbf{2} are treated in detail in the present work. We actually conjecture that all these examples satisfies GR in a certain form. Let us now comment on each of them. 


\begin{description}
\item[3] We denote $\Delta_\sigma$ the Laplacian operator with Dirichlet boundary conditions on the bordered surface $\Sigma$ equipped with hyperbolic metric $\sigma \in \mathcal{T}_\Sigma$. The test function $f : \mathbb{R} \ra \mathbb{C}$ is assumed to decay sufficiently fast, so that $f(\Delta_\sigma)$ is trace-class. We expect that an extension of GR to surfaces with corners on their boundary, which we are in the process of writing up, will treat this example. The reason for this expectation is that, by the Selberg trace formula, $f(\Delta_\Sigma)$ is expressed as a sum over all closed simple geodesics, thus involving not just simple curves, but also self-intersecting curves; cutting along such curves yields surfaces with corners.
\item[4] The Weil--Petersson symplectic form $\omega_{\text{WP}}$ on Teichm\"{u}ller space of bordered surfaces with fixed boundary lengths. We expect that example {\bf 1} should be applied in combination with specific $2$-forms on a certain torus bundle over the Teichm\"{u}ller space of a pair of pants and on a circle bundle over the Teichm\"{u}ller space of a one-holed torus to build the right initial data and one will then get forms over torus bundles over Teichm\"{u}ller spaces, which caries the natural analog of the Weil--Petersson symplectic form for bordered surfaces.
\item[5] Bers complex structure $I_{\text{Bers}}$ on  Teichm\"{u}ller space. Again, we expect that one should need to consider similar bundles over Teichm\"{u}ller spaces as in example {\bf 4}.  
\item[6] Closed forms on Teichm\"{u}ller space representing non-trivial cohomology classes on the moduli space of curves ${\mathcal M}(\Sigma)= \mathcal{T}_\Sigma/\Gamma_\Sigma^{\partial}$ are for most classes not known explicitly. We expect that such formulae can be given via GR applied to the same torus bundles over Teichm\"{u}ller spaces as in example {\bf 4}.  
\item[7]  The Fock--Rosly Poisson structure \cite{FR}, $P_{\text{FR}}$ on moduli spaces of flat connections $M_G(\Sigma)$, where $G$ is any semi-simple Lie group either complex or real. We expect that the constructions from example {\bf 1} can be combined with Fock--Rosly Poisson structures for a pair of pants and a one holed torus to give the needed initial data for this example. More precisely, one should use here the idea of fibering (Section~\ref{sfiber}) and work with the target theory being the space of functions on the universal moduli space of flat connections $M_G(\Sigma) \times \mathcal{T}_{\Sigma}$. 
\item[8] The Narasimhan--Seshadri and Mehta--Seshadri complex structure \cite{NS,MS} on moduli spaces of flat connections $M_G(\Sigma,c)$. Here $G$ is any real semi-simple Lie group and $c$ is an assignment of conjugacy classes to each boundary components of $\Sigma$, in which we assume the holonomy around each boundary component is contained. As for \textbf{7.}  (and for all further examples involving moduli space of flat connections) the idea of fibering over Teichm\"uller space should be used.
\item[9] Ricci potentials on the moduli spaces of flat connections $M_G(\Sigma,c)$. By the work of Takhtajan and Zograf \cite{TZ}, it is given by a certain expression involving sums over simple geodesics and thus we expect it will lead to a situation similar to example {\bf 3}, captured by the extension of GR to surfaces with corners on their boundary. See also \cite{AP} where a similar formula for the universal Ricci potential appears.
\item[10] Hitchin's hyperk\"{a}hler structure \cite{Hitchin, Nit, Ko} on moduli spaces of parabolic Higgs bundles. Here  $G$ is a complex semi-simple Lie group and $c$ is as before. 
\item[11] Representations of mapping class groups $\rho_\Sigma :\Gamma_\Sigma \ra \text{Aut}(V)$. Here $V$ is a finite dimensional vector space and $\rho$ is equivalent information to a $\Gamma_\Sigma$-invariant flat connection $\nabla^\rho$ in the trivial $V$-bundle over ${\mathcal T}_\Sigma$, thus given by a $u_{\rho_\Sigma} \in \Omega^1(\mathcal{T}, \mathcal{T}\times \text{End}(V))^{\Gamma_\Sigma}$. 
\item[12] The Masur--Veech measure on the bundle over the Teichm\"uller space of meromorphic quadratic differentials. 
\item[13] In this example we let $E(\Sigma)$ be the algebraic dual of the free vector space generated by the set of Heegaard diagrams on $\Sigma$. Any invariant $I_3$ of closed oriented $3$-manifolds will give an element  in $E(\Sigma)^{\Gamma_\Sigma}$ by the assignment indicated in the table above \cite{Hempel}. The interesting question is which invariants of closed 3-manifolds satisfies GR.
\item[14] In this example we let $E(\Sigma)$ be the algebraic dual of the free vector space generated by the set of tri-section diagrams on $\Sigma$. Any invariant $I_4$ of smooth closed oriented $4$-manifolds will give an element  in $E(\Sigma)^{\Gamma_\Sigma}$ by the assignment indicated in the table above \cite{GK}. The interesting question is which invariants of closed smooth 4-manifolds satisfies GR.
\item[15] Closed forms representing cohomology classes of relevance for the Gromov--Witten invariants of a symplectic manifold, which a possible hyperbolic geometric approach to this theory might establish a GR construction of.
\item[16]  Amplitudes in closed string theory. These are briefly discussed in Section \ref{ASFT} where indication is given as to why they might satisfy GR.
\item[17]  This is the space of functions on $\mathcal{T}_{\Sigma}$ valued in Radon measures on collections of self-avoiding loops on $\Sigma$. The conformal loop ensemble ${\rm CLE}_{\kappa}$ on bordered surfaces, depending on a parameter $\kappa$, are specified by probability measures $\mu_{{\rm CLE}_{\kappa}}$ and take their origin in the work of \cite{SW10}. The Markovian property of these ensembles leads to a factorisation property conditionally to having a given multicurve part of the system of loops. By a mechanism similar to \textbf{2} (see the proof of Theorem~\ref{Genegin}) the GR sum should allow removing this conditioning, as for every multicurve one can find a pair of pants whose interior is not crossed by this multicurve.
\end{description}

In future applications, one may have to introduce variants of the formalism of this article, e.g. excising all kinds of embedded surfaces instead of only pairs of pants, or only those having specific properties with respect to an extra structure, assigning different initial data to different mapping class group orbits $\mathcal{O}$ instead of assigning $C$ to all orbits inside $\mathcal{P}_{\Sigma}^{\emptyset}$, or adapting the theory to a category of surfaces with boundaries and corners. Yet, the core idea of GR remains to consider non-local glueing in the form of (countable, absolutely convergent) sums over homotopy classes of excisions of subsurfaces assembling into a recursion for the functorial assignment $\Omega_{\Sigma}$. This ``non-local glueing'' is quite different from exact factorisation rules met in topological quantum field theories \cite{At2}, or in the Burghalea--Friedlander--Kappeler formula for the spectral determinant \cite{Friedlander}, and in fact allows more flexibility. For GR valued in functions on Teichm\"uller space, we take in Section~\ref{Boundbeh} the first step towards understanding the relation between exact glueing and non-local glueing, namely we establish in Proposition~\ref{pinchth} that the GR amplitudes factorise asymptotically in the limit where a maximal number of curves are pinched, provided the initial data decay fast enough. In other words, the GR sum reduces to a single term in this asymptotic regime.

Example \textbf{11} is closely related to conformal blocks in conformal field theories, which we also expect to satisfy GR. A hint that it may be the case is that conformal blocks satisfy bootstrap, that is, they can be expressed as a series in the neighborhood of the boundary of the moduli space obtained by glueing conformal blocks. In particular conformal blocks satisfy an asymptotic factorisation when some cycles are pinched. The fact that the conformal blocks form a representation of the mapping class group is not clear from the conformal bootstrap. We believe that GR is a formalism which could reconciliate these two aspects. As we have just mentioned, under some conditions the GR sums reduce asymptotically to one term (asymptotic factorisation) while taking into account all terms in GR sums exactly guarantees the resulting object is mapping class group invariant.

\subsection{Relation with the topological recursion}

Let $\Gamma^{{\rm 1d}}_{\Sigma}$ be the subgroup of permutations of $\pi_0(\partial \Sigma)$ preserving the map $\pi_0(\partial \Sigma) \rightarrow \{-,+\} \times \pi_0(\Sigma) \times \mathbb{N}^2$ which encodes for each boundary component, its orientation versus the orientation of the surface, the connected component of the surface to which it belongs, and the genus $g$ and number of boundary components $n$ of this connected component. For connected surfaces $\Gamma^{{\rm 1d}}_{\Sigma}$ is just the permutation group of $\pi_0(\partial_{+}\Sigma)$.  We have a natural projection $\Gamma_{\Sigma} \rightarrow \Gamma^{{\rm 1d}}_{\Sigma}$ which only records the effect of a mapping class on the boundary components. 

The geometric recursion is only relevant for target theories in which the action of the mapping class group is non-trivial. If indeed the action of $\Gamma_{\Sigma}$ on $E(\Sigma)$ factors through $\Gamma_{\Sigma}^{{\rm 1d}}$, admissible initial data do not exist since infinitely many terms in the series \eqref{firstseries} are equal. The meaningful construction in that case is simpler and summarised in the following result. As it only requires finite sums, there is no need to address questions of topology on the spaces $E(\Sigma)$ and we can work with a weaker axiomatic setting which we call pre-target theory.

\begin{theorem}
\label{thad}Assume $F:\,\Surf_{1} \rightarrow \mathcal{C}$ is a pre-target theory (Definition~\ref{PRETARGETDEF}) with glueing maps $\Theta$ and disjoint union maps $\sqcup$, such that the action of $\Gamma$ factors through $\Gamma^{{\rm 1d}}$. Assume we are given initial data
\[
A_{P} \in F(P)^{\Gamma^{{\rm 1d}}_{P}},\qquad B_{P}^{b} \in F(P),\qquad C_{P} \in F(P)^{\Gamma^{{\rm 1d}}_{P}},\qquad D_{T} \in F(T).
\]
Then there exists a unique functorial assignment
\[
\Sigma \mapsto W_{\Sigma} \in F(\Sigma)
\]
specified by
\[
W_{\emptyset} = 1_{F(\emptyset)},\qquad W_{P} = A_{P},\qquad W_{T} = D_{T},\qquad W_{\Sigma_1 \cup \Sigma_2} = W_{\Sigma_1} \sqcup W_{\Sigma_2},
\] 
and for any connected $\Sigma$ of Euler characteristic
\[
W_{\Sigma} = \sum_{\substack{b \in \pi_0(\partial_+\Sigma) \\ \{P\} = \mathcal{P}^{b}(\Sigma)/\Gamma_{\Sigma}}} \Theta_{P}(B_{P}^{b},W_{\Sigma - P}) + \tfrac{1}{2} \sum_{P \in \mathcal{P}^{\emptyset}(\Sigma)/\Gamma_{\Sigma}} \Theta_{P}(C_{P},W_{\Sigma - P}).
\] 
\end{theorem} 
We are summing here over the set of embedded pairs of pants $P \mapsto \Sigma$ up to diffeomorphisms. This set is finite and in bijection with the terms in \eqref{TRintross}: its elements are characterised by the topology of the surface $\Sigma - P$. We denote by $F$-valued \emph{topological recursion} the construction of Theorem~\ref{thad}, also abbreviated TR. We have seen in Theorem~\ref{ntinghg} how the geometric recursion valued in $E(\Sigma) = {\rm Mes}(\mathcal{T}_{\Sigma})$ is intertwined with the topological recursion valued in $F(\Sigma) = {\rm Mes}(\mathbb{R}_{+}^{\pi_0(\partial \Sigma)})$, via the integration map. This is an example of a more general phenomenon.

\begin{theorem} 
Let $E$ be a target theory, and $F$ be a pre-target theory such that the action of $\Gamma$ factors through $\Gamma^{{\rm 1d}}$. Assume that $\eta:\,E \Rightarrow F$ is a (perhaps, only partially defined) natural transformation of pre-target theories. The following commutative diagram holds in the domains where all the arrows make sense
\begin{figure}[ht!]
\begin{center}
\begin{tikzpicture}[node distance = 2cm,auto]
\node (Eini) {\small{\begin{tabular}{c} $E$-valued, admissible \\  initial data \end{tabular}}} ;
\node (OM) [right= 2cm of Eini] {$\,\,\Omega\,\,$} ;
\node (W) [below of =OM] {$\,\,W\,\,$} ;
\node (Fini) [below of =Eini] {\small{\begin{tabular}{c} $F$-valued  \\ initial data \end{tabular}}} ;
\draw[->] (Eini) to node {GR} (OM) ;
\draw[->] (Fini) to node {TR} (W) ;
\draw[->] (Eini) to node {$\eta$} (Fini) ;
\draw[->] (OM) to node {$\eta$} (W) ;
\end{tikzpicture}
\end{center}
\end{figure}
\end{theorem}
In practice, we must allow $\eta$ to be only partially defined: for instance in Theorem~\ref{ntinghg}, we can only integrate functions which are already mapping class group invariant and integrable.

We briefly explain in Section~\ref{EOfmum} that the topological recursion of Chekhov, Eynard and Orantin \cite{EORev} based on spectral curves $S$ is equivalent to the construction of Theorem~\ref{thad} for a pre-target theory $F(\Sigma)$ which is a space of meromorphic multidifferentials on $S^{\pi_0(\partial \Sigma)}$. The solution to a wealth of problems in enumerative geometry is given by the TR amplitudes associated with well-chosen spectral curves. They can, for instance, give the correlation functions of semi-simple cohomological field theories, including as particular cases the Gromov--Witten invariants of toric Calabi--Yau $3$-folds, the stationary Gromov--Witten invariants of $\mathbb{P}^1$, and the Weil-Petersson volumes of the moduli space of bordered Riemann surfaces. We found that the TR amplitudes coming from spectral curves can always be realised as the result of integration of GR amplitudes over the moduli space of bordered Riemann surfaces, followed by a Laplace transform with respect to the boundary lengths; this will be discussed elsewhere. 

Compared to the topological recursion of \cite{EORev}, the construction of Theorem~\ref{thad} allows pre-target theories where $F(\Sigma)$ has a richer dependence on the topology of the surface, \textit{i.e.} does not only depend on the number of boundary components. For instance, we believe that it should be possible to construct pre-target theories from modular operads \cite{GetzlerKapranov}.

\medskip

\vspace{0.2cm}

\subsection*{Acknowledgments}

We thank Caltech and Berkeley Mathematics Departments, the Department of Mathematics and Statistics at University of Melbourne and Paul Norbury, the organisers of the thematic month ``Topological recursion and modularity'' in the Matrix@Creswick, the Sandbjerg Estate and the QGM at Aarhus University, the EPFL and the Bernoulli Center (in particular Cl\'ement Hongler and Nicolas Monod), the Max Planck Institute f\"ur Mathematik in Bonn, hospitality at various stages of this project in excellent conditions. We thanks Curtis McMullen for many enlightening comments about Teichm\"{u}ller spaces of bordered surfaces, Hugo Parlier for the reference \cite{Parlierreduce} and for communicating to us the proof of Theorem~\ref{ThPar}, Owen Gwilliam and Peter Teichner for categorical discussions, Stavros Garoufalidis and Alessandro Giacchetto for comments, Barton Zwiebach for discussions regarding string field theory and the town of \"Olgii for its hospitality in the infancy of this project. This work was supported in part by the Danish National Research Foundation grant DNRF95, ``Centre for Quantum Geometry of Moduli Spaces'' and by the ERC-SyG project, Recursive and Exact New Quantum Theory (ReNewQuantum) which received funding from the European Research Council (ERC) under the European Union's Horizon 2020 research and innovation programme under grant agreement No 810573. G.B. benefited from the support of the Max-Planck-Gesellschaft.

\vfill

\section{Two-dimensional geometry background}
\label{S1}
\subsection{Categories of surfaces}
\label{S11}
\begin{definition}
A bordered surface $\Sigma$ is a smooth, oriented, compact $2$-dimensional manifold with an orientation of its boundary. If $\Sigma$ is non-empty, we require that each connected component has a non-empty boundary.
\end{definition}

We remark that a bordered surface has a decomposition of its boundary $\partial \Sigma$ into two subsets $\partial_{+}\Sigma$ (resp. $\partial_{-}\Sigma$) consisting of those boundary components whose orientation agree (resp. disagree) with the orientation induced by $\Sigma$. The interior of $\Sigma$ is denoted $\Sigma^{\circ}$. We use the following notation for the components of $\Sigma$
\[
\Sigma = \bigcup_{a\in \pi_0(\Sigma)} \Sigma(a).
\]
The Euler characteristic is denoted $\chi_{\Sigma}$. We say that $\Sigma$ is of type $(g,n)$ if it is connected, has genus $g$ and $n$ boundary components.

\begin{definition}
A bordered surface $\Sigma$ is stable if $\Sigma = \emptyset$ or for any $a \in \pi_0(\Sigma)$ we have $\chi_{\Sigma(a)} < 0$. It is unstable otherwise. 
\end{definition}

\begin{definition}
Let $\Surf$ be the category whose objects are stable bordered surfaces, and morphisms are isotopy classes of orientation preserving diffeomorphisms, which also preserve the prescribed orientation of the boundary.
\end{definition}

The mapping class group $\Gamma_{\Sigma}$ is the automorphism group (in the category $\Surf$) of a stable bordered surface $\Sigma$. We denote by $\Gamma_{\Sigma}^{\partial}$ the pure mapping class group, \textit{i.e.} the subgroup consisting of the mapping classes inducing the identity permutation of $\pi_0(\partial \Sigma)$.

We shall often consider surfaces with a distinguished boundary component in each connected component. This can be achieved by requiring that the orientation of the distinguished boundary components disagrees with the orientation of the surface, while the orientation of the other boundary components agrees.

\begin{definition}
Let $\Surf_{1}$ be the full subcategory of $\Surf$ consisting of bordered surfaces for which the natural map $\pi_0(\partial_{-}\Sigma) \rightarrow \pi_0(\Sigma)$ is a bijection.
\end{definition}
If $\Sigma$ is a connected object in $\Surf_{1}$, we often denote by $b_1$ the boundary component with disagreeing orientation, \textit{i.e.} $\partial_{-}\Sigma = b_1$. We remark that for a surface in $\Surf_{1}$, the mapping classes in $\Gamma_{\Sigma}$ must leave $b_1$ stable.

\subsection{Multicurves and cutting}

As the notions of homotopy, isotopy and ambient isotopy are the same for $1$-dimensional closed submanifolds of bordered surfaces \cite{Epstein}, we will not distinguish between them and just refer to them jointly as homotopy.

\begin{definition}
A multicurve is the homotopy class of a (possibly empty) one-dimensional closed submanifold $c$ of a bordered surface $\Sigma$, with no component null-homotopic or homotopic to a boundary component of $\Sigma$. A multicurve is primitive when we add the requirement that two components of $c$ cannot be homotopic to each other. A simple closed curve is a primitive multicurve with a single component.

We denote $M_{\Sigma}$ the set of multicurves, $M_{\Sigma}'$ the set of primitive multicurves, $S_{\Sigma}^{\circ}$ the set of simple closed curves, and $S_{\Sigma} = S_{\Sigma}^{\circ} \cup \pi_0(\partial \Sigma)$.
\end{definition}

We have the operation of cutting a stable bordered surface $\Sigma$ along a primitive multicurve $\gamma$, producing a new stable bordered surface $\Sigma^{\gamma}$. If $c$ is a representative of the ambient isotopy class $\gamma$, the new boundary components in the bordered surface $\Sigma^c$ (obtained from $\Sigma$ by doubling the components of $c$) receive the orientation induced from the orientation of the surface. If $c'$ is another representative of $\gamma$, then an ambient isotopy between $c$ and $c'$ determines a diffeomorphism between $\Sigma^c$ and $\Sigma^{c'}$, whose isotopy class is independent of the choice of the ambient isotopy. To be precise, $\Sigma^{\gamma}$ will stand for any object in $\Surf$ of the form $\Sigma^c$ and we note that any two such objects are canonically isomorphic in $\Surf$.

\subsection{Excising a pair of pants}
\label{glueinpant}
The cutting operation which is particularly relevant for us is the excision of a pair of pants (\textit{i.e.} a genus $0$ surface with $3$ boundary components) from a connected surface $\Sigma$ in $\Surf_{1}$ of type $(g,n)$. For this purpose we assume $2g + n \geq 0$.  

Let us fix a connected bordered surface $P$ of genus $0$ with $3$ \emph{ordered} boundary components, \textit{i.e.} a fixed bijection between $\pi_0(\partial P)$ and $\{1,2,3\}$. We turn it into an object of $\Surf_{1}$ by declaring that the orientation on the first boundary component disagrees with the one of $P$ while the orientation of the second and the third boundary component agrees with the one of $P$. Equivalently, $P$ is a pair of pants in $\Surf_{1}$ together with an ordering of $\pi_0(\partial_{+}\Sigma)$.

Let $f:\,P \longrightarrow \Sigma$ be an embedding of a pair of pants with ordered boundary components such that $f(\partial_{-}P) = b_1 = \partial_{-}\Sigma$ and for any $\beta \in \pi_0(\partial_{+}P)$ either $f(\beta)$ coincides with a boundary component of $\Sigma$ or $f(\beta)$ is neither null homotopic relative to $\partial \Sigma$ nor boundary parallel. Let us cut $\Sigma$ along the boundary components of $f(P)$ which are not boundaries of $\Sigma$, \textit{i.e.} along the multicurve $\gamma_P = f(\partial P) \cap \Sigma^{\circ}$. This multicurve has either one or two components. In the latter case we denote $\gamma_P^1$ and $\gamma_P^2$ the two components respecting the order in which they appear in $\pi_0(\partial_{+}P)$. There are finitely many possible topologies for
\[
\Sigma^{\gamma} = P \,\cup\, \overline{\Sigma \setminus f(P)}.
\]
\begin{itemize}
\item[\textbf{I}] $(g,n) \neq (1,1)$, $g \geq 1$, $\gamma_P$ has two components and $\overline{\Sigma \setminus f(P)}$ is connected. Then $\overline{\Sigma \setminus f(P)}$ must have type $(g - 1,n + 1)$.
\item[\textbf{I'}] $(g,n) \neq (0,3)$, $f(P)$ has two  boundary components in common with $\Sigma$, namely $b_1$ and some $b \subseteq \partial_{+}\Sigma$. In this case we require that $f$ is such that $b$  is the second boundary component of $f(P)$. Then $\overline{\Sigma \setminus f(P)}$ has type $(g,n - 1)$.
\item[\textbf{II}] $\gamma_P$ has two components and $\overline{\Sigma \setminus f(P)}$ is not connected. Then $\overline{\Sigma \setminus f(P)}$ must have two ordered connected components. The $i$th one has genus $g_i$, contains $\gamma_P^i$ and a subset $N_i \subseteq \partial_{+}\Sigma$ of other boundary components. We must have $g_1 + g_2 = g$ and $N_1 \cup N_2 = \pi_0(\partial_{+}\Sigma)$ such that
\[
\forall i \in \{1,2\},\qquad (g_i,|N_i|) \neq (0,1),(0,2).
\]
\end{itemize}

We change the choice of orientation of the boundary components in $\overline{\Sigma \setminus f(P)}$ in the following way. In case \textbf{I}, we reverse the orientation of $\gamma_P^1$ so that it disagrees with the orientation of the surface, while the orientation of $\gamma_P^2$ remains in agreement with the orientation of the surface. In case \textbf{I'}, we reverse the orientation of $\gamma_P$ (which is a simple closed curve). In case \textbf{II}, we reverse the orientation of $\gamma_P^1$ and of $\gamma_P^2$. In this way, $\overline{\Sigma \setminus f(P)}$ becomes a surface in $\Surf_{1}$.

By an argument similar to our description of cutting, if $f':\,P \longrightarrow \Sigma$ belongs to the same homotopy class of embedding of pairs of pants as $f$, we have a canonical isomorphism in $\Surf_{1}$ between $\overline{\Sigma \setminus f(P)}$ and $\overline{\Sigma \setminus f'(P)}$.

\begin{definition}
We let $\mathcal{P}_{\Sigma}$ be the set of homotopy classes, relative to the boundary, of embedded pairs of pants $f:\,P \longrightarrow \Sigma$ as described above. If $d \in \{\emptyset\} \cup \pi_0(\partial_{+}\Sigma)$, $\mathcal{P}_{\Sigma}^{d}$ is the subset of $[f:\,P \rightarrow \Sigma] \in \mathcal{P}_{\Sigma}$ such that $f(\partial_{+}P) \cap \partial_{+}\Sigma = d$.
\end{definition}
For conciseness, we use the notation $[P]$ for elements of $\mathcal{P}_{\Sigma}$. We also write $\Sigma - P$ for any object in $\Surf_{1}$ which is obtained by excision of a representative of the homotopy class $[P]$ from $\Sigma$.

The mapping class group $\Gamma_{\Sigma}$ acts on $\mathcal{P}_{\Sigma}$, with a finite number of orbits partitioned according to the cases \textbf{I}, \textbf{I'} and \textbf{II}. 
In particular, for each $b \in \pi_0(\partial_{+}\Sigma)$ the set $\mathcal{P}^{b}_{\Sigma}$ is itself a $\Gamma_{\Sigma}^{\partial}$-orbit.

\vspace{0.2cm}

\noindent \textsc{Alternative description of $\mathcal{P}_{\Sigma}$}

\vspace{0.2cm}

\begin{figure}[ht!]
\begin{center}
\includegraphics[width=0.45\textwidth]{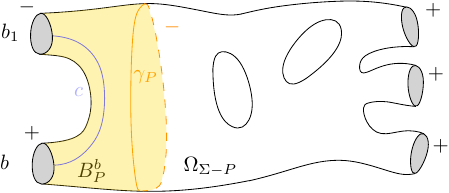}
\end{center}
\caption{\label{GRBFIG} The $B$ orbits (case \textbf{I'}).}
\end{figure}

If $b$ is a component in $\partial_{+}\Sigma$, the set $\mathcal{P}_{\Sigma}^{b}$ is in bijection with the set of homotopy classes of simple curves $c$ in $\Sigma$ from a point in $b_1$ to a point in $b$. Given such a $c$, we can indeed consider the free homotopy class $[\gamma]$ of a  simple closed curve $\gamma$ obtained by following $c$, going around $b$, following $c$ in the opposite direction and going around $b_1$ (for the two boundary components we here used the orientation induced by the surface). We can always find a representative $\gamma$ of this class which is simple. The component of $\Sigma^{\gamma}$ which contains $b_1$ is a pair of pants, with its ordered boundary components $(b_1,b,\gamma)$. We therefore obtain an embedding $f:\,P \longrightarrow \Sigma$ such that $f(P)$ contains $c$ as shown in Fig.~\ref{GRBFIG}, which is of type  \textbf{I'}. The homotopy class of this embedding is independent of the choice of representatives of $[c]$, $[\gamma]$ and $P$. The reciprocal bijection is obtained by associating to an embedding $f:\,P \longrightarrow \Sigma$ as in case \textbf{I'}  the homotopy class relative to the boundary of the curve $c$ shown in Figure~\ref{GRBFIG}.

\begin{figure}[ht!]
\begin{center}
\includegraphics[width=0.5\textwidth]{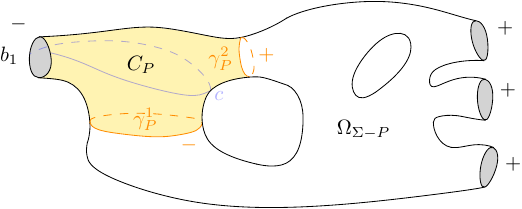}
\end{center}
\caption{\label{GRC1FIG} The $C$ orbits (case \textbf{I}).}
\end{figure}

\begin{figure}[ht!]
\begin{center}
\includegraphics[width=0.5\textwidth]{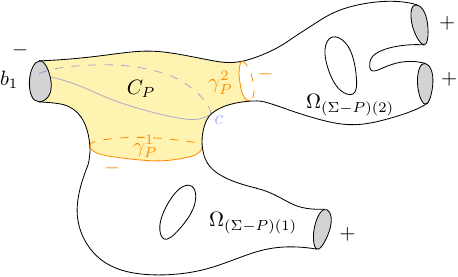}
\end{center}
\caption{\label{GRC2FIG} The $C$ orbits (case \textbf{II}).}
\end{figure}

Likewise, $\mathcal{P}_{\Sigma}^{\emptyset}$ is in bijection with the set of homotopy class of simple closed oriented curves $c$ in $\Sigma$ with one point in $b_1$, which are not homotopic to any composition of the form $\tilde{c}^{-1} \cdot b \cdot \tilde{c}$ for some component $b$ in  $\partial_{+}\Sigma$ and some $\tilde{c}$ in the equivalent description we just gave of $\mathcal{P}_{\Sigma}^{b}$. Here we use $\cdot$ for the composition of curves in the path groupoid of $\Sigma$. The data of such a $[c]$ indeed determines two free homotopy class relative to the boundary of closed curves $[\gamma^1]$ and $[\gamma^2]$, where $[\gamma^1]$ is to the left of $c$. After we pick simple representatives $\gamma^1$ and $\gamma^2$, we obtain an embedding $f:\,P \longrightarrow \Sigma$ (of the type \textbf{I} or \textbf{II}) of the component of $\Sigma^{\gamma^1 \cup \gamma^2}$ which contains $b_1$. The curve $c$ appears in the embedded pair of pants as it is shown in Figures \ref{GRC1FIG} and \ref{GRC2FIG}.

\section{The geometric recursion and its main properties}
\label{S3}
Our goal is to construct, by means of sums over pair of pants decompositions, for any surface $\Sigma$ in $\Surf_{1}$ an element $\Omega_{\Sigma}$ of a ``space'' $E(\Sigma)$ attached to $\Sigma$, which is a functorial assignment. First of all we need to specify what is the space $E(\Sigma)$ to which $\Omega_{\Sigma}$ will belong. For the construction to be meaningful, it should incorporate a compatible set of union and glueing maps, and a mechanism making sense of the countable sums over pair of pants decompositions. This data will be called a \emph{target theory}. More precisely, it will be a functor $E$ from $\Surf_{1}$ to some multicategory $\mathcal{C}$, satisfying a list of axioms  specified in Section~\ref{STarget}. The definition of $\Omega_{\Sigma}$ will depend on a small amount of \emph{initial data}, specifying what happens for a pair of pants, for a torus with one boundary, and a way to inductively increase the complexity of the surfaces (Section~\ref{SIniData}). 

We denote the inductive process leading to $\Omega_{\Sigma}$ from such initial data \emph{geometric recursion} (with values in the chosen target theory $E$).

\subsection{The category of projective systems of vector spaces}
\label{STarget}

Let $\mathbb{K}$ be $\mathbb{R}$ or $\mathbb{C}$. We will choose $\mathcal{C}$ to be the category Pro\,-$\mathcal{V}$, where $\mathcal{V}$ be the category of Hausdorff, complete, locally convex topological vector spaces over $\mathbb{K}$. Recall that an object $V$ of $\mathcal{V}$ comes with a (possibly uncountable) family of seminorms $(|\cdot|_\alpha)_{\alpha \in \mathscr{A}}$ which induces the topology of the vector space $V$ such that it is Hausdorff and complete as a locally convex topological vector space. We remark that $\mathcal{V}$ admits projective limits \cite{Taylor}. Pro\,-$\mathcal{V}$ is a category defined using projective families of objects in $\mathcal{V}$ which we now explain in detail.

\begin{definition}
An object ${\mathbf V}$ of Pro\,-$\mathcal{V}$ is a directed set $\mathscr{I}$ and an inverse system over $\mathscr{I}$ of objects 
\[
(V^i, (|\cdot|_{i,\alpha})_{\alpha \in \mathscr{A}^{i}})_{i\in \mathscr{I}}
\]
 of $\mathcal{V}$, together with an object $V$ of $\mathcal{V}$ which is the projective limit of $(V^{i})_{i \in \mathscr{I}}$.
 \end{definition}
The base field $\mathbb{K}$ can be considered as an object in Pro\,-$\mathcal{V}$ with set of indices $\mathscr{I} = \{1\}$ and $\mathscr{A}^1 = \{1\}$ with single norm $|\cdot|_{1,1} = |\cdot|$. By definition, for any $i \in \mathscr{I}$ we have a restriction linear map $\rho^i:\,V \rightarrow V^i$ which is continuous. Using this map, elements of $V$ can be silently considered as element of $V^i$ for any $i \in \mathscr{I}$.

\begin{definition}
A morphism ${\mathbf \Phi}:\,\mathbf{V}_1 \rightarrow \mathbf{V}_2$ in Pro\,-$\mathcal{V}$ is an inverse system of continuous linear maps 
\[
\Phi^{i,j} : V^{i}_1 \longrightarrow V^{j}_2
\]
indexed by $(i,j) \in \mathscr{I}_{1} \times \mathscr{I}_{2}$ such that $i \geq h_{j}$ for some $h_{j} \in \mathscr{I}_{1}$. This  data determine continuous linear maps $\Phi^{j} := \Phi^{i,j} \circ \rho^i:\,V_1 \rightarrow V^j_{2}$ which do not depend on the choice of $i \geq h_{j}$ by the very definition of an inverse system.
\end{definition}

The category Pro\,-$\mathcal{V}$ can be considered as a multicategory using the cartesian product of objects of $\mathcal{V}$ and adopting the following definition.

\begin{definition}
\label{multilinearmor}Given three objects $\mathbf{V}_1,\mathbf{V}_2,\mathbf{V}_3$ of Pro\,-$\mathcal{V}$, a bilinear morphism
\[
{\mathbf \Phi} : {\mathbf V}_1\times {\mathbf V}_2 \longrightarrow {\mathbf V}_3
\]
is an inverse system of continuous bilinear maps
\[
\Phi^{i,j,k} : V^{i}_1 \times V^{j}_2 \longrightarrow V^{k}_3,
\]
indexed by $(i,j,k) \in \mathscr{I}_1\times \mathscr{I}_2 \times \mathscr{I}_{3}$ such that $(i,j) \geq h_{k}$ for some $h_{k} \in \mathscr{I}_{1} \times \mathscr{I}_{2}$ -- where we use the lexicographic order on $\mathscr{I}_{1} \times \mathscr{I}_{2}$. This data determine continuous bilinear maps $\Phi^{k} = \Phi^{i,j,k} \circ (\rho_{i} \times \rho_{j}):\,V_1 \times V_2 \ra V_3^k$ which do not depend on the choice of $(i,j) \geq h_{k}$.

Multilinear morphisms are defined in a similar fashion. 
\end{definition}

We do not consider tensor products of topological vector spaces, which usually turn bilinear maps into linear ones, since we want to (and we can) avoid discussing completions of the algebraic tensor product. This makes our setup more flexible.

This multicategory is suited to treat in a uniform way spaces of functions, forms, distributions, sections of vector bundles on many different topological spaces with various structures associated to surfaces, such as various moduli spaces of different kinds of structures on surfaces, possibly coupled to Teichm\"uller space, etc. It seems to contain the minimal structure needed to make sense of GR. Applications may justify other choices of categories. For instance, one could use the pro-category of the category of chain complexes of objects in $\mathcal{V}$.

Hereafter, we set once for all $\mathcal{C} =$ Pro\,-$\mathcal{V}$. To keep lighter notations, we refrain from the use of bold letters for objects and morphisms in $\mathcal{C}$. For instance, an object $\mathbf{V} = (V^i,(|\cdot|_{i,\alpha})_{\alpha \in \mathscr{A}^i})_{i \in \mathscr{I}})$ will simply be denoted $V$, and the context will leave no room for confusion with the letter $V$ also used for the projective limit of $(V^i)_{i \in \mathscr{I}}$. 

\subsection{The notion of a target theory}
\label{STargetT}

\begin{definition}
\label{PRETARGETDEF} A ($\mathcal{C}$-valued) pre-target theory is a functor $E$ from $\Surf_{1}$ to the category $\mathcal{C}$, together with functorial extra structures specified by the union and excision axioms below. 
\end{definition}

\noindent \textsc{Union --} For any two surfaces $\Sigma_1$ and $\Sigma_2$, we ask for a bilinear morphism in $\mathcal{C}$
\[
\sqcup\,:\,E(\Sigma_1) \times E(\Sigma_2) \longrightarrow E(\Sigma_1 \cup \Sigma_2)
\]
compatible with commutativity and associativity of cartesian products and of unions. We require $E(\emptyset) = \mathbb{K}$ and the union map $\sqcup\,:\,E(\emptyset) \times E(\Sigma) \rightarrow E(\Sigma)$ is specified by $1\sqcup v = v$.

\vspace{0.2cm}

\noindent \textsc{Excision --} For any connected surface $\Sigma$ with $\chi_{\Sigma} \leq -2$ and $[P] \in \mathcal{P}_{\Sigma}$, we ask for the data of a bilinear morphism in $\mathcal{C}$
\[
\Theta_{P}:\,E(P) \times E(\Sigma - P) \longrightarrow E(\Sigma).
\]

\vspace{0.2cm}

From the union axiom, one can define multilinear morphisms corresponding to any finite union of surfaces. It can be rephrased by saying that $E$ is a multifunctor between the multicategories $(\Surf_{1},\cup)$ and $(\mathcal{C},\times)$. The morphisms $\Theta_{P}$ in the excision axiom will be called the glueing morphisms.

\vspace{0.2cm}

\vspace{0.2cm}

\begin{definition}
\label{TARGETDEF} A target theory is a pre-target theory $E$ together with, for any surface $\Sigma$ in $\Surf_{1}$, a functorial collection of functions defined on the set of simple closed curves and boundary components
\[
l_{i,\alpha}:\, S_{\Sigma} \longrightarrow \mathbb{R}_{\geq 0}.
\]
Here $i \in \mathscr{I}_{\Sigma}$ and $\alpha \in \mathscr{A}^{i}_{\Sigma}$ respectively index the spaces $E^{i}(\Sigma)$ and their semi-norms $|\cdot|_{i,\alpha}$ provided by the pre-target theory. We call $l_{i,\alpha}$ the length functions and we require that they satisfy the following three properties.
\end{definition}

\vspace{0.2cm}

\noindent \textsc{Polynomial growth --} For any $i\in \mathscr{I}_{\Sigma}$, there exists functorial constants $N_i,d_i > 0$ such that for any $\alpha \in \mathscr{A}^{i}_{\Sigma}$ and $L > 0$,
\[
\#\big\{ \gamma \in M_{\Sigma}' \quad \big| \quad l_{i,\alpha}(\gamma) \leq L \big\} \leq N_i\,(1 + L)^{d_i}
\]
Here, we extended the definition of $l_{i,\alpha}$ to non primitive multicurves by setting
\[
l_{i,\alpha}(\gamma) = \sum_{\beta \in \pi_0(\gamma)} l_{i,\alpha}(\beta).
\]

\vspace{0.2cm}

\noindent \textsc{Small pair of pants --} For any $i \in \mathscr{I}_{\Sigma}$, there exists a functorial $Q_i > 0$ such that for any $\alpha \in \mathscr{A}^{i}_{\Sigma}$ and $b \in \pi_0(\partial_{+}\Sigma)$
\begin{itemize}
\item[$\bullet$] the number of $[P] \in \mathcal{P}_{\Sigma}^{b}$ such that $l_{i,\alpha}(\gamma) \leq l_{i,\alpha}(b_1) + l_{i,\alpha}(b)$ is less than $Q_{i}$ ;
\item[$\bullet$] the number of $[P] \in \mathcal{P}_{\Sigma}^{\emptyset}$ such that $l_{i,\alpha}(\gamma_1) +l_{i,\alpha}(\gamma_2) \leq l_{i,\alpha}(b_1)$ is less than $Q_{i}$;
\end{itemize}
with the notations introduced in Section~\ref{glueinpant}.

The functoriality of the constants $d_i,D_i,Q_i$ prosaically means that they only depend on the topology of $\Sigma$ -- taking into account that the ordered sets $\mathscr{I}_{\Sigma}$ to which $i$ belongs are canonically isomorphic for surfaces $\Sigma$ of the same topology.

\vspace{0.2cm} 

\noindent \textsc{Restriction --} We assume that for any $k \in \mathscr{I}_{\Sigma}$, there exists $K_{k} > 0$ functorial, such that for any $[P] \in \mathcal{P}_{\Sigma}$, there exists $j_k \in \mathscr{I}_{\Sigma - P}$ functorial, such that for any $\alpha \in \mathscr{A}^k_{\Sigma}$ and $j \geq j_k$, there exists $\alpha' \in \mathscr{A}_{\Sigma - P}^{j}$ functorial such that for any $\gamma \in S_{\Sigma - P}$, we have
\[
l_{j,\alpha'}(\gamma) \leq K_{k}\,l_{k,\alpha}(\iota_*\gamma),
\]
where $\iota:\,(\Sigma - P) \rightarrow \Sigma$ is the natural inclusion.

\vspace{0.2cm}

This axiom means that we can control the lengths of a curve in an excised surface by the lengths of the same curve in the surface before excision.

We define pseudo-norms in $E(\Sigma)$, indexed by $i \in \mathscr{I}_{\Sigma}$ and $\varsigma > 0$
\[
||v||_{i,\varsigma} = \sup_{\alpha \in \mathscr{A}_{\Sigma}^i} |v|_{i,\alpha} \prod_{\beta \in \pi_0(\partial \Sigma)} (1 + l_{i,\alpha}(\beta))^{-\varsigma},
\]
and introduce the (not necessarily closed) subspace of bounded elements
\[
{}^{\flat}E(\Sigma) := \big\{v \in E(\Sigma)\quad \big|\quad  \forall i \in \mathscr{I}_{\Sigma}\quad  \exists \varsigma_i \geq 0\quad   ||v||_{i,\varsigma_i} < +\infty\big\}.
\]

\subsection{Initial data}
\label{SIniData} 

\begin{definition}
\label{DEFDECAY0} \emph{Initial data} for a given target theory $E$ are functorial assignments
\begin{itemize} 
\item[$\bullet$] of $A_{P},C_{P} \in E(P)^{\Gamma_{P}}$ for any pair of pants $P$.
\item[$\bullet$] of $B_{P}^{b} \in E(P)$ for any pair of pants $P$ in which some $b \in \pi_0(\partial_+P)$ has been selected.  
\item[$\bullet$] of $D_{T} \in E(T)^{\Gamma_{T}}$ for any torus with one boundary $T$.
\end{itemize}
\end{definition}
It is enough to make such choices $(A_{P},B_{P}^{b},C_{P},D_{T})$ for a reference $P$ and $T$ in such a way that
\[
A_{P},C_{P} \in E(P)^{\Gamma_{P}},\qquad D_{T} \in E(T)^{\Gamma_{T}}.
\]

If $P'$ is another pair of pants and $X \in E(P)^{\Gamma_{P}}$, we can use a morphism $f:\,P \rightarrow P'$ to transport $X' = E(f)(X) \in E(P')^{\Gamma_{P'}}$. As $X$ is invariant under the mapping class group, $X'$ is independent of the choice $f$. We can define uniquely $A_{P'}$ and $C_{P'}$ for all pair of pants $P'$ in $\Surf_{1}$. Likewise one defines unambiguously for any $b' \in \pi_0(\partial_{+}\Sigma)$
\[
B_{P'}^{b'} = E(f)(B^{f^{-1}(b')}_{P})
\]
The same argument defines $D_{T'}$ for all tori $T'$ with one boundary in $\Surf_{1}$.

\begin{definition}
\label{DEFDECAY}  An initial data is called \emph{admissible} if $A_P$ and $D_{T}$ are bounded and $B_{P}^{b}$ and $C_{P}$ satisfy the decay properties stated below, where for $x \in \mathbb{R}$ we use the notation $[x]_{+} = \max(x,0)$. 
\end{definition}
 
\vspace{0.2cm}  
 
\noindent \textsc{Decay --}  For any connected surface $\Sigma$ in $\Surf_{1}$ and any $[P] \in \mathcal{P}_{\Sigma}$, we require that for any $k \in \mathscr{I}_{\Sigma}$, there exist $t_{k} \geq 0$ and $j_k \in \mathscr{I}_{\Sigma - P}$ functorial such that for any $s > 0$, there exists $M_{k,s} > 0$ functorial such that for any $\alpha \in \mathscr{A}_{\Sigma}^{k}$, any $j \geq j_{k}$ and any $\alpha' \in \mathscr{A}_{\Sigma - P}^{j}$, we have for any $v \in E(\Sigma - P)^{\Gamma_{\Sigma - P}}$
\begin{itemize} 
\item[$\bullet$] if $b \in \pi_0(\partial_{+}\Sigma)$ and $[P] \in \mathcal{P}_{\Sigma}^{b}$,
\begin{equation}
\label{decB} \big| \Theta^{k}_P(B_{P}^b,v)\big|_{k,\alpha} \leq M_{k,s}\,|v|_{j,\alpha'}\,\frac{\prod_{\beta \in \pi_0(\partial \Sigma) \cup \{\gamma_P\}} (1 + l_{k,\alpha}(\beta))^{t_k}}{\big(1 + [l_{k,\alpha}(\gamma_P) - l_{k,\alpha}(b) - l_{k,\alpha}(b_1)]_{+}\big)^{s}}.
\end{equation}
\item[$\bullet$] if $[P] \in \mathcal{P}_{\Sigma}^{\emptyset}$, 
\begin{equation}
\label{decC} \big|\Theta^{k}_P(C_{P},v)\big|_{k,\alpha} \leq M_{k,s}\, |v|_{j,\alpha'} \frac{\prod_{\beta \in \pi_0(\partial \Sigma) \cup \{\gamma_P^1,\gamma_P^2\}} (1 + l_{k,\alpha}(\beta))^{t_k}}{\big(1 + [l_{k,\alpha}(\gamma_P^1) + l_{k,\alpha}(\gamma_P^2) - l_{k,\alpha}(b_1)]_{+}\big)^{s}}.
\end{equation}
\end{itemize}

\subsection{Definition of the geometric recursion}

Let $(A,B,C,D)$ be admissible initial data for a target theory $E$.

\begin{definition} 
\label{GRdef} We define the GR amplitudes as follows.

\begin{itemize}
\item[$\bullet$] We put $\Omega_{\emptyset} := 1 \in E(\emptyset) = \mathbb{K}$.
\item[$\bullet$]  If $P$ is a pair of pants in $\Surf_{1}$, we put $\Omega_{P} = A_{P}$.
\item[$\bullet$] If $T$ is a torus with one boundary in $\Surf_{1}$, we put $\Omega_{T} = D_{T}$.
\item[$\bullet$] For disconnected surfaces $\Sigma$ in $\Surf_{1}$, we declare $\Omega_{\Sigma} := \bigsqcup_{a\in \pi_0(\Sigma)} \Omega_{\Sigma(a)}.$
\item[$\bullet$]  When $\Sigma$ is a connected surface with $\chi_{\Sigma} \leq -2$, we seek to inductively define
\begin{equation}
\label{defGR}\Omega_{\Sigma} := \sum_{b \in \pi_0(\partial_{+}\Sigma)} \sum_{[P] \in \mathcal{P}_{\Sigma}^{b}} \Theta_{P}(B_{P}^{b},\Omega_{\Sigma - P}) + \frac{1}{2} \sum_{[P] \in \mathcal{P}_{\Sigma}^{\emptyset}} \Theta_P(C_{P},\Omega_{\Sigma - P})
\end{equation}
as an element of $E(\Sigma)$.
\end{itemize}
\end{definition}

\begin{theorem}
\label{welldefGR} \label{converge} The assignment $\Sigma \mapsto \Omega_{\Sigma}$ is well-defined for any surface in $\Surf_{1}$. More precisely, the series \eqref{defGR} converges absolutely for any of the seminorms $|\cdot|_{i,\alpha}$ to a unique limit. This limit is an element of ${}^{\flat}E(\Sigma)$ which is functorial. In particular, $\Omega_{\Sigma}$ is mapping class group invariant.
\end{theorem}

\noindent \textbf{Proof.} The result holds true in the case $\chi_{\Sigma} = -1$ by the assumptions on initial data. Assume it holds for all connected surfaces of Euler characteristic strictly smaller than some $\chi_0 \leq -2$. Let $\Sigma$ be a surface in  $\Surf_{1}$ with Euler characteristic $\chi_0$. For any $[P] \in \mathcal{P}_{\Sigma}$ the induction hypothesis applies to $\Sigma - P$. In particular $\Omega_{\Sigma - P}$ is functorial and belongs to ${}^{\flat}E(\Sigma - P)$.

Since $\Theta_{P}$, $B_{P}^{b}$ and $\Omega_{\Sigma - P}$ are functorial, the value of $\Theta_{P}(B_{P}^{b},\Omega_{\Sigma - P}) \in E(\Sigma)$ is independent of the embedding $f:\,P \rightarrow \Sigma$ representing a given class $[P] \in \mathcal{P}_{\Sigma}^{b}$. Likewise, $\Theta_{P}(C_{P},\Omega_{\Sigma - P}) \in E(\Sigma)$ is independent of the embedding $f:\,P \rightarrow \Sigma$ representing a given class $[P] \in \mathcal{P}_{\Sigma}^{\emptyset}$ due to the functoriality of $C_{P}$.

Let us consider a $\Gamma^{\partial}_{\Sigma}$-orbit $\mathcal{O} \subseteq \mathcal{P}_{\Sigma}$. By functoriality of $E$ and the glueing morphisms we can find
\begin{itemize}
\item[$\bullet$] an ordered set $\mathscr{I}_{\mathcal{O},1}$, to which the ordered sets $\mathscr{I}_{P}$ for representatives of $[P] \in \mathcal{O}$ are canonically identified.
\item[$\bullet$] an ordered set $\mathscr{I}_{\mathcal{O},2}$, to which the ordered sets $\mathscr{I}_{\Sigma - P}$ for representatives of $[P] \in \mathcal{O}$ are canonically identified.
\item[$\bullet$] an object $E_{\mathcal{O},2}$ in $\mathcal{C}$, to which the objects $E(\Sigma - P)$ for representatives of $[P] \in \mathcal{O}$ are (perhaps non-canonically) isomorphic to. However, the functorial elements $\Omega_{\Sigma - P}$ are all identified with a unique mapping class group invariant element $\Omega_{\mathcal{O}} \in E_{\mathcal{O},2}$. We denote $\mathscr{A}_{\mathcal{O},2}^{j}$ the set indexing the semi-norms in $E_{\mathcal{O},2}^j$ for $j \in \mathscr{I}_{\mathcal{O},2}$.
\end{itemize}

We first examine the case of $\mathcal{O} = \mathcal{P}_{\Sigma}^{b}$. Let $k \in \mathscr{I}_{\Sigma}$ and $\alpha \in \mathscr{A}_{\Sigma}^{k}$ and use the decay axiom for $B^{b}$. It ensures the existence of $t_k$ and $j_k \in \mathscr{I}_{\mathcal{O},2}$ such that, for any $s > 0$ there exists $M_{k,s} > 0$ satisfying for any $j \geq j_k$, any $[P] \in \mathcal{P}_{\Sigma}^{b}$ and $\alpha' \in \mathscr{A}^{j}_{\mathcal{O},2}$
\begin{equation}
\label{boundThetaB} \big|\Theta_{P}^{k}(B_{P}^{b},\Omega_{\Sigma - P})\big|_{k,\alpha} \leq M_{k,s}\,|\Omega_{\mathcal{O}}|_{j,\alpha'}  \frac{(1 + l_{k,\alpha}(\gamma_P))^{t_k}}{\big(1 + [l_{k,\alpha}(\gamma_P) - l_{k,\alpha}(b) - l_{k,\alpha}(b_1)]_{+}\big)^{s}}\,\prod_{\beta \in \pi_0(\partial \Sigma)} (1 + l_{k,\alpha}(\beta))^{t_{k}}.
\end{equation}
By the induction hypothesis, $\Omega_{\mathcal{O}}$ belongs to ${}^{\flat}E_{\mathcal{O},2}$, so there exists constants $\varsigma_{\mathcal{O}} > 0$ and $W_{\mathcal{O},j} > 0$ such that
\begin{equation}
\label{thenineq}|\Omega_{\mathcal{O}}|_{j,\alpha'} \leq W_{\mathcal{O},j}\,\prod_{\beta \in \pi_0(\partial(\Sigma - P))} (1 + l_{j,\alpha'}(\beta))^{\varsigma_{\mathcal{O}}}.
\end{equation}
We now pick $j$ large enough compared to $k$ and $\alpha'$ so that the restriction axiom on lengths function can be used, and deduce there exists a constant $K'_{k} > 0$ such that
\begin{equation}
\label{Omegsdgg}|\Omega_{\mathcal{O}}|_{j,\alpha'} \leq W_{\mathcal{O}}\,K'_k\,(1 + l_{k,\alpha}(\gamma_P))^{\varsigma_{\mathcal{O}}} \prod_{\substack{\beta \in \pi_0(\partial \Sigma) \\ \beta \neq b_1,b}} (1 + l_{k,\alpha}(\beta))^{\varsigma_{\mathcal{O}}},
\end{equation}
where we simply denoted $W_{\mathcal{O}} = W_{\mathcal{O},j}$ for the chosen $j$, and wrote down separately the contribution of $\gamma_P$ and the contribution of the other boundary components of $\Sigma - P$ to the product in \eqref{thenineq}. We then insert this inequality in \eqref{boundThetaB} and would like to perform the sum over all $[P] \in \mathcal{P}_{\Sigma}^{b}$. We distinguish the contribution of small pairs of pants, i.e. those satisfying $l_{k,\alpha}(\gamma_P) \leq l_{k,\alpha}(b_1) + l_{k,\alpha}(b)$ from the non-small ones.

For small pairs of pants, we can bound
\[
\big|\Theta_{P}^{k}(B_{P}^{b},\Omega_{\Sigma - P})\big|_{k,\alpha} \leq M_{k,s}\,W_{\mathcal{O}}\,K'_k\,\prod_{\beta \in \pi_0(\partial \Sigma)} (1 + l_{k,\alpha}(\beta))^{t_{k} + \varsigma_{\mathcal{O}}},
\]
and the eponymous axiom says there are less than $Q_k$ such terms.

The set of non-small $[P] \in \mathcal{P}_{\Sigma}^{b}$ (with respect to $l_{k,\alpha}$) injects in the set of simple closed curves $\gamma$ in the interior of $\Sigma$ such that $l_{k,\alpha}(\gamma_P) \geq l_{k,\alpha}(b_1) + l_{k,\alpha}(b)$. Setting $M'_{k,s,\mathcal{O}} = M_{k,s}K'_k W_{\mathcal{O}}$, we deduce that the sum of $|\Theta_{P}^{k}(B_{P}^b,\Omega_{\Sigma - P})|_{k,\alpha}$ over non-small pairs of pants admits for upper bound
\begin{equation*}
\begin{split}
& \quad M'_{k,s,\mathcal{O}}\bigg( \sum_{L \in \mathbb{N} + l_{k,\alpha}(b) + l_{k,\alpha}(b_1)} \!\!\! \frac{(2 + L)^{t_k + \varsigma_{\mathcal{O}}}\, \#\big\{\gamma \in S_{\Sigma}^{\circ}\quad |\quad  L < l_{k,\alpha}(\gamma) \leq L + 1\big\}}{\big(1 + [L - l_{k,\alpha}(b) - l_{k,\alpha}(b_1)]_{+}\big)^{s}}\bigg)\prod_{\substack{\beta \in \pi_0(\partial \Sigma) \\ \beta \neq b_1,b}} (1 + l_{k,\alpha}(\beta))^{t_k + \varsigma_{\mathcal{O}}} \\
& \leq N_k\,M'_{k,s,\mathcal{O}} \bigg(\sum_{L > 0} \frac{(1 + L + l_{k,\alpha}(b_1) + l_{k,\alpha}(b))^{t_k + \varsigma_{\mathcal{O}} + d_k}}{L^{s}}\bigg) \prod_{\substack{\beta \in \pi_0(\partial \Sigma) \\ \beta \neq b_1,b}} (1 + l_{k,\alpha}(\beta))^{t_{k} + \varsigma_{\mathcal{O}}},
\end{split} 
\end{equation*}
where we have used the polynomial growth axiom between the second and third line. We can choose $s = t_k + \varsigma_{\mathcal{O}} + d_{k} + 2$ to make the series between the bracket convergent, and there is a constant $M''_{k,\mathcal{O}} > 0$ (also dependent on $ \varsigma_{\mathcal{O}}$) such that the previous expression is bounded by
\[
M''_{k,\mathcal{O}} \prod_{\beta \in \pi_0(\partial \Sigma)} (1 + l_{k,\alpha}(\beta))^{t_k + \varsigma_{\mathcal{O}} + d_{k}}.
\]
Therefore, there exists a constant $M'''_{k,\mathcal{O}} > 0$ (also dependent on $ \varsigma_{\mathcal{O}}$) such that
\begin{equation}
\label{theineq2222}\sum_{[P] \in \mathcal{P}_{\Sigma}^b} \big|\Theta_{P}^{k}(B_{P}^b,\Omega_{\Sigma - P})\big|_{k,\alpha} \leq M'''_{k,\mathcal{O}} \prod_{\beta \in \pi_0(\partial \Sigma)} (1 + l_{k,\alpha}(\beta))^{t_k + \varsigma_{\mathcal{O}} + d_{k}}.
\end{equation}
The other $\Gamma_{\Sigma}^{\partial}$-orbits $\mathcal{O} \subseteq \mathcal{P}_{\Sigma}^{\emptyset}$ can be treated in a similar fashion, replacing $B_{P}^{b}$ with $C_{P}$ and using the second part of the small pair of pants axiom and the decay assumption \eqref{decC}. The result is
\begin{equation}
\label{theineq3222}\sum_{[P] \in \mathcal{O}} \big|\Theta_{P}^{k}(C_P,\Omega_{\Sigma - P})\big|_{k,\alpha} \leq M'''_{k,\mathcal{O}} \prod_{\beta \in \pi_0(\partial \Sigma)} (1 + l_{k,\alpha}(\beta))^{t_k + \varsigma_{\mathcal{O}} + d_k}
\end{equation}
for some constant $M'''_{k,\mathcal{O}} > 0$.

In particular, for any orbit $\mathcal{O}$ in $\mathcal{P}_{\Sigma}$, the series $\sum_{[P] \in \mathcal{P}_{\Sigma}^{b}} \Theta_{P}^k(X_{P}^{\mathcal{O}},\Omega_{\Sigma - P})$ (where $X$ is $B^b$ or $C$ depending on the orbit $\mathcal{O}$), is absolutely convergent in $E^{k}(\Sigma)$. Let us denote $w_{\mathcal{O}}^{k}$ the limit.  Since the maps we have used are compatible with the restriction morphisms over $k \in \mathscr{I}_{\Sigma}$, there exists a unique $w_{\mathcal{O}} \in E(\Sigma)$ such that $w_{\mathcal{O}}^k = \rho^{k}(w_{\mathcal{O}})$. The inequalities \eqref{theineq2222}-\eqref{theineq3222} imply that
\[
\|w_{\mathcal{O}}\|_{k,t_k + \varsigma_{\mathcal{O}} + d_k} < + \infty,
\]
thus $w_{\mathcal{O}} \in {}^{\flat}E(\Sigma)$. Summing over the finitely many $\Gamma_{\Sigma}^{\partial}$-orbits we have a well-defined element
\begin{equation}
\label{GRorbit} \Omega_{\Sigma} = \bigg(\sum_{\mathcal{O} \in \mathcal{P}_{\Sigma}/\Gamma_{\Sigma}^{\partial}} w_{\mathcal{O}}\bigg) \in {}^{\flat}E(\Sigma),
\end{equation}
which is by definition the limit of the series in \eqref{defGR}.

Let $\varphi:\,\Sigma \rightarrow \Sigma'$ be a morphism in $\Surf_{1}$. Since the linear map $E(\varphi)$ is continuous, $E(\varphi)(\Omega_{\Sigma})$ coincides with the series \eqref{defGR} to which one applies $E(\varphi)$ term by term. For $[P] = [f:\,P \rightarrow \Sigma] \in \mathcal{P}_{\Sigma}$, $\varphi$ induces a morphism $\tilde{\varphi}$ in $\Surf_{1}$ between $P \cup \overline{\Sigma \setminus f(P)}$ and $P \cup \overline{\Sigma \setminus \varphi \circ f(P)}$. We denote $[\tilde{\varphi}(P)]$ the class of $[\varphi \circ f:\,P \rightarrow \Sigma'] \in \mathcal{P}_{\Sigma'}$. The functoriality of the glueing morphisms implies that 
\[
E(\varphi)\big(\Theta_{P}(X_{P},\Omega_{\Sigma - P})\big) = \Theta_{\tilde{\varphi}(P)}\big(E(\tilde{\varphi})(X_{P},\Omega_{\Sigma - P})\big).
\]
We apply this formula to $X_{P} = B_{P}^b$ or $C_{P}$ which are functorially attached to pairs of pants. The induction hypothesis guarantees that $\Omega_{\Sigma - P}$ is also functorial. Therefore,
\begin{equation}
\label{indj}E(\varphi)\big(\Theta_{P}(X_{P},\Omega_{\Sigma - P})\big) = \Theta_{\tilde{\varphi}(P)}\big(X_{\tilde{\varphi}(P)},\Omega_{\Sigma - \tilde{\varphi}(P)}\big).
\end{equation}
Note that if $X = B^{b}$, by $X_{\tilde{\varphi}(P)}$ we actually mean $B_{\tilde{\varphi}(P)}^{\tilde{\varphi}(b)}$ taking into account that the mapping class $\varphi$ can permute the components in $\partial_{+}\Sigma$.

Through \eqref{indj} we observe a bijection between the terms in $E(\varphi)(\Omega_{\Sigma})$ and the terms in $\Omega_{\Sigma'}$. Thanks to absolute convergence of the series for each seminorm, we deduce that
\[
E(\varphi)(\Omega_{\Sigma}) = \Omega_{\Sigma'}\,,
\]
and thus we have proved functoriality. In particular, taking $\Sigma' = \Sigma$ gives the property of mapping class group invariance $E(\varphi)(\Omega_{\Sigma}) = \Omega_{\Sigma}$. We conclude the proof by induction.
\hfill $\Box$

\subsection{Natural transformations of target theories}
\label{MainProp1}
 
Let $\mathbf{1}:\,\Surf_{1} \rightarrow \mathcal{C}$ be the symmetric monoidal functor which assigns $\mathbb{K}$ to any surface in $\Surf_{1}$. A functorial assignment in $E$ is equivalent to the data of a natural transformation from the functor $\mathbf{1}$ to the functor $E$. For any admissible data, the GR amplitudes provide such a natural transformation.

More generally, let $E$ and $\tilde{E}$ be two target theories, and $\eta:\,E \Longrightarrow \tilde{E}$ be a natural transformation compatible with the union morphisms, the glueing morphisms and the length functions. If
\[
\mathcal{I} := (A_P,B_P^{b},C_{P},D_{T})
\]
is an admissible initial data for $E$, it is clear that
\[
\eta(\mathcal{I}) := (\eta_{P}(A_P),\eta_{P}(B_P^b),\eta_P(C_P),\eta_T(D_T))
\]
are admissible initial data for $\tilde{E}$. Denote $\Omega^{\mathcal{I}}$ and $\tilde{\Omega}^{\eta(\mathcal{I})}$ the $E$-valued (respectively $\tilde{E}$-valued) outcome of GR from these two set of initial data.

\begin{proposition}
\label{GRtoGR} For any surface $\Sigma$ in $\Surf_{1}$, we have $\tilde{\Omega}_{\Sigma}^{\eta(\mathcal{I})} = \eta_{\Sigma}(\Omega_{\Sigma}^{\mathcal{I}})$.
\end{proposition}
\textbf{Proof.} This is directly implied by the fact that $\eta_{\Sigma}:\,E(\Sigma) \longrightarrow \tilde{E}(\Sigma)$ is continuous, linear, and that the assignment $\Sigma \longmapsto \eta_{\Sigma}$ is compatible with union and glueing morphisms and with the length functions. \hfill $\Box$

\subsection{Inducing initial data for tori with one boundary}
\label{toriD}

In quantum field theories, defining correlation functions for tori with one boundary by a cutting procedure often involves a renormalisation procedure to get rid of infinities. We avoid addressing these potential problems in the construction of target theories since we did not include self-glueings in our axioms.

Imagine that we are nevertheless given a functorial morphism, 
\begin{equation}
\label{Xiglue} \Xi:\,E(P) \longrightarrow E(T)\,,
\end{equation}
where $P$ is a pair of pants seen as an object in $\Surf_{1}$, and $T$ is the torus with one boundary obtained by glueing the two boundary components in $\partial_{+}P$.

If $T$ is now an arbitrary object in $\Surf_{1}$ which is a torus with one boundary, for any simple closed curve $\gamma \in S_{T}^{\circ}$, we obtain a homotopy class of pairs of pants $[P_{\gamma}]$ by cutting $T$ along $\gamma$. By definition of cutting/glueing, this $T$ is obtained by self-glueing on $P_{\gamma}$, and we denote by $\Xi_{\gamma}$ the self-glueing morphism coming from \eqref{Xiglue}.

Assume that we are given a functorial assignment $P \mapsto C_{P} \in E(P)$, in particular $C_{P}$ must be $\Gamma_P$-invariant. We remark that, due to the assumed invariance of $C_{P}$ under braiding of the two boundary components of $\partial_+P$, and the assumed functoriality of the self-glueing morphism, $\Xi_{\gamma}(C_{P_{\gamma}})$ does not depend on an ordering of the two last boundary components. We may seek to define
\[
\Omega_{T} := \sum_{\gamma \in S_{T}^{\circ}} \Xi_{\gamma}(C_{P_{\gamma}})  \in E(T).
\]
For this purpose, we introduce an additional axiom for $C$.

\vspace{0.2cm}

\noindent \textsc{Decay for self-glueing.} Assume that we are given a functorial self-glueing morphism for pairs of pants as above. For any object $T$ in $\Surf_{1}$ which is a torus with one boundary component, for any $\gamma \in S_{T}^{\circ}$, for any $k \in \mathscr{I}_{T}$, there exists $t_k \geq 0$ and $j_k \in \mathscr{I}_{P_{\gamma}}$ functorial such that for any $s > 0$, there exists $M_{k,s} > 0$ functorial such that for any $\alpha \in \mathscr{A}_{T}^{k}$ , we have
\begin{equation}
\label{Csefl} \big| \Xi_{\gamma}^{k}(C_{P_{\gamma}})\big|_{k,\alpha} \leq M_{k,s} \frac{(1 + l_{k,\alpha}(\gamma))^{2t_k}(1 + l_{k,\alpha}(\partial T))^{t_k}}{\big(1 + [2l_{k,\alpha}(\gamma) - l_{k,\alpha}(\partial T)]_{+}\big)^{s}}\,.
\end{equation}

\vspace{0.2cm}

By an argument of absolute convergence similar to the one detailed in the proof of Theorem~\ref{converge}, we find
\begin{lemma}
\label{sdddd} Assume $(A,B,C)$ satisfy the properties listed in Definition~\ref{DEFDECAY} and the axiom of decay for self-glueing. Then
\begin{equation}
\label{DTFOR}D_{T} = \sum_{\gamma \in S_{T}^{\circ}} \Xi_{\gamma}(C_{P_{\gamma}})
\end{equation}
is a well-defined, functorial assignment for any object $T$ in $\Surf_{1}$ which is a torus with one boundary component, and $(A,B,C,D)$ is an admissible initial data.
\end{lemma}

Note that the $(A,B,C)$ parts of the initial data are sufficient to define the GR amplitudes $\Omega_{\Sigma}$ for surfaces $\Sigma$ of genus $0$. The $D$ part is necessary to extend this definition to positive genus. Lemma~\ref{sdddd} gives a way to induce a $D$ if we have self-glueing morphisms for pairs of pants and if $C$ satisfying the self-glueing decay axiom. Note that we could replace $C$ by $A$ in the discussion. Note that $B_{P}$ does not a priori have the invariance under braiding of the two boundary components of $\partial_{+}P$. If we insist in using a $B$ to induce a $D$, we should replace \eqref{DTFOR} with the formula
\[
\sum_{\gamma \in S_{T}^{\circ,\rm or}} \Xi_{\gamma}(B_{P_{\gamma}}^{\gamma_+})\,,
\]
where $S_{T}^{\circ,\rm or}$ now enumerates oriented simple closed curves and $\gamma_+$ is the boundary component of $P_{\gamma}$ created to the left of $\gamma$ when we cut $T$.

\section{Teichm\"uller theory background}
\label{SSS21}

In order to build our first family of applications of our general theory presented above, we need to recall a few facts from Teichm\"{u}ller theory.

\subsection{Teichm\"uller spaces}

Let $\Sigma$ be a stable bordered surface. If we consider the space of smooth metrics modulo conformal equivalence on $\Sigma$, we obtain a space of conformal classes of metrics on $\Sigma$, which is a $\text{Diff}_0(\Sigma,\partial\Sigma)$ fiber bundle over Teichm\"{u}ller space $\mathcal{T}_{\Sigma}$. Here $\text{Diff}_0(\Sigma,\partial \Sigma)$ is the group of diffeomorphisms of $\Sigma$, which are isotopic to the identity relatively to the boundary (see e.g. \cite{EarleEells}).

Then, $\mathcal{T}_{\Sigma}$ is in bijection with the set of hyperbolic metrics on $\Sigma$, for which the boundaries are geodesic, modulo diffeomorphisms which are isotopic to the identity relatively to the boundary. We further recall that when $\Sigma$ has type $(g,n)$, $\mathcal{T}_{\Sigma}$ is a smooth manifold of dimension $6g-6+3n$.

We denote $\mathcal{T}_{\partial \Sigma} = \mathbb{R}_{+}^{\pi_0(\partial \Sigma)}$ and $p_{\partial}:\,\mathcal{T}_{\Sigma} \rightarrow \mathcal{T}_{\partial \Sigma}$ the perimeter map. We stress that $\mathbb{R}_{+}$ is the real positive axis, excluding $0$. If $L$ is an assignment of positive real numbers to the boundary components of $\Sigma$, we denote $\mathcal{T}_{\Sigma}(L) = p_{\partial}^{-1}(L)$.

We can of course also describe the Teichm\"uller space in terms of complex structures on $\Sigma$, e.g. $\mathcal{T}_{\Sigma}$ is also the set of equivalence classes of diffeomorphisms $\mu$ from $\Sigma$ to a bordered Riemann surface $\mathcal{S}$. If $\mu_i:\,\Sigma \rightarrow \mathcal{S}_{i}$ with $i = 1,2$ are two diffeomorphisms as above, we declare them equivalent if there exists a biholomorphic map $\Phi:\,\mathcal{S}_{1} \rightarrow \mathcal{S}_{2}$, such that $\mu_{2}^{-1}\circ\Phi\circ\mu_1$ is isotopic to the identity on $\Sigma$ among such diffeomorphisms, which are the identity on the boundary. When $\Sigma$ is stable, the pure mapping class group $\Gamma_{\Sigma}^{\partial}$ acts on the Teichm\"uller spaces $\mathcal{T}_{\Sigma}$ and $\mathcal{T}_{\Sigma}(L)$ properly discontinuously, possibly with finite stabilisers: the quotients are the moduli spaces $\mathcal{M}_{\Sigma}$ and $\mathcal{M}_{\Sigma}(L)$.

\subsection{Bounds on the number of multicurves}

If $\sigma$ is a hyperbolic metric on $\Sigma$ and $\gamma$ is a simple closed curve, there is a unique shortest geodesic in the homotopy class of $\gamma$, and we denote $\ell_{\sigma}(\gamma)$ its length. The length of a multicurve is by definition the sum of lengths of its components.

The length spectrum is the sequence $(l_{\sigma,i})_{i \geq 1}$ of lengths of isotopy classes of simple closed curves (which are not isotopic to the boundary) in $\Sigma$, in weakly increasing order. The systole ${\rm sys}_{\sigma} := l_{\sigma,1}$ is the shortest of these lengths. We will exploit the following result.
\begin{lemma} 
\label{ThPar} Let $\epsilon \in (0,1)$. Let $\sigma$ be a hyperbolic metric on a connected bordered surface $\Sigma$ with non-zero boundary lengths. For any $t \in [0,\ell_{\sigma}(\partial \Sigma))$ there exists another hyperbolic metric $\tilde{\sigma}$ on $\Sigma$ such that
\begin{itemize}
\item[$\bullet$] $t = \ell_{\tilde{\sigma}}(\partial \Sigma) < \ell_{\sigma}(\partial \Sigma)$.
\item[$\bullet$] $l_{\tilde{\sigma},i} \leq l_{\sigma,i}$ for any $i \geq 1$.
\item[$\bullet$] if $l_{\sigma,i} < \epsilon$, then $l_{\tilde{\sigma},i} = l_{\sigma,i}$.
\item[$\bullet$] if $l_{\sigma,i} \geq \epsilon$, then $l_{\tilde{\sigma},i} \geq \epsilon$.
\end{itemize} 
\end{lemma} 
\noindent \textbf{Proof.} We use a slight modification of the proof of \cite[Theorem 3.3]{Parlierreduce}, which Hugo Parlier communicated to us. By the collar lemma, two simple closed geodesics of length $\epsilon \leq 1 \leq 2\,\ln(1 + \sqrt{2})$ cannot intersect. Let $\sigma(0)$ be a hyperbolic metric on $\Sigma$. We introduce the set $\mathcal{S}_{0}$ of $\gamma \in S_{\Sigma}^{\circ}$ such that $\ell_{\sigma(0)}(\gamma) \leq \epsilon$, and denote $\Sigma_0 := \Sigma$. We denote $\Sigma'_{0}$ the surface obtained by cutting $\Sigma_0$ along each of the curves which belong to $\mathcal{S}_0$. We are going to construct for each $t \in [0,\ell_{\sigma}(\partial \Sigma)]$ a finite set of simple closed curves $\mathcal{S}_{t}$ and hyperbolic metric $\sigma(t)$ on $\Sigma$ such that
\begin{equation}
\label{condsigm}\ell_{\sigma(t)}(\partial \Sigma) = \ell_{\sigma}(\partial \Sigma) - t,\qquad \forall i \geq 1,\quad l_{\sigma(t),i} \leq l_{\sigma(0),i}.
\end{equation}
In this process we will always denote $\Sigma'_{t}$ the surface $\Sigma_{t}$ cut along the geodesics that represent the elements of $\mathcal{S}_{t}$. We start to decrease the length as done in \cite[Theorem 3.3.]{Parlierreduce} on each boundary component of the connected components of $\Sigma'_0$ which do not belong to $\mathcal{S}_0$, so that the length spectrum of $\Sigma_{0}$ decreases continuously. This defines a new hyperbolic metric $\sigma(t)$ satisfying \eqref{condsigm} for $t \geq 0$ small enough which keep $\mathcal{S}_{t} = \mathcal{S}_{0}$. If for $t \in [0,\ell_{\sigma}(\partial \Sigma)]$ all simple closed geodesics which are not in $\mathcal{S}_{t}$ have $\sigma(t)$-length $> \epsilon$, the algorithm terminates. Otherwise, there exists a minimal $t^* \in (0,\ell_{\sigma(0)}(\partial \Sigma))$ such that some simple closed geodesic in $\Sigma$ has $\sigma(t^*)$-length $\epsilon$. We add all of those curves to $\mathcal{S}_{t}$ ($= \mathcal{S}_{0}$ if $t < t^*$) to form an updated set $\mathcal{S}_{t^*}$. We then continue to decrease the length of the boundary components of $\Sigma'_{t}$ which are not elements of $\mathcal{S}_{t}$, continuously decreasing the $\sigma(t)$-length spectrum, repeating our update each time we meet simple closed geodesics of $\sigma(t)$-length exactly $\epsilon$. For any $t$, the curves collected in $\mathcal{S}_{t}$ cannot intersect due to the previous observation, and there are at most $3g - 3 + n$ non-intersecting simple closed geodesics in the interior of a surface of genus $g$ with $n$ boundary components. Therefore, $\#\mathcal{S}_{t} \leq 3g - 3 + n$ for any $t$, and $\mathcal{S}$ can only be updated finitely many times. The construction makes sure that for all $t$, the two last requirements hold, while  \eqref{condsigm} meets the two first requirements.  \hfill $\Box$  

\medskip

We will need a well-known estimate on the number of multicurves of bounded length, which we can make uniform in the boundary length using the previous lemma. We also state a similar estimate for multicurves not intersecting a given one, which will be used in Section~\ref{Boundbeh}.

\begin{theorem}
\label{Mirzath} Let $\Sigma$ be a surface of type $(g,n)$ and $\mu$ a fixed primitive multicurve in $\Sigma$ with $k$ components. For any $\varepsilon \in (0,1)$, there exists $N_{\epsilon} > 0$ depending only on $(g,n,\epsilon)$ such that, for any $\sigma \in \mathcal{T}_{\Sigma}$ and $L > 0$
\[
\#\big\{c \in M_{\Sigma} \quad \big| \quad \ell_{\sigma}(c) \leq L\,\,{\rm and}\,\,c \cap \mu = \emptyset\big\} \leq N_{\epsilon}\,L^{6g - 6 + 2(n + k)} \prod_{\substack{\gamma \in S_{\Sigma^{\mu}}^{\circ} \\ \ell_{\sigma}(\gamma) \leq \epsilon}} \frac{1}{\ell_{\sigma}(\gamma)}.
\]
In particular, if $\mu = \emptyset$ and ${\rm sys}_{\sigma} \geq \varepsilon$, the right-hand side is equal to $N_{\epsilon}\,L^{6g - 6 + 2n}$.
\end{theorem} 
\noindent \textbf{Proof.} We first assume $\mu = \emptyset$. We start from the result of Mirzakhani \cite[Proposition 3.6]{Mirzakhanigrowth} stating that, if $(\Sigma,\tilde{\sigma})$ only has punctures (i.e. boundaries of length $0$), for any $\varepsilon \in (0,1)$, there exists a constant $N'_{\epsilon} > 0$ depending only on $(\epsilon,g,n)$ such that
\begin{equation}
\label{BSigmaell} \#\big\{c \in M_{\Sigma} \quad \big| \quad \ell_{\tilde{\sigma}}(c) \leq L\big\} \leq N'_{\epsilon} \,L^{6g - 6 + 2n} \prod_{\substack{\gamma \in S_{\Sigma} \\ \ell_{\tilde{\sigma}}(\gamma) \leq \varepsilon}} \frac{1}{\ell_{\tilde{\sigma}}(\gamma)}.
\end{equation}
If $\sigma$ is a hyperbolic metric on $\Sigma$ with non-zero boundary lengths, Theorem~\ref{ThPar} provides us with $\tilde{\sigma}$ having punctures, such that $l_{i,\tilde{\sigma}} \leq l_{i,\sigma}$ for any $i \geq 1$, and
\[
\prod_{\substack{\gamma \in S_{\Sigma}^{\circ} \\  \ell_{\tilde{\sigma}}(\gamma) \leq \epsilon}} \ell_{\tilde{\sigma}}(\gamma) \geq \epsilon^{3g - 3 + n}\,\prod_{\substack{\gamma \in S_{\Sigma}^{\circ} \\ \ell_{\sigma}(\gamma) \leq \epsilon}} \ell_{\sigma}(\gamma).
\]
The prefactor of $\epsilon$ represent the short curves with respect to $\tilde{\sigma}$ which were not short for $\sigma$: by design they cannot have length shorter than $\epsilon$ and there cannot be more than $3g - 3 + n$ of them. Therefore
\[
\#\big\{c \in M_{\Sigma} \quad \big| \quad \ell_{\sigma}(c) \leq L\big\} \leq N_{\epsilon}\,L^{6g - 6 + 2n} \prod_{\substack{\gamma \in S_{\Sigma} \\ \ell_{\sigma}(\gamma) \leq \epsilon}} \frac{1}{\ell_{\sigma}(\gamma)}
\]
for another constant $N_{\epsilon}$ depending only on $(\epsilon,g,n)$. When $\mu \neq \emptyset$, we apply this result for $\Sigma^{\mu}$ to conclude. \hfill $\Box$

\subsection{Bounds on the number of small pairs of pants}

 \begin{definition}
Given a connected $\Sigma$ in $\Surf_{1}$, we say that $[P] \in \mathcal{P}_{\Sigma}$ is (the homotopy class of) a $\sigma$-small pair of pants if
\[
\ell_{\sigma}(\partial P \cap \Sigma^{\circ}) \leq \ell_{\sigma}(\partial P \cap \partial \Sigma).
\]
In other words, if $[P] \in \mathcal{P}_{\Sigma}^{b}$ for some $b \in \pi_0(\partial_{+}\Sigma)$, this means $\ell_{\sigma}(\gamma_P) \leq \ell_{\sigma}(b_1) + \ell_{\sigma}(b)$, while if $[P] \in \mathcal{P}_{\Sigma}^{\emptyset}$ this means $\ell_{\sigma}(\gamma_P^1) + \ell_{\sigma}(\gamma_{P}^2) \leq \ell_{\sigma}(b_1)$.
\end{definition}

In this paragraph we prove a finiteness result for small pairs of pants.
\begin{proposition}
\label{smallpp} For any $\epsilon \in (0,1)$ and $\sigma \in \mathcal{T}_{\Sigma}$ such that ${\rm sys}_{\sigma} \geq \epsilon$,  the number of $\sigma$-small pairs of pants in a surface $\Sigma$ of type $(g,n)$ with respect to $\sigma$ is upper bounded by a constant $Q_{g,n,\epsilon} > 0$.
\end{proposition}

The proof, to which Parlier contributed, relies on the following lemma which provides a uniform estimate on the number of closed curves, and is a variant of \cite[Lemma 6.6.4]{Buser} proved in \cite[Lemma 2.4]{Parlier1}.
\begin{lemma}
\label{EffeL} Let $(\Sigma,\sigma)$ be a hyperbolic surface of genus $g \geq 2$ without boundaries. The number of primitive closed geodesics of length smaller than $L > 0$ is bounded by $(g - 1)e^{L + 6}$.
\end{lemma}
\noindent \textbf{Proof of Proposition~\ref{smallpp}.}
Let us fix $\sigma \in \mathcal{T}_{\Sigma}$. For each $[P] \in \mathcal{P}_{\Sigma}$, there is a unique representative of $[P]$ embedded with geodesic boundaries in $\Sigma$ of respective lengths $L_1,L_2,L_3$. This hyperbolic pair of pants can be obtained by glueing two hyperbolic right-angled hexagons, whose boundary arcs have successive lengths $(L_1/2,d_{1,2},L_2/2,d_{2,3},L_3/2,d_{1,3})$, with separate identification of the second, fourth and sixth boundary arcs -- which are called ``seams''. Then $d_{i,j}$ is the distance between the $i$-th and $j$-th boundaries in $P$.

Fix $b \in \pi_0(\partial_+\Sigma)$. For $[P] \in \mathcal{P}_{\Sigma}^{b}$, the ordering of boundary components of $P$ is such that $L_1 = \ell_{\sigma}(b_1)$, $L_2 = \ell_{\sigma}(b)$ and $L_3 = \ell_{\sigma}(\gamma_{P})$. We recall from Section~\ref{glueinpant} that $[P]$ is uniquely determined by the free homotopy class (in $\Sigma$) of the seam $c_{P}$ between $b_1$ and $b$ (that is, between the first and the second boundary of $P$), which has length $d_{1,2}$. Hyperbolic trigonometry (see e.g. \cite{Buser}) gives
\begin{equation}
\label{Hch1} {\rm cosh}(d_{1,2}) = \frac{{\rm cosh}(L_3/2) + {\rm cosh}(L_1/2){\rm cosh}(L_2/2)}{{\rm sinh}(L_1/2){\rm sinh}(L_2/2)}.
\end{equation}
Let us fix $\epsilon < 1$ and assume $L_1,L_2 \geq \epsilon$. Using the elementary bounds
\begin{equation*}
\begin{split}
\tfrac{1}{2}e^{x} \leq {\rm cosh}(x) \leq e^{x} & \qquad {\rm for}\,\,x \geq 0, \\
{\rm sinh}(x/2) \geq \tfrac{\epsilon}{2}e^{x/2} & \qquad {\rm for}\,\,x \geq \epsilon, \\
{\rm cotanh}(x) \leq \tfrac{2}{\epsilon} & \qquad {\rm for}\,\,x \geq \epsilon,
\end{split}
\end{equation*}
we deduce that
\[
\frac{e^{d_{1,2}}}{2} \leq \frac{1 + e^{(L_3 - L_2 - L_1)/2}}{(\epsilon/2)^2}.
\]
For $\sigma$-small pairs of pants, we have $L_3 \leq L_2 + L_1$ hence
\[
d_{1,2} \leq 2\ln(4/\epsilon)
\]
since $L_3 \geq {\rm sys}_{\sigma} \geq \epsilon$.
By considering the image of $c_{P}$ in the double of $\Sigma$,  which is a closed surface of genus $g^{d} = 2g + n -1$, the subset of $\sigma$-small pairs of pants in $\mathcal{P}_{\Sigma}^{b}$ injects into the set of primitive closed geodesics of length bounded by $4\ln(4/\epsilon)$. Therefore, there are less than $(2g - 2 + n)e^{6}(4/\epsilon)^{4}$ $\sigma$-small pairs of pants of this type.

For $[P] \in \mathcal{P}_{\Sigma}^{\emptyset}$, the ordering of the boundary is such that $L_1 = \ell_{\sigma}(b_1)$, $L_2 = \ell_{\sigma}(\gamma_P^1)$ and $L_3 = \ell_{\sigma}(\gamma_P^2)$. Likewise, $[P]$ is uniquely determined by the free homotopy class (in $\Sigma$) of the curve $c_P$ that starts from $b_1$, travel in the first hexagon to reach the seam in between $\gamma_P^1$ and $\gamma_P^2$ (the second and third boundaries in $P$) and comes back to $b_1$ through the second hexagon. By hyperbolic trigonometry in one of the pentagon cut out by $c_P$ in the first hexagon, we have
\begin{equation*}
\begin{split}
{\rm cosh}(\ell_{\sigma}(c_P)/2) & = {\rm sinh}(d_{1,2}){\rm sinh}(L_2/2)  \\
& =  \frac{\sqrt{{\rm sinh}^2(L_1/2) + {\rm cosh}^2(L_3/2) + {\rm cosh}^2(L_2/2) + 2\,{\rm cosh}(L_1/2){\rm cosh}(L_2/2){\rm cosh}(L_3/2)}}{{\rm sinh}(L_1/2)}.
\end{split}
\end{equation*}
With elementary bounds we get under the assumption $L_1 \geq \epsilon$
\[
\frac{e^{\ell_{\sigma}(c_P)}}{4} \leq 1 + \frac{e^{L_3 - L_1} + e^{L_2 - L_1} + 2\,e^{(L_3 + L_2 - L_1)/2}}{(\epsilon/2)^2}.
\]
For $\sigma$-small pairs of pants, we have $L_2 + L_3 \leq L_1$ and a fortiori $L_3 \leq L_1$ and $L_2 \leq L_1$. Since $1 < 4/\epsilon$, we obtain
\[
\ell_{\sigma}(c_P) \leq \ln(128/\epsilon^2).
\]
Repeating the previous argument we deduce that there are less than $(2g - 2 + n)e^{6}(128/\epsilon^2)^2$ $\sigma$-small pairs of pants in $\mathcal{P}_{\Sigma}^{\emptyset}$ whenever $\ell_{\sigma}(b_1) \geq \epsilon$, and a fortiori when ${\rm sys}_{\sigma} \geq \epsilon$. \hfill $\Box$

\subsection{Comparing different hyperbolic metrics}

If $\Sigma$ does not have boundaries, the Teichm\"uller distance $d_{T}$ makes $\mathcal{T}_{\Sigma}$ a complete metric space.  We can use $d_{T}$ to compare lengths of curves with respect to different metrics on $\Sigma$, due to the following result of Wolpert.

\begin{theorem} \label{Wolpertl} \cite[Lemma 3.1]{Wolpertlength}.
Let $\Sigma$ be a surface without boundaries, and $\gamma$ a non-null homotopic simple closed curve. Then for any two hyperbolic metrics $\sigma,\sigma'$ on $\Sigma$
\[
e^{-2d_{T}(\sigma,\sigma')} \leq \frac{\ell_{\sigma'}(\gamma)}{\ell_{\sigma}(\gamma)} \leq e^{2d_{T}(\sigma,\sigma')}.
\]
\end{theorem}
The result is also true if $\Sigma$ has boundaries, provided we use the Teichm\"uller distance on $\mathcal{T}_{\Sigma^d}$, where recall that $\Sigma^d$ is the surface without boundary obtained by doubling $\Sigma$ along $\partial\Sigma$.

\begin{remark}
\label{remiso} If $\varphi:\,\Sigma \rightarrow \Sigma'$ is a morphism between surfaces in $\Surf$, then $\varphi$ induces a continuous map $\mathcal{T}_{\Sigma} \rightarrow \mathcal{T}_{\Sigma'}$ which is an isometry for the respective Teichm\"uller distances.
\end{remark}

\section{Functions on Teichm\"uller space}
\label{SfunSctin}
The first non-trivial example of target theories comes from spaces of functions on Teichm\"uller spaces. If $X$ is a topological space, $\mathscr{F}(X)$ denotes either the space  of all functions defined on $X$, the space ${\rm Mes}(X)$ of all measurable functions on $X$, or the space $\mathscr{C}^0(X)$ of all continuous functions on $X$.

\subsection{Target theory and geometric recursion}
\label{TargeTeich} 
If $\epsilon > 0$, we denote
\[
\mathcal{T}_{\Sigma}^{\epsilon} = \big\{\sigma \in \mathcal{T}_{\Sigma}\,\,\,|\,\,\,\forall \gamma \in S_{\Sigma}\,\,\,\,\,\ell_{\sigma}(\gamma) \geq \epsilon\big\}
\]
the $\varepsilon$-thick part of the Teichm\"uller space. In other words, $\sigma$ belongs to the $\epsilon$-thick part if the systole and the length of each boundary component are bounded below by $\epsilon$.

\vspace{0.2cm}

\noindent \textsc{Vector spaces and topology --} We take $\mathscr{I}_{\Sigma} = (0,1)$ and for any $\varepsilon \in (0,1)$, we let $E^{\epsilon}(\Sigma) = \mathscr{F}(\mathcal{T}_{\Sigma}^{\varepsilon})$. If $\sigma \in \mathcal{T}_{\Sigma}^{\epsilon}$, we set
\begin{equation}
\label{SetsF}F_{\epsilon,\sigma} = \bigg\{\sigma' \in \mathcal{T}_{\Sigma}^{\epsilon}\,\,\,\colon\,\,\,\sup_{\gamma \in S_{\Sigma}} \bigg|\ln \frac{\ell_{\sigma'}(\gamma)}{\ell_{\sigma}(\gamma)}\bigg| \leq 2 \bigg\}
\end{equation}
By Theorem~\ref{Wolpertl}, $F_{\epsilon,\sigma}$ contains the ball of radius $1$ centered at $\sigma$ for the Teichm\"uller distance. We equip $E^{\epsilon}(\Sigma)$ with the semi-norms indexed by $\mathscr{A}_{\Sigma}^{\varepsilon} = \mathcal{T}_{\Sigma}^{\epsilon}$ with
\begin{equation}
\label{SeminormsF} |f|_{\epsilon,\sigma} = \sup_{\sigma' \in F_{\epsilon,\sigma}} |f(\sigma')|.
\end{equation}
This makes $E^{\epsilon}(\Sigma)$ a locally convex, complete Hausdorff topological vector space which is functorial in $\Sigma$. The projective limit of these spaces is $E(\Sigma) = \mathscr{F}(\mathcal{T}_{\Sigma})$. Since any compact can be covered by finitely many balls of radius $1$ for the Teichm\"uller distance, the topology on $E(\Sigma)$ is the topology of convergence on any compact, and the map $\rho_{\varepsilon}\,:\,E(\Sigma) \rightarrow E^{\varepsilon}(\Sigma)$ consists in restricting the domain of a function on Teichm\"uller space to the $\varepsilon$-thick part. 

\vspace{0.2cm}

\noindent \textsc{Union morphisms --} Using the canonical homeomorphism $\mathcal{T}_{\Sigma_1 \cup \Sigma_2} \cong \mathcal{T}_{\Sigma_1} \times \mathcal{T}_{\Sigma_2}$, we define the union morphism for $f_i \in E(\Sigma_i)$ and $(\sigma_1,\sigma_2) \in \mathcal{T}_{\Sigma_1 \cup \Sigma_2}$ by
\[
(f_1 \sqcup f_2)(\sigma_1,\sigma_2) = f_1(\sigma_1)f_2(\sigma_2).
\]

\vspace{0.2cm}

\noindent \textsc{Glueing morphisms --} If $[P] \in \mathcal{P}_{\Sigma}$ and $\sigma \in \mathcal{T}_{\Sigma}$, there is a unique representative $f:\,P \rightarrow \Sigma$ of $[P]$ such that $f(P)$ has geodesic boundaries with respect to $\sigma$. We denote $\sigma|_{P}$ and $\sigma|_{\Sigma - P}$ the restriction of the hyperbolic metric $\sigma$ to $f(P)$ and $\overline{\Sigma \setminus f(P)}$, which are indeed hyperbolic metrics for which the boundary components of the two subsurfaces are geodesic. We take as the glueing of $f_1 \in E(\Sigma)$ and $f_2 \in E(\Sigma - P)$ evaluated at $\sigma \in \mathcal{T}_{\Sigma}$
\[
\Theta_{P}(f_1,f_2)(\sigma) = f_1(\sigma|_{P})f_2(\sigma|_{\Sigma - P}).
\]

\begin{lemma}
$\sqcup$ and $\Theta_{P}$ are bilinear morphisms in the category ${\rm Pro}$-$\mathcal{V}$.
\end{lemma}
\noindent \textbf{Proof.} As the case of $\sqcup$ follows similar steps but is simpler, we only discuss the glueing morphism. Due to the inclusion $S_{P} \cup S_{\Sigma - P} \subseteq S_{\Sigma}$, we know that if $\varepsilon \in (0,1)$ and $\sigma$ is in the $\varepsilon$-thick part of $\mathcal{T}_{\Sigma}$, the restrictions $\sigma|_{P}$ and $\sigma|_{\Sigma - P}$ are also in the $\varepsilon$-thick part of their respective Teichm\"uller spaces. We can therefore take $h_{\varepsilon} = (\varepsilon,\varepsilon)$ in the definition~\ref{multilinearmor} of bilinear morphisms and it suffices to check the continuity of the map
\[
\Theta_{P}^{\varepsilon,\varepsilon,\varepsilon}\,:\,E^{\epsilon}(P) \times E^{\varepsilon}(\Sigma - P) \rightarrow E^{\varepsilon}(\Sigma).
\]
If $\sigma \in \mathscr{A}_{\Sigma}^{\epsilon}$, we note that that the image of $F_{\epsilon,\sigma}$ via the restriction maps $q_{P}:\,\mathcal{T}_{\Sigma} \rightarrow \mathcal{T}_{P}$ and $q_{\Sigma - P}:\,\mathcal{T}_{\Sigma} \rightarrow \mathcal{T}_{\Sigma - P}$ are respectively included in $F_{\epsilon,\sigma|_{P}}$ and $F_{\epsilon,\sigma|_{\Sigma - P}}$. Therefore
\[
|\Theta_{P}^{\varepsilon,\varepsilon,\varepsilon}(f_1,f_2)|_{\varepsilon,\sigma} \leq |f_1|_{\varepsilon,\sigma|_{P}}\,|f_2|_{\varepsilon,\sigma|_{\Sigma - P}}
\]
which ensures continuity. \hfill $\Box$

\vspace{0.2cm}

So far, we have obtained a pre-target theory. We use the hyperbolic length to induce length functions turning it into a target theory.

\vspace{0.2cm}

\noindent \textsc{Length functions --} For any $\varepsilon \in (0,1)$ and $\sigma \in \mathscr{A}_{\Sigma}^{\varepsilon}$, we let
\begin{equation}
\label{lengthC0} l_{\varepsilon,\sigma}(\gamma) = \ell_{\sigma}(\gamma).
\end{equation}
As a direct consequence of Theorem~\ref{Mirzath} and Proposition~\ref{smallpp}, we have

\begin{lemma}
The length functions satisfy the polynomial growth, the small pair of pants and the restriction axioms of Section~\ref{STargetT}, and turn $E(\Sigma) = \mathscr{F}(\mathcal{T}_{\Sigma})$ into a target theory.
\end{lemma}

\vspace{0.2cm}

\noindent \textsc{Initial data and admissibility --} An initial data for this target theory amounts to a quadruple $(A,B,C,D)$ where $A,B,C \in \mathscr{F}(\mathcal{T}_{P}) \cong \mathscr{F}(\mathbb{R}_{+}^3)$ such that
\[
X(L_1,L_2,L_3) = X(L_1,L_3,L_2),\qquad X \in \{A,C\},
\]
and $D_T \in \mathscr{F}(\mathcal{T}_{T})$ is a functorial assignment for tori $T$ with one boundary component. By construction of our seminorms \eqref{SetsF}-\eqref{SeminormsF}, the condition of admissibility amounts to requiring, for any $\epsilon > 0$, the existence of $t \geq 0$ such that for all $s \geq 0$, there exists $M_{\epsilon,s} > 0$ for which, for any $L_1,L_2,L_3,\ell,\ell' \geq \epsilon$ and $\sigma \in \mathcal{T}_{T}^{\epsilon}$
\begin{equation}
\label{theinininini} \begin{split}
|A(L_1,L_2,L_3)| & \leq M_{\epsilon,0}\,(1 + L_1)^{t}(1 + L_2)^{t}(1 + L_3)^{t}, \\
|B(L_1,L_2,\ell)| & \leq M_{\epsilon,s}\,\frac{\big((1 + L_1)^{t}(1 + L_2)^{t}}{\big(1 + [L_1 + L_2 - \ell]_{+}\big)^{s}} ,\\
|C(L_1,\ell,\ell')| & \leq M_{\epsilon,s}\,\frac{(1 + L_1)^{t}}{\big(1 + [L_1 - \ell - \ell']_{+}\big)^{s}}, \\
|D_T(\sigma)| &\leq  M_{\epsilon,0}\,\big(1 + \ell_{\sigma}(\partial T)\big)^{t}.
\end{split}
\end{equation}
\textsc{Geometric recursion --} Let us specialise our main result Theorem \ref{converge} to this context. We seek to define a functorial assignment $\Sigma \mapsto \Omega_{\Sigma} \in \mathscr{F}(\mathcal{T}_{\Sigma})$ for any object $\Sigma$ in $\mathcal{B}_{1}$. For pairs of pants $P$ and tori with one boundary component $T$ we let
\[
\Omega_{P}(\sigma) = A(\vec{\ell}_{\sigma}(\partial P)),\qquad \Omega_{T} = D_{T},
\]
where $\vec{\ell}_{\sigma}(\partial P)$ is the triple of boundary lengths of $P$ in which $\partial_-P$ appears first. For disconnected surfaces
\[
\Omega_{\Sigma_1 \cup \cdots \cup \Sigma_k}(\sigma_1,\ldots,\sigma_k) = \prod_{i = 1}^k \Omega_{\Sigma_i}(\sigma_i).
\] 
For connected surfaces $\Sigma$ with $\chi_{\Sigma} \leq -2$ we let
\begin{equation}
\label{GRseriesagain}\Omega_{\Sigma}(\sigma) = \sum_{b \in \pi_0(\partial_{+}\Sigma)} \sum_{[P] \in \mathcal{P}_{\Sigma}^{b}} B(\vec{\ell}_{\sigma}(\partial P))\,\Omega_{\Sigma - P}(\sigma|_{\Sigma - P}) + \frac{1}{2} \sum_{[P] \in \mathcal{P}_{\Sigma}^{\emptyset}} C(\vec{\ell}_{\sigma}(\partial P))\,\Omega_{\Sigma - P}(\sigma|_{\Sigma - P}),
\end{equation}
where $\vec{\ell}_{\sigma}(\partial P)$ is the ordered triple of boundary lengths of $P$.
Theorem~\ref{welldefGR} specialises to
\begin{corollary}
If $(A,B,C,D)$ is an admissible initial data in the above sense, $\Sigma \mapsto \Omega_{\Sigma}$ is a well-defined functorial assignment. More precisely, the series \eqref{GRseriesagain} converge absolutely and uniformly on any compact of $\mathcal{T}_{\Sigma}$, and there exists $u \geq 0$ depending only on the topological type of $\Sigma$, such that for any $\epsilon > 0$, we have a constant $K_{\epsilon} > 0$ depending only on $\epsilon$ and the topological type of $\Sigma$ such that
\begin{equation}
\label{theboundLLL}\forall \sigma \in \mathcal{T}_{\Sigma}^{\epsilon}\qquad  |\Omega_{\Sigma}(\sigma)| \leq K_{\epsilon} \prod_{b \in \pi_0(\partial \Sigma)} (1 + \ell_{\sigma}(b))^{u}.
\end{equation}
\end{corollary} 

This result can be proved in a slightly simpler way than the general Theorem~\ref{welldefGR}. As the strategy to prove it in fact inspired the general scheme of target theories and Theorem~\ref{welldefGR}, it is worth recalling the main steps. The idea is to first prove absolute convergence of the series pointwise in Teichm\"uller space. First, we split off from the sum \eqref{GRseriesagain} the contribution of the small pairs of pants, which are finitely many due to Proposition~\ref{smallpp}. For the remaining sum we know by Theorem~\ref{Mirzath} there exists $O(\lambda^{d})$ terms for which the hyperbolic lengths of curves along which we cut is $\leq \lambda$, for some $d > 0$. Since $B$ and $C$ decay faster than any power law with respect to the lengths of curves along which we cut, we have a uniform control for $\Omega_{\Sigma - P}$ by induction hypothesis, we obtain absolute convergence of the series. As the estimate is uniform for $\sigma$ in any given compact (to see this we use Theorem~\ref{Wolpertl}), we deduce uniform convergence of the GR series on any compact. One can bound the absolute value of the series by a constant depending only on the systole (due to Theorem~\ref{Mirzath}) and the boundary lengths, and it is important in the proof to make sure that the bound \eqref{theboundLLL} is stable under induction.

\subsection{Behavior at the boundary of \texorpdfstring{Teichm\"uller}{Teichmueller} space}

\label{Boundbeh}
We analyse the behaviour of GR amplitudes $\Omega_{\Sigma}(\sigma)$ when $\sigma$ approaches the boundary of $\mathcal{T}_{\Sigma}$, \textit{i.e.} when some curves are pinched. For this we need to assume some uniform control of the initial data over $\mathcal{T}_{\Sigma}$. We first show it results in a uniform control of the GR amplitudes over $\mathcal{T}_{\Sigma}$. If $\epsilon \in (0,1)$, we denote
\[
S_{\Sigma,\sigma}^{\circ,\epsilon} := \big\{\gamma \in S_{\Sigma}^{\circ} \quad \big|\quad \ell_{\sigma}(\gamma) \leq \epsilon\big\}.
\]
Notice that the curves in this set cannot intersect, in particular this set is finite. 

\begin{definition}
\label{strongedeff}We say that initial data $(A,B,C,D)$ are strongly admissible if there exists $t \geq 0$ and for any $s \geq 0$, there exists $M_{s} > 0$ for which for any $L_1,L_2,L_3,\ell,\ell' > 0$, we have
\begin{equation}
\label{ABCstrongs}\begin{split}
|A(L_1,L_2,L_3)| & \leq M_{0}\,(1 + L_1)^{t}(1 + L_2)^{t}(1 + L_3)^{t} ,\\
|B(L_1,L_2,\ell)| & \leq M_{s}\,\frac{(1 + L_1)^{t}(1 + L_2)^{t}}{\big(1 + [\ell - L_1 - L_2]_{+}\big)^{s}} ,\\
|C(L_1,\ell,\ell')| & \leq M_{s}\,\frac{(1 + L_1)^{t}}{\big(1 + [\ell + \ell' - L_1]_{+}\big)^{s}} .
\end{split}
\end{equation}
and there exists $u_{1,1} \geq 0$ such that for any $\epsilon \in (0,1)$, there exists $K_{1,1,\epsilon} > 0$ for which for any $\sigma \in \mathcal{T}_{T}$, we have
\begin{equation}
\label{theDstrong}
|D_{T}(\sigma)| \leq K_{1,1,\epsilon}\,\frac{\big(1 + \ell_{\sigma}(\partial T))^{u_{1,1}}}{\prod_{\gamma \in S_{T,\sigma}^{\circ,\epsilon}} \ell_{\sigma}(\gamma)}.
\end{equation}
We notice that either $S_{T,\sigma}^{\circ,\epsilon}$ is empty or contains a single curve.
\end{definition}

Under certain conditions, we can specify $D$ from the data of $C$.
\begin{lemma}
\label{lemiaugnfgiu}If $(A,B,C)$ satisfy the conditions of \eqref{ABCstrongs}, then setting
\begin{equation}
\label{DseriesC} D_{T}(\sigma) := \sum_{\gamma \in S_{T}^{\circ}} C\big(\ell_{\sigma}(\partial T),\ell_{\sigma}(\gamma),\ell_{\sigma}(\gamma)\big)
\end{equation}
yields a well-defined $\Gamma_{T}$-invariant function on $\mathcal{T}_{T}$ which satisfies \eqref{theDstrong}.
\end{lemma}
\noindent \textbf{Proof.}
Let $L_1 = \ell_{\sigma}(\partial T)$. We write for any $s > 0$
\begin{equation} 
\begin{split}
|D_{T}(\sigma)| & \leq M_s\,(1 + L_1)^{t}\cdot \#\big\{\gamma \in S_{T}^{\circ} \quad \big| \quad 2\ell_{\sigma}(\gamma) \leq L_1\big\} \\
& \quad + M_s\,\sum_{L \in L_1/2 + \mathbb{N}/2} \frac{(1 + L_1)^{t}\cdot \#\big\{\gamma \in S_{T}^{\circ} \quad \big |\quad L \leq \ell_{\sigma}(\gamma) < L + 1/2\big\}}{(1 + 2L - L_1)^{s}}.
\end{split}  
\end{equation}
With Theorem~\ref{Mirzath} for $\mu = \emptyset$, we deduce
\begin{equation}
\label{DTbounfigufngf}
|D_T(\sigma)| \leq \frac{N_{\epsilon}M_s\,(1 + L_1)^{t}}{\prod_{\gamma \in S_{T,\sigma}^{\circ,\epsilon}} \ell_{\sigma}(\gamma)}\bigg\{\frac{L_1^2}4 + \sum_{L > 0} \frac{(L_1 + L)^{2}}{4\,L^{s}}\bigg\}.
\end{equation}
Choosing $s = 4$ makes the right-hand side of \eqref{DTbounfigufngf} convergent, and we deduce the claimed bound. In particular, the series \eqref{DseriesC} is absolutely convergent on any compact of $\mathcal{T}_{T}$. \hfill $\Box$

\vspace{0.2cm}

Later, when we will say that $(A,B,C)$ is a strongly admissible initial data, it will be implicit that we choose $D$ equal to \eqref{DseriesC}. Lemma~\ref{lemiaugnfgiu} shows that this choice indeed makes $(A,B,C,D)$ a strongly admissible initial data in the sense of Definition~\ref{strongedeff}.

\begin{lemma}
\label{lembounddun}Let $(A,B,C,D)$ be a strongly admissible initial data. Then, for any $\Sigma$ in $\Surf_1$ there exists $u \geq 0$ depending only on the topological type of $\Sigma$ such that for $\epsilon \in (0,1)$,  there exist $K_{\epsilon} > 0$ depending only on $\epsilon$ and the topological type of $\Sigma$, such that for any $\sigma \in \mathcal{T}_{\Sigma}$
\[
|\Omega_{\Sigma}(\sigma)| \leq K_{\epsilon}\,\frac{\prod_{b \in \pi_0(\partial \Sigma)} (1 + \ell_{\sigma}(b))^{u}}{\prod_{\gamma \in S_{\Sigma,\sigma}^{\circ,\epsilon}} \ell_{\sigma}(\gamma)}.
\]
\end{lemma}
\noindent \textbf{Proof.} It is enough to prove the result for connected $\Sigma$. The case of pairs of pants and tori with one boundary component is clear from the $A$ and $D$ bounds. In general, consider $\Sigma$ an object of $\Surf_1$ of type $(g,n)$ with $2g - 2 + n \geq 2$, and assume the result holds for all surfaces of Euler characteristic $< 2 - 2g - n$. By the GR formula \eqref{GRseriesagain}, we have for $\sigma \in \mathcal{T}_{\Sigma}$
\[
|\Omega_{\Sigma}(\sigma)| \leq \sum_{[P] \in \mathcal{P}_{\Sigma}} |X_{P}(\vec{\ell}_{\sigma}(\partial P))|\,|\Omega_{\Sigma - P}(\sigma|_{\Sigma - P})|
\]
for $X_{P}$ equal to $B$ or $C$ depending on the type of $[P]$. We will analyse separately various contributions to this sum.

We first remark for any $[P] \in \mathcal{P}_{\Sigma}$, we can write
\begin{equation} 
\label{indhuynrg}|\Omega_{\Sigma - P}(\sigma|_{\Sigma - P})| \leq K'_{\epsilon} \frac{\prod_{\beta \in \pi_0(\partial(\Sigma - P))}  (1 + \ell_{\sigma}(\beta))^{u'}}{\prod_{\gamma \in S_{\Sigma - P,\sigma|_{\Sigma - P}}^{\circ,\epsilon}} \ell_{\sigma}(\gamma)} \leq K'_{\epsilon} \frac{\prod_{\beta \in \pi_0(\partial(\Sigma - P))}  (1 + \ell_{\sigma}(\beta))^{u'}}{\prod_{\gamma \in S_{\Sigma ,\sigma}^{\circ,\epsilon}} \ell_{\sigma}(\gamma)},
\end{equation}
where we can choose constants $u' \geq 0$ depending only on the topological type of $\Sigma$, and $K'_{\epsilon}$ depending only on $\epsilon \in (0,1)$ and the topological type of $\Sigma$ and we have used
\[
S_{\Sigma - P,\sigma|_{\Sigma - P}}^{\circ,\epsilon} \subseteq S_{\Sigma,\sigma}^{\circ,\epsilon},
\] 
and since for all $\gamma \in S_{\Sigma ,\sigma}^{\circ,\epsilon} - S_{\Sigma - P,\sigma|_{\Sigma - P}}^{\circ,\epsilon}$ we have that $\ell_\sigma(\gamma) \leq \epsilon$, so $1<\epsilon^{-1} \leq \ell_\sigma(\gamma)^{-1}$.


We can then use Theorem~\ref{Mirzath} to estimate the number of multicurves with bounded length. This will be used several times in the following inequalities, together with the induction hypothesis \eqref{indhuynrg} and the assumptions on $B$ and $C$.

If $b \in \pi_0(\partial_+\Sigma)$, decomposing the sum over $[P] \in \mathcal{P}_{\Sigma}^{b}$ into small and non-small pairs of pants, we obtain for any $s > 0$ that there exist constant $N_{\varepsilon}>0$ depending as well on $(g,n)$ such that
\begin{equation}
\label{nfgsigunfgi}\begin{split}
& \quad \sum_{[P] \in \mathcal{P}_{\Sigma}^{b}} |B(\ell_{\sigma}(b_1),\ell_{\sigma}(b),\ell_{\sigma}(\gamma_P))|\,|\Omega_{\Sigma - P}(\sigma|_{\Sigma - P})| \\ 
& \leq K'_{\varepsilon} \prod_{\substack{\beta \in \pi_0(\partial \Sigma) \\ \beta \neq b_1,b}} (1 + \ell_{\sigma}(\beta))^{u'} \cdot (1 + \ell_{\sigma}(b_1))^{t}(1 + \ell_{\sigma}(b))^{t} \cdot \frac{N_{\varepsilon}}{\prod_{\gamma \in S_{\Sigma,\sigma}^{\circ,\epsilon}} \ell_{\sigma}(\gamma)}  \\
& \quad \cdot M_s \bigg\{(\ell_{\sigma}(b_1) + \ell_{\sigma}(b))^{6g - 6 + 2n}(1 + \ell_{\sigma}(b_1) + \ell_{\sigma}(b))^{u'} + \sum_{L \in \ell_{\sigma}(b_1) + \ell_{\sigma}(b) + \mathbb{N}} \frac{(L + 1)^{6g - 6 + 2n}(L + 2)^{u'}}{\big(1 + [L  - \ell_{\sigma}(b_1) - \ell_{\sigma}(b)]_{+}\big)^{s}}\bigg\}.
\end{split}
\end{equation}

Choosing a value $s > (6g - 6 +2n + u') + 1$ makes the sum over $L$ convergent and the right-hand side of \eqref{nfgsigunfgi} is bounded by
\begin{equation}
\label{thegoodbund} K_{\epsilon}''\, \frac{\prod_{\beta \in \pi_0(\partial \Sigma)} (1 + \ell_{\sigma}(\beta))^{u}}{\prod_{\gamma \in S_{\Sigma,\sigma}^{\circ,\epsilon}} \ell_{\sigma}(\gamma)} 
\end{equation}
for $u = u' + t + (6g - 6 + 2n)$ and some $K''_{\epsilon} > 0$ which only depends on $(g,n,\varepsilon)$. A similar argument leads to bound the sum over $[P] \in \mathcal{P}_{\Sigma}^{\emptyset}$ by \eqref{thegoodbund} for a perhaps larger constant $K_{\epsilon}$. Thus we get that $|\Omega_{\Sigma}(\sigma)|$ is bounded by \eqref{thegoodbund} for a perhaps larger constant $K_{\epsilon}$. This completes the proof by induction. \hfill $\Box$
 
 \vspace{0.2cm}

We would like to isolate the dominant contribution of the GR sums when some curves are pinched. As we mainly want to illustrate the mechanism, we will impose stronger assumptions on the initial data that facilitate the analysis. First, we will consider that $D$ comes from a $C$ as in \eqref{DseriesC}. If it was not the case, we could easily extend our analysis by adopting suitable assumption for $D$. Let us introduce a notation for the set of ordered pairs of pants decomposition that appear by unfolding the GR sum. Namely, we introduce the functorial assignment $\Sigma \longmapsto \mathcal{Q}_{\Sigma}$ from $\Surf_1$ to the category of countable sets that is uniquely defined by the following properties
\begin{itemize}
\item[$\bullet$] $\mathcal{Q}_{P} = \{P\}$ and $\mathcal{Q}_{T}$ is the set of homotopy classes of oriented simple closed curves in $T^{\circ}$.
\item[$\bullet$] disjoint unions of surfaces are sent to cartesian products of sets.
\item[$\bullet$] if $\Sigma$ is connected with Euler characteristic $\leq -2$, we have
\[
\mathcal{Q}_{\Sigma} = \bigsqcup_{[P] \in \mathcal{P}_{\Sigma}} \mathcal{Q}_{\Sigma - P}.
\]
\end{itemize}
Any element $Q \in \mathcal{Q}_{\Sigma}$ determines
\begin{itemize}
\item[$\bullet$] a decomposition $(Q_1,\ldots,Q_{3g - 3 + n})$ of $\Sigma$ into pairs of pants with ordered boundary components.
\item[$\bullet$] a primitive multicurve $c_{Q}$ with ordered components, consisting of the boundary components of the $Q_i$s that are not boundary components of $\Sigma$, with their order prescribed by their appearance in the recursive construction of $Q$.
\item[$\bullet$] a type map $X \colon \{1,\ldots,3g - 3 + n\} \rightarrow \{A,B,\tfrac{1}{2}C\}$.
\end{itemize}
The GR formula can then be unfolded into
\begin{equation}
\label{GRsumfggggun}
\Omega_{\Sigma}(\sigma) = \sum_{Q \in \mathcal{Q}_{\Sigma}} \prod_{i = 1}^{3g - 3 + n} X_{i}(\sigma|_{Q_i}).
\end{equation}
If $\mu \in M'_{\Sigma}$ we denote $\mathcal{Q}_{\Sigma}[\mu]$ the set of $Q \in \mathcal{Q}_{\Sigma}$ such that $\mu$ does not intersect $c_{Q}$, which is tantamount to saying that $\mu \subseteq c_{Q}$.

\begin{proposition}
\label{pinchth} Let $(A,B,C,D)$ be initial data and $\daleth\, \colon\, \mathbb{R}_{\geq 0} \rightarrow \mathbb{R}_{+}$ be a decreasing function such that
\begin{equation}
\label{subexpdaleth} \daleth(x + y) \leq \daleth(x)\daleth(y),\qquad \daleth(1) < 1.
\end{equation}
We assume there exists $t,M \geq 0$ such that, for any $L_1,L_2,L_3,\ell,\ell' > 0$
\begin{equation}
\label{strongcontrol}
\begin{split}
|A(L_1,L_2,L_3)| & \leq M\,(1 + L_1)^{t}(1 + L_2)^{t}(1 + L_3)^{t}, \\
|B(L_1,L_2,\ell)| & \leq M\,(1 + L_1)^{t}(1 + L_2)^{t}\,\daleth(\ell), \\
|C(L_1,\ell,\ell')| & \leq  M\,(1 + L_1)^{t}\,\daleth(\ell)\daleth(\ell'),
\end{split}
\end{equation} 
and that $D$ is specified by \eqref{DseriesC}.

Let $\Sigma$ be an object of type $(g,n)$ in $\Surf_1$, $\mu$ be a primitive multicurve in $\Sigma$ and introduce
\[
\Omega_{\Sigma}^{[\mu]}(\sigma) = \sum_{Q \in \mathcal{Q}_{\Sigma}[\mu]}  \prod_{i = 1}^{3g - 3 + n} X_{i}(\sigma|_{Q_i}).
\]
 There exists $u$ depending only on $(g,n)$ and for each $\varepsilon \in (0,1)$ a constant $K_{\epsilon} > 0$ depending only on $(\epsilon,g,n)$ such that, for any $\sigma \in \mathcal{T}_{\Sigma}$ such that all simple closed geodesics on $\Sigma^{\circ}_\sigma$ with length $\leq \epsilon$ must be components of $\mu$, we have that
\begin{equation}
\label{thedesiredbound}
|\Omega_{\Sigma}^{[\mu]}(\sigma)| \leq K_{\epsilon}\,\prod_{b \in \pi_0(\partial \Sigma)} (1 + \ell_{\sigma}(b))^{u},
\end{equation}
and
\begin{equation}
\label{thedesirederror}
\big|\Omega_{\Sigma}(\sigma) - \Omega_{\Sigma}^{[\mu]}(\sigma)\big| \leq K_{\epsilon} \prod_{b \in \pi_0(\Sigma)} (1 + \ell_{\sigma}(b))^{t} \cdot \bigg\{\sum_{\nu \subset \mu} \prod_{m \in \pi_0(\mu - \nu)} \big(1 + 4\ln(1/\ell_{\sigma}(m))\big)^{6g - 6 + 2n} \frac{\daleth^4\big(\ln(1/\ell_{\sigma}(m))\big)}{\ell_{\sigma}(m)}\bigg\}.
\end{equation}
\end{proposition}

Equation~\eqref{subexpdaleth} in particular implies that $\daleth(l)$ is bounded for $l \geq 0$, and if $\daleth(1) < e^{-1/4}$, then
\[
\forall s \geq 0 \qquad \lim_{l \rightarrow 0} \frac{\daleth^{4}\big(\ln(4/\ln(l))\big)}{l}\,\ln^{s}(1/l) = 0,
\]
so the error bound in \eqref{thedesirederror} tends to $0$ when $\ell_{\sigma}(\mu) \rightarrow 0$.

When $\mu$ has $3g - 3 + n$ components, it determines a pair of pants decomposition and $\Omega_{\Sigma}^{(\mu)}$ contains only finitely many terms of \eqref{GRsumfggggun}, which can only differ by the ordering of the pairs of pants in this decomposition and the ordering of their boundary components -- these orderings are determining the type $(A,B,C)$ of each factor. In the situation of Proposition~\ref{pinchth}, one can therefore observe that  GR produces a mapping class group invariant function globally defined on $\mathcal{T}_{\Sigma}$ and which has prescribed asymptotic behavior, when a maximal number of curves are pinched, completely specified by the functions $(A(L_1,L_2,L_3),B(L_1,L_2,0),C(L,0,0))$.

\vspace{0.2cm}

\noindent \textbf{Proof of Proposition~\ref{pinchth}.} For each $Q \in \mathcal{Q}_{\Sigma}$, the corresponding term in \eqref{GRsumfggggun} is a function of the lengths of the boundary components of $Q_1,\ldots,Q_{3g - 3 + n}$, which is bounded by the product of polynomials factors in the length of the $\partial Q_i$s and $\daleth$ factors. We observe that each component $\delta$ of $c_{Q}$ bounds at most two pairs of pants $Q_i$ and $Q_j$ for some $i,j \in \{1,\ldots,3g - 3 + n\}$. If $i = j$, then the length of $\delta$ only appears in a factor $C(\ell_{\sigma}(\delta'),\ell_{\sigma}(\delta),\ell_{\sigma}(\delta))$ for some other $\delta'$, and its absolute value is bounded by
\[
M\,(1 + \ell_{\sigma}(\delta'))^{t}\,\daleth(\ell_{\sigma}(\delta))^2 \leq M\,\daleth(0)\,(1 + \ell_{\sigma}(\delta'))^{t}\,\daleth(\ell_{\sigma}(\delta)).
\]
If $i \neq j$, we can always assume $i < j$, and $\ell_{\sigma}(\delta)$ appears only as a variable of $X_i$ resulting into a factor of $\daleth(\ell_{\sigma}(\delta))$ in the upper bound, and as a variable of $X_j$ resulting into a factor $(1 + \ell_{\sigma}(\delta))^{t}$ in the upper bound. As a result, we have
\begin{equation}
\label{thepolytrnigug}
\bigg|\prod_{i = 1}^{3g - 3 + n} X_{i}(\sigma|_{Q_i})\bigg| \leq M'\,\prod_{b \in \pi_0(\partial \Sigma)} (1 + \ell_{\sigma}(b))^{t} \prod_{\delta \in \pi_0(c_Q)} (1 + \ell_{\sigma}(\delta))^{t}\,\daleth(\ell_{\sigma}(\delta))
\end{equation}
for some constant $M' > 0$ depending only on $g$ and $n$.

Let $\nu \subseteq \mu$. If $Q \in \mathcal{Q}_{\Sigma}[\nu]$, we consider the multicurve $c'_Q[\nu] = \sum_{\delta \in \pi_0(c_Q \setminus \nu)} \delta$. Introduce now the set
$$ \mathcal{Q}_\Sigma[\nu,\mu] = \{ Q \in \mathcal{Q}_{\Sigma}[\nu] \mid c_Q \cap m \neq \emptyset, \ \ \forall m\in \pi_0(\mu \setminus \nu)\}$$
of pairs of pants decompositions containing $\nu$ but none of the components of $\mu \setminus \nu$. Assume from now on  that $Q\in \mathcal{Q}_\Sigma[\nu,\mu].$
If $\delta \in \pi_0(c'_Q[\nu])$ and $m \in \pi_0(\mu \setminus \nu)$ we denote $I_{\delta,m}$ the number of intersections of $\delta$ and $m$, and the aforementioned constraint means that
\begin{equation}
\label{thewtwot}\sum_{\delta \in \pi_0(c'_Q[\nu])} I_{\delta,m} \geq 2.
\end{equation}
By the collar lemma we must have for any $\delta \in  \pi_0(c'_Q[\nu])$ the upper bound
\[
\ell_{\sigma}(\delta) \geq \sum_{m \in \pi_0(\mu - \nu)} I_{\delta,m}\,w(\ell_{\sigma}(m)),\qquad w(l) = 2\,{\rm arcsinh}\bigg(\frac{1}{{\rm sinh}(l/2)}\bigg).
\] 
We have $w(l) \geq 2\ln(4/l)$ for all $l > 0$. In particular
\[
\ell_{\sigma}(\delta) \geq \Lambda_{\sigma,\delta} := \sum_{m \in \pi_0(\mu \setminus \nu)} 2\,I_{\delta,m}\,\ln(4/\ell_{\sigma}(m)),
\]
and
\[
\ell_{\sigma}(c'_Q[\nu]) \geq \Lambda_{\sigma}[\nu] := \sum_{m \in \pi_0(\mu \setminus \nu)} 4\ln(4/\ell_{\sigma}(m)).
\]
Let $k = \#\pi_0(\nu)$ and $d' = 3g - 3 + n - k$. We have
\[
\bigg|\sum_{Q \in \mathcal{Q}_{\Sigma}[\nu, \mu]} \prod_{i = 1}^{3g - 3 + n} X_{i}(\sigma|_{Q_i})\bigg| \leq M'\,2^{kt} \prod_{b \in \pi_0(\partial \Sigma)} (1 + \ell_{\sigma}(b))^{t} \!\!\!\!\!\!\!\!\!\!\!\! \sum_{\substack{ L_1,\ldots,L_{d'} \geq 0 \\ \delta_1,\ldots,\delta_{d'} \\ L_i \leq \ell_{\sigma}(\delta_i) - \Lambda_{\sigma,\delta_i} < L_i + 1}} \!\!\!\!\!\!\! \prod_{i = 1}^{d'} (2 + \Lambda_{\sigma,\delta_i} + L_i)^{t}\,\daleth\big(\Lambda_{\sigma,\delta_i} + L_i\big),
\]
where we sum over multicurves $\delta$ with ordered boundary components $\delta_1,\ldots,\delta_k$ which do not intersect $\nu$. The prefactor $2^{kt}$ comes from the polynomial factor in \eqref{thepolytrnigug} and the fact that components of $\nu \subseteq \mu$ have length $\leq \epsilon \leq 1$. We then use the sub-exponentiality of $\daleth$ and \eqref{thewtwot}, and Theorem~\ref{Mirzath} to estimate the number of multicurves of length $\leq \sum_{i = 1}^{d'} (1 + L_i + \Lambda_{\sigma}[\nu])$, noticing that the ordering of the multicurve only contribute to an extra factor that depends only on $(g,n)$. The result is an upper bound of the form
\begin{equation*}
\begin{split}
\bigg|\sum_{Q \in \mathcal{Q}_{\Sigma}[\nu,\mu]} \prod_{i = 1}^{3g - 3 + n} X_{i}(\sigma|_{Q_i})\bigg| & \leq M''\,N_{\epsilon}\,\prod_{b \in \pi_0(\partial \Sigma)} (1 + \ell_{\sigma}(b))^{t} \prod_{m \in \pi_0(\mu \setminus \nu)} \frac{\daleth^{4}\big(\ln(4/\ell_{\sigma}(m))\big)}{\ell_{\sigma}(m)}\ \\
& \quad \cdot \bigg\{ \sum_{L_1,\ldots,L_{d'} \geq 0}\,\, \prod_{i = 1}^{d'} \daleth(L_i)\big(2 + L_i + \Lambda_{\sigma}[\nu])^{t}\Big(\sum_{i = 1}^{d'} (\Lambda_{\sigma}[\nu] + 1 + L_i)\Big)^{6g - 6 + 2n}\bigg\}
\end{split}
\end{equation*}
for some constant $M''$ depending only on $(g,n)$. We observe that, since the subexponentiality and the condition $\daleth(1) < 1$ implies that for any $s \geq 0$
\[
\sum_{L > 0} \daleth(L) L^{s} \leq \sum_{L > 0} \daleth(1)^{L}\,L^{s} < +\infty.
\]
Hence, using that $\daleth$ is decreasing to replace $4/\ell_{\sigma}(m)$ by $1/\ell_{\sigma}(m)$,
\[
\bigg|\sum_{Q \in \mathcal{Q}_{\Sigma}[\nu]} \prod_{i = 1}^{3g - 3 + n} X_{i}(\sigma|_{Q_i})\bigg|  \leq K_{\epsilon}\,(1 + \Lambda_{\sigma}[\nu])^{6g - 6 + 2n}\,\prod_{b \in \pi_0(\partial \Sigma)} (1 + \ell_{\sigma}(b))^{t} \prod_{m \in \pi_0(\mu - \nu)} \frac{\daleth^4\big(\ln(1/\ell_{\sigma}(m))\big)}{\ell_{\sigma}(m)}
\]
for some constant $K_{\epsilon}$ depending only on $(g,n,\epsilon)$. In particular, if $\nu = \mu$, since $\mathcal{Q}_{\Sigma}[\nu, \mu] = \mathcal{Q}_{\Sigma}[\mu]$ and $\Lambda_{\sigma}[\mu] = 0$, we get the desired bound \eqref{thedesiredbound}, and summing over the strict subsets $\nu \subset \mu$ gives the desired error bound \eqref{thedesirederror} after we use the crude upper bound
\[
\bigg(1 + \sum_{m \in \pi_0(\mu - \nu)} 4\ln(4/\ell_{\sigma}(m))\bigg)^{6g - 6 + 2n} \leq (1 + 4d'\ln 4)^{6g - 6 + 2n} \prod_{m \in \pi_0(\mu \setminus \nu)} (1 + \ln(1/\ell_{\sigma}(m)))^{6g - 6 + 2n}.
\]
\hfill $\Box$

\subsection{Integration and topological recursion}

\label{Sinrgng}
In this paragraph, we study the integration of GR amplitudes valued in $E(\Sigma) = {\rm Mes}(\mathcal{T}_{\Sigma})$ over the moduli space of bordered Riemann surfaces with fixed boundary lengths, with respect to the Weil--Petersson measure $\mu_{{\rm WP}}$. This is the volume form associated with Weil--Petersson symplectic form $\omega_{{\rm WP}}$. Our normalisation convention is $\omega_{{\rm WP}} = \sum_{j} \dd l_j \wedge \dd \tau_j $ in Fenchen--Nielsen coordinates \cite{WolpertWPFN}. For convenience of notations, if $\Sigma \mapsto \Omega_{\Sigma}$ is a functorial assignment, we will write $\Omega_{g,n}$ for the (uniquely defined) function on $\mathcal{M}_{g,n} \cong \mathcal{M}_{\Sigma_{g,n}}$ that is induced by $\Omega_{\Sigma_{g,n}}$ where $\Sigma_{g,n}$ is a bordered surface of genus $g$ with $n$ boundary components labelled $(\partial_i\Sigma)_{i = 1}^n$, which is an object of $\Surf_1$ for the choice $\partial_-\Sigma = \partial_1\Sigma$.

Here we adopt slightly weaker assumptions than the strong admissibility of Definition~\ref{Boundbeh}, which is useful for applications \cite{MVpaper}. The difference is that we allow some mild divergence when $\ell,\ell' \rightarrow 0$.
\begin{definition}
Let $\eta \in [0,2)$. We say that $(A,B,C,D)$ is $\eta$-strongly admissible if there exists $t \geq 0$ such that for any $s > 0$, there exists $M_s > 0$ for which, for any $L_1,L_2,L_3,\ell,\ell' > 0$, we have
\begin{equation}
\begin{split}
|A(L_1,L_2,L_3)| & \leq M_0 (1 + L_1)^t (1 + L_2)^t (1 + L_3)^t, \\
|B(L_1,L_2,\ell)| & \leq \frac{M_s}{\ell^{\eta}} \frac{ (1 + L_1)^{t}(1 + L_2)^{t}}{\big(1 + [\ell - L_1 - L_2]_+\big)^s}, \\
|C(L_1,\ell,\ell')| & \leq \frac{M_s}{(\ell\ell')^{\eta}} \frac{(1 + L_1)^t}{(\ell\ell')^{\eta} \big(1 + [\ell + \ell' - L_1]_+\big)^{s}} ,
\end{split}
\end{equation}
and there exists $u_{1,1}$ such that for any $\varepsilon \in (0,1)$, there exists a constant $K_{1,1,\varepsilon}> 0$ for which for any $\sigma \in \mathcal{T}_T$ we have:
$$
|D_T(\sigma)|  \leq K_{1,1,\varepsilon} \frac{\big(1 + \ell_{\sigma}(\partial T)\big)^{u_{1,1}}}{\prod_{\gamma \in S_{T,\sigma}^{\circ,\varepsilon} }\big(\ell_{\sigma}(\gamma)\big)^{\eta/2}} .
$$
\end{definition}

\begin{theorem}
\label{intTRHTT} Let $(A,B,C,D)$ be measurable, $\eta$-strongly admissible initial data, and consider the corresponding GR amplitudes $\Omega_{\Sigma}$. For any $g \geq 0$ and $n \geq 1$ such that $2g - 2 + n > 0$ and $L \in \mathbb{R}_+^{n}$, the function $\Omega_{g,n}$ is integrable over $\mathcal{M}_{g,n}(L)$, and the integrals denoted
\[
V\Omega_{g,n}(L) = \int_{\mathcal{M}_{g,n}(L)} \Omega_{g,n}\,\dd\mu_{{\rm WP}}.
\]
satisfy the topological recursion
\begin{equation}
\label{theteatthtuhtu}\begin{split}
& \quad V\Omega_{g,n}(L_1,\ldots,L_n) \\
& = \sum_{m = 2}^{n} \int_{\mathbb{R}_{+}} B(L_1,L_m,\ell)\,V\Omega_{g,n - 1}(\ell,L_2,\ldots,\widehat{L_m},\ldots,L_n) \ell\, \dd \ell\, \\
& + \frac{1}{2} \int_{\mathbb{R}_{+}^2}C(L_1,\ell,\ell')\bigg(V\Omega_{g - 1, n + 1}(\ell,\ell',L_2,\ldots,L_n) +\!\!\!\!\!\!\!\!\! \!\!\!\sum_{\substack{h + h' = g \\ J \cup J' = \{L_2,\ldots,L_n\}}} \!\!\!\!\!\!\!\!\!\!\!\! V\Omega_{h',1 + \#J}(\ell,J) V\Omega_{h',1 + \#J'}(\ell',J')\bigg) \ell\ell'\, \dd \ell \dd \ell'\,
\end{split}
\end{equation}
with initial data
\[
V\Omega_{0,3}(L_1,L_2,L_3) = A(L_1,L_2,L_3),\qquad V\Omega_{1,1} = V\!D(L) := \int_{\mathcal{M}_{1,1}(L)} D\,\dd\mu_{{\rm WP}}.
\]
\end{theorem}

\vspace{0.2cm}

\noindent \textbf{Proof of Theorem~\ref{intTRHTT}.} For strongly admissible initial data, the integrability follows from the uniform bound of Lemma~\ref{lembounddun} and the argument of Mirzakhani's \cite[p111--112]{Mirzakhanigrowth}. We now sketch the proof under the weaker assumptions of $\eta$-strong admissibility, which is based on similar ideas and a generalisation of the integration lemma in \cite{Mirzakhani}. We first remark that, if we replace $(A,B,C,D)$ with their absolute value, it is easy to prove by induction on $2g - 2 + n > 0$ that, under our assumptions, the formulas \eqref{theteatthtuhtu} produce well-defined, finite functions that are uniformly bounded by a polynomial in $L_1,\ldots,L_n$. This comes from the fact that the measure of integration is $\ell\,\dd \ell$, making a function bounded by $\ell^{-\eta}$ integrable near $0$ provided $\eta < 2$. Near $\infty$ the integrability comes from the decay faster than any power $\ell^{-s}$ for fixed $L_i$s. By Tonelli's theorem, in order to prove the Theorem in general, it is sufficient to justify that the integrals of the GR amplitudes for the initial data $(|A|,|B|,|C|,|D|)$ satisfy this recursion -- without caring in the intermediate computations whether the integrals are finite of equal to $+\infty$.

Consider $\Sigma$ of type $(g,n)$ with $2g - 2 + n \geq 2$, $L \in \mathbb{R}_{+}^n$ and let $\gamma$ be a multicurve in with ordered components $\gamma_1,\ldots,\gamma_k$ in $\Sigma$. We consider the orbifold
\[
\mathcal{M}_{\Sigma}^{\gamma}(L) = \mathcal{T}_{\Sigma}(L)\Big/\bigcap_{i = 1}^k {\rm Stab}(\gamma_i),
\]
where ${\rm Stab}(\gamma_i)$ is the stabiliser of $\gamma_i$ in $\Gamma_{\Sigma}^{\partial}$. It is equipped with the Weil--Petersson symplectic structure, and with an $\prod_{i = 1}^k (\mathbb{R}/2^{-\mathfrak{t}_i}\mathbb{Z})$ action induced by the moment map $\mathcal{L}^{\gamma}(\sigma) = (\ell_{\sigma}^2(\gamma_i)/2)_{i = 1}^k$, where
\[
\mathfrak{t}_i = \left\{\begin{array}{lll} 1 & & {\rm if}\,\,\gamma_i\,\,{\rm separates}\,\,{\rm off}\,\,{\rm a}\,\,{\rm torus}\,\,{\rm with}\,\,{\rm one}\,\,{\rm boundary}\,\,{\rm component} \\
 0 & & {\rm otherwise} \end{array}\right.\,.
\]
The presence of $2^{-\mathfrak{t}_i}$ comes from the fact that the half-Dehn twist belongs to ${\rm Stab}(\gamma_i)$ when $\gamma_i$ separates off a torus with one boundary component.

If we include $\gamma$ into a pair of pants decomposition and consider the associated Fenchel--Nielsen coordinates $(l_j,\tau_j)_j$, the Weil--Petersson symplectic form on $\mathcal{T}_{\Sigma}$ takes the form
\[
\omega_{{\rm WP}} = \sum_{j} \dd l_j \wedge \dd \tau_j.
\]
With this description, one can check that for each $\ell \in \mathbb{R}_{+}^k$, the symplectic quotient $(\mathcal{L}^{\gamma})^{-1}(\ell)/\!/U(1)^{k}$ is a symplectic orbifold which is symplectomorphic to a moduli space $\tilde{\mathcal{M}}_{\Sigma^{\gamma}}(L,\ell,\ell)$, which is the cartesian product of moduli spaces of bordered Riemann surfaces associated with each connected component of $\Sigma^{\gamma}$, where boundary components that were originally in $\Sigma$ have fixed length $L$, boundary components that come from a $\gamma_i$ have length $\ell_i$, and each factor associated with a component of $\Sigma^{\gamma}$ which is a torus with one boundary component $\gamma_i$ should be the double cover $\tilde{\mathcal{M}}_{1,1}(\ell_i)$ instead of the usual moduli space $\mathcal{M}_{1,1}(\ell_i) = \tilde{\mathcal{M}}_{1,1}(\ell_i)/\langle {\rm elliptic}\,\,{\rm involution}\}$.

We deduce that if $f$ is a nonnegative measurable function on $\mathcal{M}_{\Sigma}^{\gamma}$ which is invariant under the torus action, it induces a measurable function $\tilde{f}$ on $\mathcal{M}_{\Sigma^{\gamma}}$ for each $\ell \in \mathbb{R}_{+}^k$ and we have for any $\ell \in \mathbb{R}_{+}^k$
\begin{equation}
\label{intsympequ}\int_{\mathcal{M}_{\Sigma}^{\gamma}} f\,\dd \mu_{{\rm WP}} = \int_{\mathbb{R}_{+}^k} \bigg( \int_{\mathcal{M}_{\Sigma^{\gamma}}(L,\ell,\ell)} \tilde{f}\,\dd\mu_{{\rm WP}} \bigg) \prod_{i = 1}^k \ell_i\, \dd \ell_i.
\end{equation}
We use that the fibers of the moment map have volume $\prod_{i = 1}^k 2^{-\mathfrak{t}_i}\ell_i$ but the factors of $2$ disappear after we replaced the symplectic quotient spaces $\tilde{\mathcal{M}}_{\Sigma^{\gamma}}(L,\ell,\ell)$ with the usual moduli space of bordered Riemann surfaces $\mathcal{M}_{\Sigma^{\gamma}}(L,\ell,\ell)$.

We apply this discussion to the GR amplitudes associated with $(|A|,|B|,|C|,|D|)$, denoted again $\Omega_{\Sigma}$ for convenience of the proof. In each $\Gamma_{\Sigma}^{\partial}$-orbit $\mathcal{O} \subset \mathcal{P}_{\Sigma}$, let us choose arbitrarily a representative $[P_{\mathcal{O}}]$ and denote $\gamma_{\mathcal{O}}$ the multicurve with ordered components $\partial P_{\mathcal{O}} \cap \Sigma^{\circ}$. We have
\begin{equation}
\label{gnigfgfggfsfdgun}
\begin{split}
& \quad \int_{\mathcal{M}_{\Sigma}(L)} \bigg(\sum_{[P] \in \mathcal{O}} X_{\mathcal{O}}(\vec{\ell}_{\sigma}(\partial P))\,\Omega_{\Sigma - P}(\sigma|_{\Sigma - P})\bigg)\dd \mu_{{\rm WP}}(\sigma) \\
& = \int_{\mathcal{M}_{\Sigma}^{\gamma_{\mathcal{O}}(L)}} X_{P}\big(\vec{\ell}_{\sigma}(\partial P_{\mathcal{O}})\big)\,\Omega_{\Sigma - P_{\mathcal{O}}}\big(\sigma|_{\Sigma - P_{\mathcal{O}}}\big)\dd\mu_{{\rm WP}}(\sigma)  ,
\end{split}
\end{equation} 
where $X_{\mathcal{O}}$ is the function $B$ or $C$ depending on the type of $\mathcal{O}$. The integral of the GR amplitude is obtained by summing the result of integration over all $\Gamma_{\Sigma}^{\partial}$-orbits $\mathcal{O} \subseteq \mathcal{P}_{\Sigma}$. Since the function in the right-hand side of \eqref{gnigfgfggfsfdgun} does not depend on the twist along the components of $\gamma_{\mathcal{O}}$, we can apply \eqref{intsympequ}. For the orbit $\mathcal{O} = \mathcal{P}_{\Sigma}^{b}$ with $b \in \pi_0(\partial_+\Sigma) \simeq \{2,\ldots,n\}$, we have $X_{\mathcal{O}} = B$ and $\Sigma - P_{\mathcal{O}}$ has type $(g,n - 1)$, hence
\begin{equation*}
\begin{split}
& \quad \int_{\mathcal{M}_{\Sigma}(L)} \bigg(\sum_{[P] \in \mathcal{O}} X_{\mathcal{O}}(\vec{\ell}_{\sigma}(\partial P))\,\Omega_{\Sigma - P}(\sigma|_{\Sigma - P})\bigg)\dd \mu_{{\rm WP}}(\sigma) \\
& =  \int_{\mathbb{R}_{+}} B(L_1,L_b,\ell)\bigg(\int_{\mathcal{M}_{g, n -1}(\ell,L_2,\ldots,\widehat{L_m},\ldots,L_n)} \Omega_{g,n - 1}\,\dd\mu_{{\rm WP}}\bigg) \ell \, \dd \ell\,.
\end{split}
\end{equation*}
If $\mathcal{O} \subseteq \mathcal{P}_{\Sigma}^{\emptyset}$, $\gamma_{\mathcal{O}}$ has two components, $X_{\mathcal{O}} = \frac{1}{2}C$ and we rather obtain
\begin{equation*}
\begin{split}
& \quad  \int_{\mathcal{M}_{\Sigma}(L)} \bigg(\sum_{[P] \in \mathcal{O}} X_{\mathcal{O}}(\vec{\ell}_{\sigma}(\partial P))\,\Omega_{\Sigma - P}(\sigma|_{\Sigma - P})\bigg)\dd \mu_{{\rm WP}}(\sigma) \\
& = \frac{1}{2} \int_{\mathbb{R}_{+}^2}C(L_1,\ell,\ell')\bigg( \int_{\mathcal{M}_{\Sigma - P_{\mathcal{O}}(\ell,\ell',L_2,\ldots,L_n)}} \Omega_{\Sigma - P_{\mathcal{O}}}\dd\mu_{{\rm WP}}\bigg) \ell\ell' \, \dd \ell \dd \ell'\,
\end{split}
\end{equation*}
with the appropriate labelling of boundary components of $\Sigma - P_{\mathcal{O}}$. Summing over $\mathcal{O}$ reconstructs all the terms in \eqref{theteatthtuhtu} and justifies the claim.
\hfill $\Box$

\section{Revisiting Mirzakhani-McShane identities}

\label{MMId}

In this section we demonstrate that how fundamental contributions of Mirzakhani and Kontsevich fit into the framework of the geometric recursion. We work here with the target theory $E(\Sigma) = \mathscr{C}^0(\mathcal{T}_{\Sigma})$.

\subsection{Mirzakhani identities}

Let us consider the following initial data
\begin{equation}
\label{iniMirza}
\begin{split}
A^{{\rm M}}(L_1,L_2,L_3) & =  1, \\
B^{{\rm M}}(L_1,L_2,\ell) & = 1 - \frac{1}{L_1}\ln \left(\frac{\cosh\left(\frac{L_2}{2}\right) + \cosh\left(\frac{L_1 + \ell}{2}\right)}{{\cosh\left(\frac{L_2}{2}\right) + \cosh\left(\frac{L_1 - \ell}{2}\right)}}\right), \\ 
C^{{\rm M}}(L_1,\ell,\ell') & = \frac{2}{L_1} \ln\left(\frac{\exp(\frac{L_1}{2}) + \exp(\frac{\ell+\ell'}{2})}{\exp(-\frac{L_1}{2}) + \exp(\frac{\ell + \ell'}{2})}\right) , \\
D^{{\rm M}}_{T}(\sigma) & = \sum_{\gamma \in S_{T}^{\circ}} C^{{\rm M}}\big(\ell_{\sigma}(\partial T),\ell_{\sigma}(\gamma),\ell_{\sigma}(\gamma)\big).
\end{split}
\end{equation}
The functions $B^{{\rm M}}$ and $C^{{\rm M}}$ can equivalently be expressed in terms of $F(x) = 2\ln(1 + e^{x/2})$ as
\begin{equation}
\begin{split}
B^{{\rm M}}(L_1,L_2,\ell) & = \frac{1}{2L_1}\big(F(L_1 +L_2 - \ell) + F(L_1 - L_2 -\ell) - F(-L_1 + L_2 - \ell) - F(-L_1 - L_2 - \ell)\big), \\
C^{{\rm M}}(L_1,\ell,\ell') & = \frac{1}{L_1}\big(F(L_1 - \ell - \ell') - F(-L_1 -\ell - \ell')\big).
\end{split}
\end{equation}
From these formulae we immediately conclude that $(A,B,C,D)$ is strongly admissible. Mirzakhani's identities can now be reformulated as
\begin{theorem}\cite{Mirzakhani}
\label{thMirza} For any object $\Sigma$ in $\Surf_1$, $\Omega_{\Sigma}$ is the constant function $1$ on $\mathcal{T}_{\Sigma}$. 
\end{theorem}

This shows that Mirzakhani's identities should be seen as a recursion producing the constant function $1$ on the Teichm\"uller space. The way Theorem~\ref{thMirza} is proved in \cite{Mirzakhani} does not require showing a priori the convergence of the series \eqref{GRseriesagain} defining $\Omega^{{\rm M}}$. Nevertheless, if one takes for granted Theorem~\ref{thMirza} for the case of the torus with one boundary i.e. $D_{T}(\sigma) = 1$, Theorem~\ref{welldefGR} justifies the absolute convergence of these series uniformly on any compact.

In \cite{Mirzakhani}, Mirzakhani used her identities to obtain a recursive formula for the Weil--Petersson volumes of the moduli space of bordered surfaces
\[
V\Omega^{{\rm M}}_{g,n}(L) = \int_{\mathcal{M}_{g,n}(L)} 1\cdot \dd\mu_{{\rm WP}},
\]
and this is a special case of our Theorem~\ref{intTRHTT}. The initial data for this recursion is
\[
V\Omega^{{\rm M}}_{0,3}(L_1,L_2,L_3) = 1,\qquad V\!D^{{\rm M}}(L_1) = \frac{\pi^2}{12} + \frac{L_1^2}{48}.
\]

\begin{theorem} \cite{EOwp} 
\label{iuniuufgbdfginugiun}
This recursion for the Weil--Petersson volumes $V\Omega^{{\rm M}}_{g,n}$ is equivalent to the statement that 
\[
\omega_{g,n}^{{\rm M}}(z_1,\ldots,z_n) =  \bigg(\int_{\mathbb{R}_{+}^{n}}  \prod_{i = 1}^n \dd L_i\,L_i\,e^{-z_iL_i}\,V\Omega^{{\rm M}}_{g,n}(L_1,\ldots,L_n)\bigg) \,\dd z_1 \otimes \cdots \otimes \dd z_n
\]
is computed by Eynard--Orantin's topological recursion for the initial data
\begin{equation}
x(z) = z^2/2,\qquad y(z) = -\frac{\sin(2\pi z)}{2\pi},\qquad \omega_{0,2}(z_1,z_2) = \frac{\dd z_1 \otimes \dd z_2}{(z_1 - z_2)^2}.
\end{equation}
\end{theorem}

\begin{remark}  It is an interesting problem to find a mechanism which would explain \textit{a priori} -- without relying on Mirzakhani's identity/Theorem~\ref{thMirza} -- that $\Omega^{{\rm M}}_{\Sigma}$ defined by the inductive formula \eqref{GRseriesagain} is invariant under braiding of all boundary components of $\Sigma$.  We do not know if Mirzakhani initial data satisfy the four symmetry constraints of Definition~\ref{DEFGRSSSS} -- if they hold true, it would give such an explanation via our Theorem~\ref{SGRSYMTH}. On the other hand, one can check -- see e.g. \cite{TRABCD,TRlecture} -- that Mirzakhani's initial data \eqref{iniMirza} satisfy the ``averaged'' constraints \eqref{theuytbsugfngoig}, and therefore that the topological recursion for $V\Omega^{{\rm M}}_{g,n}(L_1,\ldots,L_n)$ does produce \textit{a priori} symmetric functions of its $n$ length variables.
\end{remark}

Let us review the relation between the Weil--Petersson volumes $V\Omega^{{\rm M}}_{g,n}(L)$ and intersection theory on Deligne--Mumford compactification of the moduli space $\overline{\mathfrak{M}}_{g,n}$ of genus $g$ Riemann surfaces $S$ with $n$ labelled punctures $p_1,\ldots,p_n \in S$. We recall that $\overline{\mathfrak{M}}_{g,n}$ is a complex orbifold for which Poincar\'e duality holds. Let us denote $\psi_i$ the first Chern class of the cotangent line bundle $T^*_{p_i}S$ at the $i$-th puncture. We also denote $\kappa_d$ the class of complex degree $d$ obtained by pushforward of $\psi_{n + 1}^{d + 1}$ via the morphism $\overline{\mathfrak{M}}_{g,n + 1} \rightarrow \overline{\mathfrak{M}}_{g,n}$ forgetting the last puncture. It is well-known that the cohomology class of the Weil--Petersson symplectic form on $\mathfrak{M}_{g,n}$ is $2\pi^2 \kappa_{1}$ \cite{Wolperthomology}. By examination of the symplectic reduction on the space $\mathfrak{M}_{g,n}$ with the moment map $(L_i^2/2)_{i = 1}^n$, Mirzakhani proved  
\begin{theorem}\cite{Mirzaint} 
\label{IntPsi} For $2g - 2 + n > 0$, we have
\[
V\Omega^{{\rm M}}_{g,n}(L_1,\ldots,L_n) = \int_{\overline{\mathfrak{M}}_{g,n}} \exp\bigg(2\pi^2 \kappa_1 + \sum_{i = 1}^n \frac{L_i^2}{2}\,\psi_i\bigg),
\]
where the exponential is understood by expanding in Taylor series, and keeping only the terms of cohomology degree $2d_{g,n} = 6g - 6 + 2n$. 
\end{theorem}

\subsection{Interpolation to Kontsevich's amplitudes}

Let us consider a deformation of Mirzakhani initial data, which consists in rescaling all length variables by a factor $\beta > 0$, namely
\begin{equation}
X^{\beta}(L_1,L_2,L_3) = X^{{\rm M}}(\beta L_1,\beta L_2,\beta L_3),\qquad X \in \{A,B,C\}.
\end{equation}
 One easily checks that the initial data $(A^{\beta},B^{\beta},C^{\beta})$ is strongly admissible with constants that can be chosen independently of $\beta \geq 1$. We induce $D^{\beta}$ from $C^{\beta}$, that is
\[
 D^{\beta}_{T}(\sigma) = \sum_{\gamma \in S_{T}^{\circ}} C^{{\rm M}}\big(\beta \ell_{\sigma}(\partial T),\beta \ell_{\sigma}(\gamma),\beta\ell_{\sigma}(\gamma)\big).
\]
This deformation of the Mirzakhani initial data is not induced by a flow on the Teichm\"uller spaces: except for $\Sigma = P$, there is no flow on $\mathcal{T}_{\Sigma}$ whose effect is to rescale the length of all curves by $\beta$. So, $D^{\beta}$ is \textit{a priori} not equal to the constant function $1$ on $\mathcal{T}_{T}$ as soon as $\beta \neq 1$. We denote $\Omega^{\beta}_{\Sigma}$ the corresponding GR amplitudes.

Notice that $\lim_{\beta \rightarrow \infty} \beta^{-1}F(\beta x) = [x]_{+}$. and the convergence is uniform for $x \in \mathbb{R}$. Therefore, $(A_{\beta},B_{\beta},C_{\beta})$ converges uniformly on  $[\epsilon,+\infty)^3$ for each $\epsilon > 0$, to a strongly admissible initial data which we denote $(A^{{\rm K}},B^{{\rm K}},C^{{\rm K}})$. Explicitly
\begin{equation}
\label{Kontinfigufn}\begin{split}
A^{{\rm K}}(L_1,L_2,L_3) & = 1, \\
B^{{\rm K}}(L_1,L_2,\ell) & = \frac{1}{2L_1}\big([L_1 + L_2 - \ell]_{+} + [L_1 - L_2 - \ell]_{+}  - [-L_1 + L_2 - \ell]_{+} \big), \\
C^{{\rm K}}(L_1,\ell,\ell') & = \frac{1}{L_1}\,[L_1  - \ell - \ell']_{+}.
\end{split}
\end{equation}
We call it the Kontsevich initial data and we denote $\Omega^{{\rm K}}$ the corresponding GR amplitudes. 

\begin{lemma}
For any object $\Sigma$ in $\Surf_1$, we have
\begin{equation}
\label{convtKMKM} \forall \sigma \in \mathcal{T}_{\Sigma},\qquad \lim_{\beta \rightarrow \infty} \Omega^{\beta}(\sigma) = \Omega^{{\rm K}}(\sigma),
\end{equation}
and the convergence is uniform for $\sigma$ in any compact of $\mathcal{T}_{\Sigma}$. Besides, we have for $2g - 2 + n > 0$
\begin{equation}
\label{cinfguyfgn}V\Omega^{\beta}_{g,n}(L_1,\ldots,L_n) = \beta^{-(6g - 6 + 2n)}\,V\Omega^{{\rm M}}_{g,n}(\beta L_1,\ldots,\beta L_n),
\end{equation}
and
\[
\lim_{\beta \rightarrow \infty} V\Omega^{\beta}_{g,n}(L_1,\ldots,L_n) = V\Omega^{{\rm K}}_{g,n}(L_1,\ldots,L_n) = \int_{\overline{\mathfrak{M}}_{g,n}} \exp\bigg(\sum_{i = 1}^n \frac{L_i^2}{2}\,\psi_i\bigg).
\]
In particular $V\!D^{{\rm K}}(L_1) = \frac{L_1^2}{48}$.
\end{lemma}
\begin{remark} The geometric realisation of this deformation via a flow on Teichm\"uller space and the geometric meaning of the functions $\Omega^{{\rm K}}$ are explained in \cite{WKarticle}.
\end{remark}
\noindent \textbf{Proof.} Since the constants appearing in the admissibility bounds \eqref{theinininini} for $(A^{\beta},B^{\beta},C^{\beta})$ can be chosen independently of $\beta \in [1,+\infty]$, the proof of Theorem~\ref{welldefGR} actually shows that the GR series defining $\Omega^{\beta}$ are uniformly convergent for $\sigma$ in any compact of $\mathcal{T}_{\Sigma}$ and $\beta \in [1,+\infty]$. Therefore, we can take the limit $\beta \rightarrow \infty$ term by term in the series, which yields \eqref{convtKMKM}. Likewise, since the constants appearing in the strong admissibility bounds \eqref{ABCstrongs} can be chosen independently of $\beta \in [1,+\infty]$, the proof of Theorem~\ref{intTRHTT} shows that we can take the term by term $\beta \rightarrow \infty$ limit in the topological recursion \eqref{theteatthtuhtu}, which yields
\[
\lim_{\beta \rightarrow \infty} V\Omega^{\beta}_{g,n}(L_1,\ldots,L_n) = V\Omega^{{\rm K}}_{g,n}(L_1,\ldots,L_n).
\]

We prove the equality \eqref{cinfguyfgn} by induction. The case of  $(g,n) = (0,3)$ is clear since $A^{\beta} = A^{{\rm K}} = 1$. For $(g,n) = (1,1)$, we compute for $L_1 > 0$
\begin{equation}
\begin{split}
VD^{\beta}(L_1) & = \int_{\mathcal{M}_{T}(L_1)}\bigg( \sum_{\gamma \in S_{T}^{\circ}} C^{{\rm M}}\big(\beta L_1,\beta \ell_{\sigma}(\gamma),\beta \ell_{\sigma}(\gamma)\big)\bigg)\dd\mu_{{\rm WP}} = \frac{1}{2} \int_{\mathbb{R}_{+}} \, C^{{\rm M}}(\beta L_1,\beta \ell,\beta \ell)\ell  \, \dd \ell\,\\
& = \frac{1}{2\beta^2} \int_{\mathbb{R}_{+}} \, C^{{\rm M}}(\beta L_1,\tilde{\ell},\tilde{\ell})\tilde{\ell}  \, \dd \tilde{\ell}\, = \frac{V\!D^{{\rm M}}(L_1)}{\beta^2},
\end{split}
\end{equation}
where to get the second line we used the same tricks as in the proof of Theorem~\ref{intTRHTT}, and the $1/2$ takes into account the fact that $\Gamma_{T}^{\partial}$ contains the half-Dehn twist along $\gamma \in S_{T}^{\circ}$. Assume the claim holds for $(g',n')$ such that $2g' - 2 + n' < 2g - 2 + n$. We examine the terms in the topological recursion formula \eqref{theteatthtuhtu} for $V\Omega^{\beta}_{g,n}$. It contains for instance the terms
\begin{equation}
\begin{split}
& \quad \sum_{m = 2}^n \int_{\mathbb{R}_{+}} B(\beta L_1,\beta L_m,\beta \ell)\,V\Omega^{\beta}_{g,n - 1}(\ell,L_2,\ldots,\widehat{L_m},\ldots,L_n) \ell\,\dd \ell \\
& =  \sum_{m = 2}^n \beta^{-(6g - 6 + 2(n - 1))} \cdot \beta^{-2} \int_{\mathbb{R}_{+}}\,B(\beta L_1,\beta L_m,\tilde{\ell})\,V\Omega^{{\rm M}}_{g,n - 1}(\beta \tilde{\ell},\beta L_2,\ldots,\widehat{L_m},\ldots,\beta L_n)  \tilde{\ell} \dd\tilde{ \ell}\,,
\end{split}
\end{equation}
which is $\beta^{-(6g - 6 + 2n)}$ times the corresponding terms in the recursion for $V\Omega^{{\rm M}}(\beta L_1,\ldots,\beta L_n)$. The $C$-terms are handled in the same way, and this completes the proof of \eqref{cinfguyfgn}.

Finally, the equality
\[
\lim_{\beta \rightarrow \infty} \beta^{-(6g - 6 + 2n)}\,V\Omega^{{\rm M}}_{g,n}(\beta L_1,\ldots,\beta L_n) = \int_{\overline{\mathfrak{M}}_{g,n}} \exp\bigg(\sum_{i = 1}^{n} \frac{L_i^2}{2}\,\psi_i\bigg)
\]
is immediate once we observe in Theorem~\ref{IntPsi} that $V\Omega^{{\rm M}}_{g,n}(\beta L_1,\ldots,\beta L_n)$ is a polynomial of degree $6g - 6 + 2n$ in $\beta$: taking the limit $\beta \rightarrow \infty$ extracts the top degree coefficient and yields the right-hand side. \hfill $\Box$

\vspace{0.2cm}

The Laplace transform of the topological recursion \eqref{theteatthtuhtu} for the Kontsevich initial data \eqref{Kontinfigufn} was first computed in \cite{Safnukwp} and shows that Eynard--Orantin's topological recursion for the spectral curve
\[
x(z) = z^2/2,\qquad y(z) = -z,\qquad \omega_{0,2}(z_1,z_2) = \frac{\dd z_1\dd z_2}{(z_1 - z_2)^2}.
\]
is computing
\begin{equation}
\begin{split}
\omega_{g,n}^{{\rm K}}(z_1,\ldots,z_n) & = \bigg(\int_{\mathbb{R}_{+}^n} e^{-z_iL_i} V\Omega^{{\rm K}}_{g,n}(L_1,\ldots,L_n)  \prod_{i = 1}^n \,L_i \,\dd L_i\bigg)\dd z_1 \otimes \cdots \otimes \dd z_n \\
& = \sum_{\substack{m_1,\ldots,m_n \geq 0 \\ m_1 + \cdots + m_n = 3g - 3 + n}} \bigg(\int_{\overline{\mathfrak{M}}_{g,n}} \prod_{i = 1}^n \psi_i^{m_i}\bigg) \bigotimes_{i = 1}^n \frac{(2m_i + 1)!!\dd z_i}{z_i^{2m_i + 2}} .
\end{split}
\end{equation}
This is one of the many ways to prove the last equality: This statement is equivalent to the Virasoro constraints conjectured by Witten in \cite{Witten} and proved by Kontsevich in \cite{Kontsevich} and it has a long history into which we do not enter here. See however \cite{WKarticle} for a geometric proof using the ideas of GR.

\section{Statistics of lengths of multicurves}
\label{S7n7n7}
\subsection{Generalisation of Mirzakhani--McShane identities}
\label{F9191}

\begin{theorem}
\label{Genegin}Let $(A,B,C,D)$ be an admissible initial data and assume that the corresponding GR amplitude $\Omega_{\Sigma}$ is invariant under braiding of all boundary components of $\Sigma$. Let $f \in \mathscr{F}(\mathbb{R}_{+})$ be such that for any $\epsilon > 0$ and $s \geq 0$ there exists $F_{s,\epsilon} > 0$ such that
\begin{equation}
\label{sifgnfdgi} \sup_{\ell \geq \epsilon} |f(\ell)|\,\ell^{s} \leq F_{s,\epsilon}.
\end{equation}
Then, for any object $\Sigma$ in $\Surf_1$, the series
\begin{equation}
\label{thedefOmgifsugn}\Omega_{\Sigma}[f](\sigma) := \sum_{c \in M'_{\Sigma}} \Omega_{\Sigma^{c}}(\sigma|_{\Sigma^{c}})\,\prod_{\gamma \in \pi_0(c)} f(\ell_{\sigma}(\gamma))
\end{equation}
converges absolutely and uniformly for $\sigma$ in any compact of $\mathcal{T}_{\Sigma}$, and it coincides with the GR amplitude associated with the twisted initial data
\begin{equation}
\label{twsitini}
\begin{split}
A[f](L_1,L_2,L_3) & = A(L_1,L_2,L_3), \\
B[f](L_1,L_2,\ell) & = B(L_1,L_2,\ell) + A(L_1,L_2,\ell)f(\ell), \\
C[f](L_1,\ell,\ell') & = C(L_1,\ell,\ell') + B(L_1,\ell,\ell')\,f(\ell) + B(L_1,\ell',\ell)\,f(\ell') + A(L_1,\ell,\ell')\,f(\ell)f(\ell'), \\
D_{T}[f](\sigma) & = D_{T}(\sigma) + \sum_{\gamma \in S_{T}^{\circ}} A(L_1,\ell_{\sigma}(\gamma),\ell_{\sigma}(\gamma))\,f(\ell_{\sigma}(\gamma)).
\end{split}
\end{equation}
\end{theorem}

For a surface of type $(g,n)$, primitive multicurves have at most $3g -3 + n$ components, so $\Omega_{\Sigma}[f]$ depends polynomially on $f$. The term of degree $0$ in $f$ corresponds to the empty multicurve belongs to $M_{\Sigma}'$ and is equal to $\Omega_{\Sigma}$. We recall that $\Sigma^{c}$ is the surface $\Sigma$ cut along $c$. Although it is not an object in $\Surf_1$, choosing arbitrarily a distinguished boundary component in each component of $\Sigma$ turn it into an object of $\Surf_1$ so that $\Omega_{\Sigma^{c}}$ makes sense independently on the choice made by the assumption of invariance under braidings of all boundary components. As a consequence, $\Omega_{\Sigma}[f]$ is also invariant under braidings of all boundary components. We call $(A[f],B[f],C[f],D[f])$ the twisting of $(A,B,C,D)$ by $f$. For instance, the twisting of Mirzakhani initial data \eqref{iniMirza} gives access to statistics of the length of multicurves
\[
\Omega^{{\rm M}}[f](\sigma) = \sum_{c \in M'_{\Sigma}} \prod_{\gamma \in \pi_0(c)} f(\ell_{\sigma}(\gamma)).
\]
The specialisation of the proof below to this case amounts to utilising Mirzakhani--McShane identities as a partition of unity to deduce the recursion for $\Omega^{{\rm M}}$. We can also see Theorem~\ref{Genegin} as a generalisation of Mirzakhani--McShane identities, in the sense that the result of GR is identified with an independently defined function on $\mathcal{T}_{\Sigma}$.

\vspace{0.2cm}

\noindent \textbf{Proof.} One easily checks that the admissibility of $(A,B,C,D)$ implies the admissibility of $(A[f],B[f],C[f],D[f])$ given the assumption \eqref{sifgnfdgi}. It is enough to prove the result for connected surfaces, and it holds for pairs of pants and tori with one boundary component, since by definition $\Omega_{P}[f]$ is given by $A = A[f]$ and we defined $D_{T}[f]$ to be equal to $\Omega_{T}[f]$. We now assume that $\Sigma$ has Euler characteristic $\leq -2$ and fix $\sigma \in \mathcal{T}_{\Sigma}$. Since all sums we consider are absolutely convergent, we can apply Fubini's theorem to cluster or exchange summations. 

\begin{figure}[ht!]
\begin{center}
\includegraphics[width=0.95\textwidth]{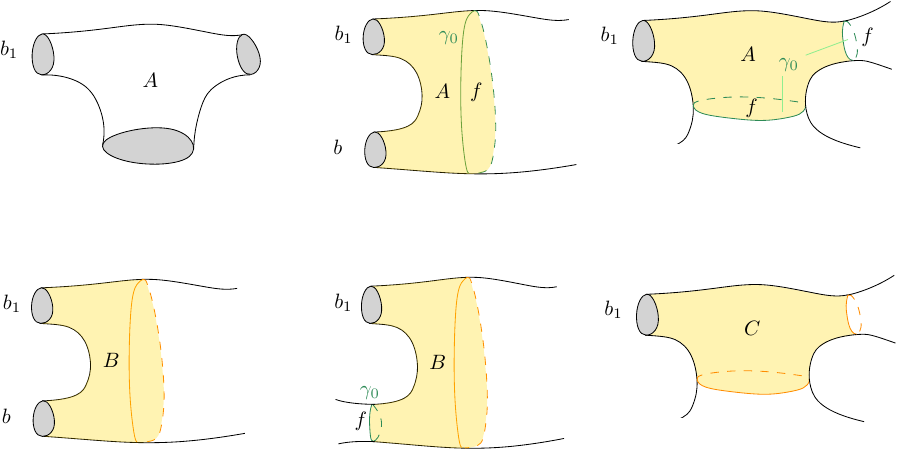}
\caption{\label{Statssss} The terms in the twisted initial data correspond to either to $\Sigma = P$ (first picture) or to the possible configurations of $[P] \in \mathcal{P}_{\Sigma^c}$ (in yellow) with respect to the multicurve $c$. We have indicated $\gamma_0 = \partial P \cap c$ in green.}
\end{center}
\end{figure}

If $c \in M'_{\Sigma}$, we have by multiplicativity
\begin{equation}
\label{oungdifggggun}\Omega_{\Sigma^{c}}(\sigma) = \prod_{S \in \pi_0(\Sigma^{c})} \Omega_{S}(\sigma|_{S}),
\end{equation}
We denote by $\Sigma^{c}_-$ the component of $\Sigma^{c}$ containing $\partial_-\Sigma$, and recall that \eqref{oungdifggggun} is valid for any choice of distinguished boundary components on each $S \neq \Sigma^{c}_-$. Since $\chi_{\Sigma} \leq -2$, $\Sigma^{c}_-$ cannot be a torus with one boundary component. The contribution to \eqref{thedefOmgifsugn} from primitive multicurves $c$ such that $\Sigma^{c}_-$ is a pair of pants is equal to 
\begin{equation}
\label{tirngifugngg}
\sum_{b \in \pi_0(\partial_+\Sigma)} \sum_{[P] \in \mathcal{P}_{\Sigma}^{b}} A(l_{P}^1,l_P^2,l_P^3)f(l_P^3)\,\Omega_{\Sigma - P}[f](\sigma|_{\Sigma - P})  + \frac{1}{2}\sum_{[P] \in \mathcal{P}_{\Sigma}^{\emptyset}} A(l_P^1,l_P^2,l_P^3)f(l_P^2)f(l_P^3) \Omega_{\Sigma - P}[f](\sigma|_{\Sigma - P}).
\end{equation}
Here we use the shorter notation $l_{P}^i = \ell_{\sigma}(\partial_iP)$. Recall that $\partial_1 P = \partial_-\Sigma$ in any case and $\partial_2 P = b$ when $[P] \in \mathcal{P}_{\Sigma}^{b}$ (see Section~\ref{glueinpant}) so that we can also denote $l^b = l_P^2$. The $1/2$ in the last sum takes into account the fact that boundary components of $[P] \in \mathcal{P}_{\Sigma}^{\emptyset}$ are ordered while components of $c$ are not.

For each $c \in M'_{\Sigma}$ such that $\Sigma^{c}_-$ is not a pair of pants, we insert the GR formula
\[
\Omega_{\Sigma^{c}_-}(\sigma) = \sum_{[P] \in \mathcal{P}_{\Sigma^{c}_-}} X_{P}(\vec{\ell}_{\sigma}(\partial P))\,\Omega_{\Sigma^{c}_- - P}(\sigma|_{\Sigma^{c}_- - P}),
\]
where $X_P$ is equal to $B$ or $C$ depending on the type of $[P]$, inside \eqref{oungdifggggun} and then \eqref{thedefOmgifsugn}. We exchange the order of the resulting summations over $c$ and $[P]$, and examine for each $[P] \in \mathcal{P}_{\Sigma}$ which $c$ can contribute. For a given $[P] \in \mathcal{P}_{\Sigma}^{b}$ for some $b \in \pi_0(\partial_+\Sigma)$, we only receive the contributions of $c \in M'_{\Sigma - P}$ and thus a total contribution
\begin{equation}
\label{OHY1}\sum_{b \in \pi_0(\partial_+\Sigma)} \sum_{[P] \in \mathcal{P}_{\Sigma}^{b}} B(\vec{\ell}_{\sigma}(\partial P)) \Omega_{\Sigma - P}[f](\sigma|_{\Sigma - P}).
\end{equation}
For a given $[P] \in \mathcal{P}_{\Sigma}^{\emptyset}$, there are three types of $c$ contributing in \eqref{thedefOmgifsugn}. When $\partial_2 P$ is a component $\gamma_0$ of $c$ but not $\partial_3 P$, we get a total contribution
\begin{equation}
\label{OHY2}\frac{1}{2} \sum_{[P] \in \mathcal{P}_{\Sigma}^{\emptyset}} B(\ell_{\sigma}(l_P^1,l_P^2,l_P^3)f(l_P^2)\,\Omega_{\Sigma - P}[f](\sigma|_{\Sigma - P}).
\end{equation}
by summing over the multicurve $c \setminus \gamma_0 \in M'_{\Sigma - P}$. From the situations where $\partial_3P$ is a component of $c$ but not $\partial_2P$, we get likewise a total contribution
\begin{equation}
\label{OHY3}\frac{1}{2} \sum_{[P] \in \mathcal{P}_{\Sigma}^P{\emptyset}} B(l_P^1,l_P^3,l_P^2)f(l_P^3)\,\Omega_{\Sigma - P}[f](\sigma|_{\Sigma - P}).
\end{equation}
From the situations where $\partial_2P$ and $\partial_3P$ are not components of $c$, we get a total contribution
\begin{equation}
\label{OHY4}\frac{1}{2} \sum_{[P] \in \mathcal{P}_{\Sigma}^{\emptyset}} C(l_P^1,l_P^2,l_P^3)\,\Omega_{\Sigma - P}[f](\sigma|_{\Sigma - P}).
\end{equation}
by summing over $c \in M'_{\Sigma - P}$. Adding up \eqref{tirngifugngg} and \eqref{OHY1}--\eqref{OHY4}, we recognise the definition of the GR amplitudes for the twisted initial data. \hfill $\Box$

\subsection{Integration over the moduli space}

We keep the notations of Section~\ref{Sinrgng} and Theorem~\ref{Genegin}.

\begin{corollary}
\label{costatunfg}Let $(A,B,C,D)$ be measurable, strongly admissible initial data and $f \colon \mathbb{R}_{+} \rightarrow \mathbb{R}$ be a measurable function for which there exists $\eta \in [0,2)$ and for any $s \geq 0$ there exists $M_s > 0$ such that for any $\ell > 0$
\begin{equation}
\label{storngagum}\big|f(\ell)\big| \leq \frac{M_s}{\ell^{\eta}(1 + \ell)^{s}}.
\end{equation}
Then, the twisted initial data \eqref{twsitini} is strongly admissible, and the twisted GR amplitudes induce for any $g \geq 0$ and $n \geq 1$ such that $2g - 2 + n > 0$ and any $L \in \mathbb{R}_{+}^n$ functions $\Omega_{g,n}[f]$ on $\mathcal{M}_{g,n}(L)$ which are integrable with respect to $\mu_{{\rm WP}}$. The integrals
\[
V\Omega_{g,n}[f](L) = \int_{\mathcal{M}_{g,n}(L)} \Omega_{g,n}[f]\,\dd\mu_{{\rm WP}}
\]
satisfy the topological recursion \eqref{theteatthtuhtu} for twisted initial data, with base cases
\begin{equation*}
\begin{split}
V\Omega_{0,3}[f](L_1,L_2,L_3) & = A(L_1,L_2,L_3), \\
V\Omega_{1,1}[f](L_1) := V\!D[f](L_1) & =  V\!D(L_1) + \frac{1}{2} \int_{\mathbb{R}_{+}} A(L_1,\ell,\ell)f(\ell)\, \ell\, \dd \ell\,.
\end{split}
\end{equation*}
\end{corollary}
\noindent \textbf{Proof.} We apply Theorem~\ref{intTRHTT} and there are only two things left to check. The first one is the strong admissibility of the twisted initial data, which is clear using the assumption \eqref{storngagum}.  The second one is the value of the integral $V\Omega_{1,1}[f](L_1)$, for which we use the technique explained in the proof of Theorem~\ref{intTRHTT} to handle cutting of tori with one boundary component. \hfill $\Box$

\vspace{0.2cm}

It is also possible to compute $V\Omega_{g,n}[f]$ by direct integration. For this purpose, we need to introduce the set $\mathcal{G}_{g,n}$ of stable graphs of type $(g,n)$.

\begin{definition}
A stable graph $G$ of type $(g,n)$ consists of the data of
\[
(\mathsf{V}_{G},\mathsf{H}_{G},\Lambda_{G},h,\mathsf{v},\mathsf{i})
\]
with the following properties.
\begin{itemize}
\item[$\bullet$] $\mathsf{V}_{G}$ is the set of vertices, equipped with a map $h\,\colon \mathsf{V}_{G} \rightarrow \mathbb{N}$.
\item[$\bullet$] $\mathsf{H}_{G}$ is the set of half-edges, $\mathsf{v}\, \colon \mathsf{H}_{G} \rightarrow \mathsf{V}_{G}$ is a map associating to each half-edge the vertex it is incident to, and $\mathsf{i}\, \colon \mathsf{H}_{G} \rightarrow \mathsf{H}_{G}$ is an involution.
\item[$\bullet$] $\mathsf{E}_{G}$ is the set of edges, consisting of the $2$-cycles of $\mathsf{i}$ in $\mathsf{H}_{G}$ (loops at vertices are allowed).
\item[$\bullet$] $\Lambda_{G}$ is the set of leaves, consisting of the fixed points of $\mathsf{i}$, and is equipped with a bijection to $\{1,\ldots,n\}$.
\item[$\bullet$] $(\mathsf{V}_{G},\mathsf{E}_{G})$ defines a connected graph, and we have the genus condition $g = \sum_{v \in \mathsf{V}_{G}} h(v) + b_1(G)$, where $b_1(G)$ is the first Betti number of the graph.
\item[$\bullet$] If $v$ is a vertex, $\mathsf{E}(v)$ (resp. $\Lambda(v)$) is the set of edges (resp. leaves) incident to $v$. Denoting $k(v) = |\mathsf{E}(v) \cup \Lambda(v)|$ the valency of $v$, we have the stability condition $2h(v) - 2 + k(v) > 0$.
\end{itemize}
An automorphism of $G$ are permutations of $\mathsf{V}_{G}$ and $\mathsf{H}_{G}$ which respect $h$, $\mathsf{v}$ and $\mathsf{i}$. We denote ${\rm Aut}(G)$ the group of automorphisms of $G$.
\end{definition}
We have a natural bijection $\mathcal{G}_{g,n} \cong M_{\Sigma}'/\Gamma_{\Sigma}^{\partial}$: For a primitive multicurve $c$, the edges of the corresponding stable graphs are the components of $c$, the vertices are the connected components of the cut surface $\Sigma^{c}$ and the map $h$ encodes their genera.

\begin{lemma}
\label{Lemstab} Under the assumptions of Theorem~\ref{costatunfg}, we have for $2g - 2 + n > 0$
\[
V\Omega_{g,n}[f](L_1,\ldots,L_n) = \sum_{G \in \mathcal{G}_{g,n}} \frac{1}{\#{\rm Aut}(G)} \int_{\mathbb{R}_{+}^{\mathsf{E}_{G}}} \! \prod_{v \in \mathsf{V}_{G}} \!\!V\Omega_{h(v),k(v)}\big((\ell_e)_{e \in \mathsf{E}(v)},(L_{\lambda})_{\lambda \in \Lambda(v)}\big)\prod_{e \in \mathsf{E}_{G}} \,f(\ell_e) \ell_e \, \dd \ell_e .
\]
\end{lemma}
\noindent \textbf{Proof.} Let $\Sigma$ be a bordered surface of type $(g,n)$ with boundary components labelled from $1$ to $n$. We will silently use Fubini's theorem to rearrange series and integrals, since the series we consider are absolutely convergent and integrable. Let $\overline{M}'_{\Sigma}$ be the set of primitive multicurves with ordered components of $\Sigma^c$, and if $c \in M'_{\Sigma}$ we denote $G(c)$ the corresponding stable graph. The fiber of the forgetful map $\overline{M}'_{\Sigma} \rightarrow M'_{\Sigma}$ above $c \in M'_{\Sigma}$ has cardinality $\#{\rm Aut}(G(c))$. From the definition \eqref{thedefOmgifsugn}, we obtain
\begin{equation*}
\begin{split}
V\Omega_{g,n}[f](L) & = \int_{\mathcal{M}_{g,n}(L)} \bigg(\sum_{c \in M'_{\Sigma}} \Omega_{\Sigma^{c}}(\sigma|_{\Sigma^{c}})\,\prod_{\gamma \in \pi_0(c)} f(\ell_{\sigma}(\gamma))\bigg)\dd\mu_{{\rm WP}}(\sigma) \\
& =  \int_{\mathcal{M}_{g,n}(L)} \bigg(\sum_{c \in \overline{M}'_{\Sigma}} \,\frac{1}{\#{\rm Aut}(G(c))}  \Omega_{\Sigma^{c}}(\sigma|_{\Sigma^{c}})\,\prod_{\gamma \in \pi_0(c)} f(\ell_{\sigma}(\gamma))\bigg)\dd\mu_{{\rm WP}}(\sigma).
\end{split}
\end{equation*}
For each $G \in \mathcal{G}_{g,n}$, we choose arbitrarily a primitive multicurve with ordered components $c(G)$  that represents $G$, and denote $\gamma_e^{G}$ its component corresponding to $e \in \mathsf{E}_{G}$. Then,
\begin{equation*}
\begin{split}
V\Omega_{g,n}[f](L) & = \sum_{G \in \mathcal{G}_{g,n}} \frac{1}{\#{\rm Aut}(G)} \int_{\mathcal{M}_{\Sigma}(L)} \bigg(\sum_{\tilde{c} \in \Gamma_{\Sigma}^{\partial}\cdot c(G)} \Omega_{\Sigma^{\tilde{c}}}(\sigma|_{\Sigma^{\tilde{c}}}) \prod_{\gamma \in \pi_0(\tilde{c})} f(\ell_{\sigma}(\beta))\bigg)\dd\mu_{{\rm WP}}(\sigma) \\
& = \sum_{G \in \mathcal{G}_{g,n}} \frac{1}{\#{\rm Aut}(G)} \int_{\mathcal{M}_{\Sigma}^{c(G)}(L)} \Omega_{\Sigma^{c(G)}}(\sigma|_{\Sigma^{c(G)}})\,\prod_{e \in \mathsf{E}_{G}} f(\ell_{\sigma}(\gamma_{e}^{G}))\, \dd \mu_{{\rm WP}}(\sigma).
\end{split}
\end{equation*}
The proof is then completed by applying the techniques explained in the proof of Theorem~\ref{intTRHTT}.
\hfill $\Box$

\vspace{0.2cm}

The results of this Section generalise straightforwardly to the target theory $E(\Sigma) = \mathscr{F}(\mathcal{T}_{\Sigma},\mathcal{A} \otimes \mathcal{A}^{* \pi_0(\partial_+\Sigma)})$ where $\mathcal{A}$ is a Frobenius algebra. In this case, we can choose $f$ to be a function from $\mathbb{R}_{+}$ to ${\rm Sym}^2(\mathcal{A})$, and in \eqref{thedefOmgifsugn} the factor $f(\ell_{\sigma}(\gamma)) \in {\rm Sym}^2(\mathcal{A})$ for each $\gamma \in \pi_0(c)$ is paired with the corresponding bivector coming from the two boundary components corresponding to $\gamma$ in $\Sigma^{c}$, using the pairing $\mathcal{A} \otimes \mathcal{A} \rightarrow \mathbb{C}$ and the corresponding isomorphism $\mathcal{A} \cong \mathcal{A}^*$.

\subsection{The systole function}

\label{sysf}

If $\epsilon > 0$ and $\sigma \in \mathcal{T}_{\Sigma}$, $S^{\circ,\epsilon}_{\Sigma,\sigma}$ is the set of simple closed curves $\gamma$ in $\Sigma$ such that $\ell_{\sigma}(\gamma) \leq \epsilon$. We denote $\tilde{S}^{\circ,\epsilon}_{\Sigma,\sigma}$ the set corresponding to the strict inequality. We first show that their cardinality can be accessed via the geometric recursion. We introduce the Heaviside functions
\[
\Theta(x) = \mathbf{1}_{[0,+\infty)}(x),\qquad \tilde{\Theta}(x) = \mathbf{1}_{(0,+\infty)}(x).
\] 

\begin{theorem} \label{sysL}
Let $\epsilon \in (0,{\rm argsinh}(1))$, $t \in \mathbb{R}$ and consider $f_{\epsilon}(\ell) = \Theta(\epsilon - \ell)$ and $\tilde{f}_{\epsilon}(\ell) = \tilde{\Theta}(\epsilon - \ell)$.  We have
\[
\Omega_{\Sigma}^{{\rm M}}[t\cdot f_{\epsilon}](\sigma) = (1 + t)^{\#S_{\Sigma,\sigma}^{\circ,\epsilon}},\qquad \Omega_{\Sigma}^{{\rm M}}[t\cdot \tilde{f}_{\epsilon}](\sigma) = (1 + t)^{\#\tilde{S}_{\Sigma,\sigma}^{\circ,\epsilon}}.
\]
In particular for $t = -1$, we obtain the indicator function of the (interior of) the $\epsilon$-thick part of Teichm\"uller space
\begin{equation}
\label{Omhnf}\Omega_{\Sigma}^{{\rm M}}[-f_{\epsilon}] =  \mathbf{1}_{(\mathcal{T}_{\Sigma}^{\epsilon})^{\circ}},\qquad \Omega_{\Sigma}^{{\rm M}}[-\tilde{f}_{\epsilon}] = \mathbf{1}_{\mathcal{T}_{\Sigma}^{\epsilon}}.
\end{equation}
\end{theorem}
\noindent \textbf{Proof.} For $\epsilon \in (0,{\rm argsinh}(1))$, the curves in $S_{\Sigma,\sigma}^{\circ,\epsilon}$ cannot intersect. Therefore $S_{\Sigma,\sigma}^{\circ,\epsilon}$ is finite and the only non-zero terms in the series \eqref{thedefOmgifsugn} defining $\Omega_{\Sigma}^{{\rm M}}[t\cdot f_{\epsilon}](\sigma)$ correspond to the primitive multicurves whose components form subsets $\mu \subseteq S_{\Sigma,\sigma}^{\circ,\epsilon}$. Hence
\[
\Omega_{\Sigma}^{{\rm M}}[t\cdot f_{\epsilon}](\sigma) = \sum_{\mu \subseteq S_{\Sigma,\sigma}^{\circ,\epsilon}} t^{\#\mu} = (1 + t)^{\#S_{\Sigma,\sigma}^{\circ,\epsilon}}.
\]
For $t = -1$ the value of this function equals $0$ when there exists a curve of length $\leq \epsilon$, and to $1$ otherwise. Equivalently, it is equal to $0$ if ${\rm sys}_{\sigma} \leq \epsilon$, and to $1$ otherwise, hence the first equality in \eqref{Omhnf}. The case of $\tilde{f}_{\epsilon}$ is proved similarly.
\hfill $\Box$

\subsection{Applications to string field theory} 

\label{ASFT}

Our Geometric Recursion construction of the systole functions has interesting applications in string field theory.
Let us here recall the bare rudiments of Zwiebach's formulation of closed string field theory (SFT), referring the reader to \cite{BZ} for a full account.

Part of this theory is the vertex Hilbert space $V$ with its inner product $\langle\cdot,\cdot\rangle$. Closed string field theory provides a family of multilinear maps (called brackets) indexed by $ g\geq 0, n\geq 0$ 
$$ [ \cdot, \ldots,  \cdot]_{g,n} :  V^ {\times n} \ra V,$$
which satisfy the quantum master equation (QME) \cite{BZ}.
These brackets are determined by the associated multi-pairings
$$ \{\cdot, \ldots, \cdot\}_{g,n} : V^{\times n} \ra {\mathbb C}$$
via the formula
$$ \forall v_1,\ldots,v_n \in V,\qquad  \{v_1, \ldots, v_n\}_{g,n} = \langle v_1, [v_2,\ldots, v_n]_{g,n-1}\rangle.$$
The latter takes the form
$$\{v_1, \ldots, v_n\}_{g,n} = \int_{V^\epsilon_{g,n}} \omega^\epsilon_{g,n}(v_1,\ldots, v_n),$$
where $\omega^\epsilon_{g,n}(v_1,\ldots, v_n)$ are certain top forms parametrised by a sufficiently small $\epsilon \in {\mathbb R}_+$ and 
where $V^\epsilon_{g,n}$ is a certain $\mathfrak{S}_n$-invariant subset of $\mathcal{M}_{g,n}$ which should satisfy the following version of the QME
$$
\partial V^\epsilon_{g,n} \cong \Bigg(\bigsqcup_{\substack{g_1+g_2 = g \\ J_1 \sqcup J_2 = \{1,\ldots,n+2\} \\ J_i \neq \emptyset}} V^\epsilon_{g_1,\# J_1} \times V^\epsilon_{g_2,\# J_2} \Bigg)\bigsqcup  V^\epsilon_{g-1, n+2},
$$
again parametrised by sufficiently small $\epsilon \in {\mathbb R}_+$. Taking into account symmetry factors, this induces at the level of chains on the orbifold $\mathcal{M}_{g,n}/\mathfrak{S}_n$ the standard form of the QME given in \cite[Equation 2.6]{BZ}. If we let $ \Omega^{\epsilon}_{g,n} = \Omega_{\Sigma}^{{\rm M}}[- f_{\epsilon}]$ be the function constructed in Section \ref{sysf}, we get that 

\begin{lemma}\label{Indic}
The function $ \Omega^{\epsilon}_{g,n}$ is the indicator function for the subset 
$$ \tilde V_{g,n}^\epsilon = \big\{ [\sigma] \in {\mathcal M}_{g,n} \mid \text{\rm{ {\small all simple interior closed geodesics on} }} [\sigma]  \text{\rm{ {\small have length at least}  }} \epsilon \big\}.$$
\end{lemma}
Letting $V_{g,n}^\epsilon = \tilde V_{g,n}^\epsilon \cap {\mathcal M}_{g,n}(\epsilon, \ldots, \epsilon)$ we arrive immediately at 
\begin{theorem}\label{sysQME}
The subsets $V_{g,n}^\epsilon$ satisfy the following quantum master equation
$$
\partial V^\epsilon_{g,n} \cong \Bigg(\bigsqcup_{\substack{g_1+g_2 = g \\ J_1 \sqcup J_2  = \{1,\ldots,n+2\} \\ J_i \neq \emptyset}} V^\epsilon_{g_1,\# J_1} \times V^\epsilon_{g_2,\# J_2} \Bigg)\bigsqcup  V^\epsilon_{g-1, n+2},
$$
thus form a topological vertex for string field theory. 
\end{theorem}

\begin{remark}

This theorem follows immediately from the fact that $  \Omega^{\epsilon}_{g,n} $ satisfies Geometric Recursion. The idea that the systole sets satisfies the QME was proposed by Kevin Costello at the workshop ``String Field Theory, BV quantisation and moduli spaces" at the Simons centre for Geometry and Physics, May 2019. We first reported on this result in \cite{JEAIHES}, where we presented the proof of Theorem \ref{sysL} and Lemma \ref{Indic} and concluded that two results implies Theorem \ref{sysQME}. Subsequently Costello and Zwiebach presented another proof of Theorem \ref{sysQME} in \cite{CZ}.

\end{remark}

Since $\Omega^{\epsilon}_{g,n}$ is the indicator function of $\tilde V^\epsilon_{g,n}$ we of course have that

$$\int_{V^\epsilon_{g,n}} \omega^{\epsilon}_{g,n}(v_1,\ldots, v_n) = \int_{{\mathcal M}_{g,n}(\epsilon, \ldots, \epsilon)}  \Omega^{\epsilon}_{g,n} \omega^{\epsilon}_{g,n}(v_1,\ldots, v_n).$$

In future work we shall seek to build $\Omega^{\epsilon}_{g,n} \omega^{\epsilon}_{g,n}(v_1,\ldots, v_n)$ via Geometric Recursion. This seems plausible since $\omega^{\epsilon}_{g,n}(v_1,\ldots, v_n)$ is build from the usual conformal field theory constructions, which satisfies factorisation and then we expect that the string brackets
$$\{v_1, \ldots, v_n\}_{g,n} = \int_{V^\epsilon_{g,n}} \omega^{\epsilon}_{g,n}(v_1,\ldots, v_n)$$
can be computed by topological recursion.

\section{The strict geometric recursion}
\label{S4strict}

Suppose we have chosen for a single reference object $\Sigma_{g,n}$ for each type $(g,n)$ in $\Surf_{1}$ or $\Surf_{\rm s}$. In the geometric recursion, the functoriality of $\Sigma \mapsto \Omega_{\Sigma}$ allows a non-ambiguous definition of
\[
w_{g,n} := \Omega_{\Sigma_{g,n}}
\]
by induction on $2g - 2 + n > 0$. Namely, it is possible to work solely with reference surfaces for each $g$ and $n$ all the way through the recursion. Passing from $\Omega_{\Sigma}$ to $w_{g,n}$, we have lost the memory of the structure of $\Surf_{1}$ or $\Surf_{{\rm s}}$, and therefore a lot of topological (and interesting) information and naturality of the construction.

In this section, we want to formalise the type of recursion satisfied by $(w_{g,n})_{g,n}$, which we call \emph{strict geometric recursion}. It is a kind of generalisation of the topological recursion of \cite{EORev,KSTR,TRABCD}, and is only based on finitely many glueing maps (Section~\ref{S41}). This construction contains less information than GR, but it can also be induced from the richer context of GR when the glueing maps are assembled from a converging series of isotopy class-dependent glueing maps as in \eqref{ThetaGT1} (see Section~\ref{S42}).

\subsection{Definition}
\label{S41}

We form a category $\overline{\Surf}_{1}$ which is a kind of strictification of $\Surf_{1}$. Its objects are finite (possibly empty) sequences $(g_i,n_i)_{i \in I}$ such that $2 - 2g_i - n_i < 0$ for all $i \in I$. Concatenation $\amalg$ of such sequences gives $\overline{\Surf}_{1}$ a monoidal structure, of which we only retain the structure of multicategory. There are no morphisms between distinct objects, and the automorphism group of $(g_i,n_i)_{i \in I}$ is the product of the permutation groups $\prod_{i \in I} \mathfrak{S}_{\llbracket 2,n_i\rrbracket}$. 

\begin{definition}
A strict target theory is a functor ${\rm e}:\, \overline{\Surf}_{1} \rightarrow \mathcal{C}$ together with functorial extra data satisfying the union and excision axioms below.
\end{definition}

\vspace{0.2cm}

\noindent \textsc{Union axiom.} For any objects $s$ and $s'$, we ask for the data of a continuous bilinear map
\[
\sqcup:\,{\rm e}(s) \times {\rm e}(s) \longrightarrow {\rm e}(s \amalg s'),
\]
compatible with commutativity and associativity of cartesian products. We require that ${\rm e}(\emptyset) = \mathbb{K}$ and the union morphism $\sqcup:\,{\rm e}(\emptyset) \times {\rm e}(s) \rightarrow {\rm e}(s)$ is specified by $1 \sqcup v = v$.

\vspace{0.2cm}

We need a few more notations before presenting the excision axiom. If $g \geq 0$ and $n \geq 1$ are such that $2g - 2 + n \geq 2$, we introduce the set $\mathcal{K}_{g,n}$ of $(0,3)$-excisions. It contains the following objects
\begin{itemize}
\item[\textbf{I}] $(g - 1,n + 1)$.
\item[\textbf{I}'] for each $j \in \llbracket 2,n \rrbracket$, a copy of the object $(g,n - 1)$ which we denote $(g,n-1)_j$.
\item[\textbf{II}] for each ordered partition $J \cup J' = \llbracket 2,n \rrbracket$ and ordered pair $(h,h')$ such that $h + h' = g$ with $2 - 2h - \#J < 0$ and $2 - 2h' - \#J' < 0$, a copy of the object $((h,1 + \#J),(h',1 + \#J'))$ which we denote $((h,1 + J),(h',1 + J'))$.
\end{itemize}
To any $\sigma \in \mathfrak{S}_{\llbracket 2,n \rrbracket}$ and $\kappa \in \mathcal{K}_{g,n}$, we assign a morphism $\tilde{\sigma}:\,\kappa \rightarrow \kappa_{\sigma}$ from $\kappa$ to the object $\kappa_{\sigma} \in \mathcal{K}_{g,n}$ in the following way.

\begin{itemize}
\item[\textbf{I}] Let $r:\,\llbracket 3,n + 1\rrbracket \rightarrow \llbracket 2,n\rrbracket$ be defined by $r(k) = k - 1$. We let $(g -1,n + 1)^{\sigma} := (g - 1,n + 1)$, and $\tilde{\sigma}$ is the morphism $r \circ \sigma \circ r^{-1}$.
\item[\textbf{I}'] Let $r$ be the unique strictly increasing map from $\llbracket 2,n - 1\rrbracket$ to $\llbracket 2,n \rrbracket \setminus\{j\}$. Then $\tilde{\sigma}$ is the morphism from the $j$-th copy of $(g,n - 1)$ to the $\sigma(j)$-th copy, and seen as an automorphism of the single object $(g,n - 1)$ in $\overline{\Surf}_{1}$, it is $r \circ \sigma \circ r^{-1}$.  
\item[\textbf{II}] If $K$ is an ordered set, let $r_{K}$ be the unique strictly increasing map from $K$ to $\llbracket 2,1 + |K|\rrbracket$. We take $((h,1 + J),(h',1 + J'))^{\sigma} := ((h,1 + \sigma(J)),(h',1 + \sigma(J')))$ and $\tilde{\sigma}$ considered as an automorphism of the single object $((h,1 + \#J),(h,1 + \#J'))$, is $(r_{\sigma(J)} \circ \sigma \circ r^{-1}_{J}) \amalg (r_{\sigma(J')} \circ \sigma \circ r^{-1}_{J'})$.  
\end{itemize}
If $\kappa = (g_i,n_i)_{i \in I}$ is an element of $\mathcal{K}_{g,n}$, we denote ${\rm e}^{\pi}(\kappa) := \prod_{i \in I} {\rm e}(g_i,n_i)$. From the previous construction, each $\sigma \in \mathfrak{S}_{\llbracket 2,n \rrbracket}$ induces morphisms ${\rm e}_{\kappa,\sigma}:\,{\rm e}^{\pi}(\kappa) \rightarrow {\rm e}(\kappa_{\sigma})$ for any $\kappa \in \mathcal{K}_{g,n}$.

\vspace{0.2cm}

\noindent \textsc{Excision axiom.} For any $(g,n)$ such that $2g - 2 + n \geq 2$ and any $\kappa \in \mathcal{K}_{g,n}$, we ask for the data of continuous multilinear morphisms
\[
\theta_{\kappa}:\,{\rm e}(0,3) \times {\rm e}^{\pi}(\kappa) \longrightarrow {\rm e}(g,n)
\]
which are compatible with the action of $\mathfrak{S}_{\llbracket 2,n\rrbracket }$ on the second factor of the source and on ${\rm e}(g,n)$.

\vspace{0.2cm}

\begin{definition}
A strict initial data is a quadruple $(A,B,C,D)$ such that
\[
A,C \in {\rm e}(0,3)^{\mathfrak{S}_{2}},\qquad B^2 \in {\rm e}(0,3),\qquad D \in {\rm e}(1,1) .
\]
\end{definition}

Let $(A,B,C,D)$ be strict initial data for a strict target theory ${\rm e}:\,\overline{\Surf}_{1} \rightarrow \mathcal{C}$.

\begin{definition}
\label{SDEFGR} We define the strict GR amplitudes as follows.

\begin{itemize}
\item[$\bullet$] We put $w_{\emptyset} := 1$, $ w_{0,3} := A$ and $w_{1,1} := D$.
\item[$\bullet$] For disconnected objects, we put $w_{(g_i,n_i)_{i \in I}} := \bigsqcup_{i \in I} w_{g_i,n_i}$.
\item[$\bullet$] If $2g - 2 + n \geq 2$, we seek to define inductively
\begin{equation}
\label{omgn}
\begin{split}
w_{g,n} & := \sum_{j = 2}^n \theta_{(g,n - 1)_j}(B^j,w_{g,n - 1}) \\
& \quad  + \frac{1}{2}\bigg(\theta_{(g - 1,n + 1)}(C,w_{g - 1,n + 1}) + \sum_{\substack{h + h'  = g \\ J \cup J' = \llbracket 2,n \rrbracket}}^{'}\theta_{((h,J),(h',J'))}(C,w_{h,|J|},w_{h',|J'|})\bigg)\,,
\end{split}
\end{equation}
where the $'$ in the last sum means we discard terms involving $(0,1)$ or $(0,2)$.
\end{itemize}
\end{definition}

\begin{proposition}
The assignment $(g_i,n_i)_{i\in I} \longmapsto \bigsqcup_{i \in I} w_{g_i,n_i} \in \prod_{i \in I} {\rm e}(g_i,n_i)$ is well-defined and functorial. In other words, $w_{g,n}$ is invariant under the action of ${\rm Aut}(g,n) = \mathfrak{S}_{\llbracket 2,n \rrbracket}$.
\end{proposition}

\noindent \textbf{Proof.} The proof proceeds by induction and is very easy. Indeed, the only terms in \eqref{omgn} which may not already be symmetric thanks to the induction hypothesis are $\theta_{(g,n - 1)_j}(B^j,w_{g_i,n_i - 1})$ indexed by $j \in \llbracket 2,n_i\rrbracket$, but the properties required for the glueing morphisms imply that the automorphism group $\mathfrak{S}_{\llbracket 2,n_i\rrbracket}$ of $(g_i,n_i)$ just permutes these terms. \hfill $\Box$

\subsection{From GR to strict GR}

\label{S42}

We can associate to any target theory $E$ a strict target theory ${\rm e}$. This is done by choosing, for each $(g,n)$, a reference surface $\Sigma_{g,n}$ in $\Surf_{1}$ and letting for $(g,n) \neq (0,3)$
\[
{\rm e}(g,n) := E(\Sigma_{g,n})^{\Gamma^{\partial}_{\Sigma_{g,n}}}.
\]
In the special case $(g,n)=(0,3)$ we specify below ${\rm e}(0,3)$ to be a certain subspace of $E(\Sigma_{0,3})$. The action of $\Gamma_{\Sigma_{g,n}}$ on $E_{\Sigma_{g,n}}$ factors into an action of $\mathfrak{S}_{\pi_0(\partial_+\Sigma_{g,n})}$ on ${\rm e}(g,n)$. Making the extra choice of an identification $\pi_0(\partial_+\Sigma_{g,n}) \simeq \llbracket 2,n \rrbracket$, this group defines the action of $\mathfrak{S}_{\llbracket 2,n \rrbracket}$ on ${\rm e}(g,n)$. We therefore obtain a functor ${\rm e}:\,\overline{\Surf}_{1} \rightarrow \mathcal{C}$.

We induce union morphisms for ${\rm e}$ from the union maps on $E$ in the obvious way. To define glueing morphisms for ${\rm e}$, we assemble the glueing morphisms for $E$. More precisely, if $\kappa \in \mathcal{K}_{g,n}$, let $\Sigma_{\kappa}$ denote, among the reference surfaces we have chosen, the one with topology $\kappa$. We let $\mathcal{P}_{g,n}(\kappa) \subseteq \mathcal{P}_{\Sigma_{g,n}}$ be the set of homotopy classes $[P]$ such that there exists a bijection $p:\,\pi_0(\Sigma_{\kappa}) \rightarrow \pi_0(\Sigma - P)$, which respects the ordering of the connected components when $\Sigma_{\kappa}$ is disconnected, together with isomorphisms $\varphi_{P}:\,\Sigma_{0,3} \rightarrow P$ in $\Surf_{1}$ and $\psi_P:=\big(\psi_{P,a}:\,\Sigma_{\kappa}(a) \rightarrow p(\Sigma_{\kappa}(a)))_{a}$. For each $[P] \in \mathcal{P}_{g,n}(\kappa)$ we pick such a $(p,\varphi_P,\psi_{P})$, and remark that any other choice is related to $(\varphi_P,\psi_{P})$ by post-composition with pure mapping classes on the connected components of $\Sigma_{0,3} \cup \Sigma_{\kappa}$.

Then we define
\begin{equation}
\label{varglumap} \theta_{\kappa} := \sum_{[P] \in \mathcal{P}_{g,n}(\kappa)} \Theta_{P} \circ \big(E(\varphi_P),E(\psi_{P,a})_{a \in \pi_0(\Sigma_{\kappa})}\big).
\end{equation}
We define ${\rm e}(0,3)$ chosen to be the subspace of $E(\Sigma_{0,3})$ on which this linear map is well defined and continuous.  we remark that the linear map $\theta_{\kappa}:\,{\rm e}(0,3) \times{\rm e}^{\pi}(\kappa) \rightarrow E(\Sigma_{g,n})$ does not depend on the choices of $([P],\varphi_P,\psi_{P})$ as above. And, since all elements in the sum \eqref{varglumap} are related by the action of the pure mapping classes of $\Sigma_{g,n}$, the morphism ${\rm e}_{\kappa}$ actually takes values in ${\rm e}(g,n) = E(\Sigma_{g,n})^{\Gamma_{g,n}}$.

Using \eqref{varglumap} we deduce the following result.
\begin{corollary}
\label{npin}Assume $E$ is a target theory, together with admissible initial data $(A,B,C,D)$, and denote $\Omega$ the outcome of GR. Construct a (non-canonical) strict target theory ${\rm e}$ as above and induce an ${\rm e}$-valued initial data by specialising $(A,B,C,D)$ to the reference $\Sigma_{0,3}$ and $D$ to the reference $\Sigma_{1,1}$. Denote $w$ the corresponding outcome of strict GR. Then for stable $(g,n)$ we have $w_{g,n} = \Omega_{\Sigma_{g,n}}$
\end{corollary}

Observe that by the decay axiom and Theorem \ref{converge}, the GR initial data does induces strict GR initial data. We note that the above construction is non-canonical as it involves the choice of reference surfaces $\Sigma_{g,n}$. However, if the target theory $E$ is such that $\Gamma_{\Sigma}^{\partial}$ acts trivially on $E(\Sigma)$, the construction is independent of these choices and therefore canonical. More generally, we have by an argument similar to Proposition~\ref{GRtoGR} that

\begin{proposition}
\label{GRtos}Let $E$ be a target theory, ${\rm \tilde{e}}$ a strict target theory, and $\eta:\,E \Longrightarrow \tilde{{\rm e}} $ a natural transformation which is compatible with the union and glueing morphisms. If $\mathcal{I} = (A,B,C,D)$ is an $E$-valued admissible initial data, then $\eta(\mathcal{I})$ is a strict initial data, and we have the relation between the corresponding GR and strict GR amplitudes
\[
w^{\eta(\mathcal{I})} = \eta(\Omega^{\mathcal{I}}) .
\] 
\end{proposition}

Noting that ${\rm e}(\Sigma) = {\rm Mes}(\mathbb{R}_{+}^{\pi_0(\partial \Sigma)})$ is a strict target theory in a natural way, Theorem~\ref{intTRHTT} can be interpreted by saying that the integration with respect the Weil--Petersson measure is a (partially defined) natural transformation between ${\rm e}$ and the $E(\Sigma) = {\rm Mes}(\mathcal{T}_{\Sigma})$.  The fact that geometric recursion implies topological recursion after integration expresses the compatibility with natural transformations stated in Proposition~\ref{GRtos}.

\subsection{Sums over fatgraphs}

\label{MainProp2}

In this paragraph, we show that the geometric recursion can be repackaged as a finite sum over fatgraphs of the appropriate topology with no reference to embeddings. The contribution of each fatgraph also satisfy a recursion under removal of the first vertex, albeit of non local nature. The similarity but also important differences with the topological recursion of \cite{EORev} will be spelled out in Section~\ref{S4strict}.

\vspace{0.4cm}

\noindent \textsc{Preliminary: fatgraphs and their doubled surface}

\vspace{0.2cm}

In this paragraph we assume $2g - 2 + n > 1$. We consider the set $\mathbb{G}^{g,n}$ of connected uni-trivalent fatgraphs with exactly $n$ univalent vertices and genus $g$. We denote by $\mathbb{G}^{g,n}_1$ the set of graphs $G$ in $\mathbb{G}^{g,n}$ equipped with a distinguished univalent vertex, called the \emph{root}.

\begin{figure}[ht!]
\begin{center} 
\includegraphics[width=0.4\textwidth]{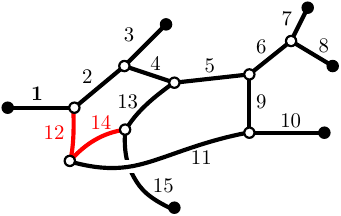}
\caption{\label{FIGBROWSE} The exploration on a fatgraph {$G \in \mathbb{G}^{2,6}_1$}. We start from the root \textbf{1}, and the numbers indicate in which order we then meet the (unoriented) edges. The spanning tree $\mathfrak{t}(G)$ is in black.} 
\end{center}
\end{figure}

We recall the construction of the \emph{canonical spanning tree} $\mathfrak{t}(G)$ for any element $G \in \mathbb{G}^{g,n}_1$ \cite{ABP} (see Figure~\ref{FIGBROWSE} for an example). We inductively define $\mathfrak{t}(G)$ by starting at the root in the direction of its incident edge, traveling on the edges of $G$ respecting the cyclic order of the half-edges at vertices, adding at each step to $\mathfrak{t}(G)$ the edge we encounter along the face we are traveling on, whenever it has not already been added and it does not create a loop. Once we come back to the root, we restart the same travel until we get to the first edge $e$, for which the face on the other side of the edge is different from the one we are traveling on currently. We then jump to this next face and continue the travel along this component, adding edges to $\mathfrak{t}(G)$ with the same rules. When we come back to the starting point of travel for this face next to the edge $e$, we iterate the process (travel till we can jump to an unexplored opposite face), until we cannot add more edges without creating a loop. The outcome is a rooted tree $\mathfrak{t}(G)$, which spans all the vertices of $G$.

\begin{figure}[ht!]
\begin{center}
\includegraphics[width=0.2\textwidth]{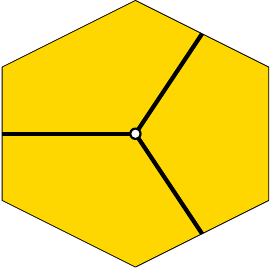}
\caption{\label{Hexafig} Hexagon attached to a trivalent vertex.}
\end{center}
\end{figure}

\begin{figure}[ht!]
\begin{center}
\includegraphics[width=0.6\textwidth]{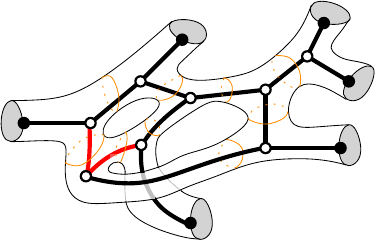}
\caption{\label{Doublefig} The double of the graph of Figure~\ref{FIGBROWSE}.}
\end{center}
\end{figure}
 
To any $G \in \mathbb{G}_{g,n}$ one can associate a bordered surface $\Sigma(G)$ of genus $g$ with $n$ boundary components, called the \emph{double} (Figure~\ref{Doublefig}). Let us briefly recall this construction. For each vertex of $v$, we take two copies $H_v^{+}$ and $H_{v}^-$ of a closed hexagon in which $v$ and its incident half-edges are embedded such that they end on three (alternating) sides of the hexagon (Figure~\ref{Hexafig}). If there is an edge from a univalent vertex $u$ to $v$, we ask that the edge from $u$ to $v$ is embedded in $H_v^{\pm}$ such that $u$ maps to a point in the interior of a side of the hexagon. We orient $H_v^+$ (respectively $H_v^-$) so as to agree (respectively disagree) with the cyclic order of half-edges around $v$. We then glue together $H_{v}^+$ and $H_v^-$ into a pair of pants $P_{v}$ by identifying pairs $(s_-,s_+)$ of sides with opposite orientation in $(H_{v}^{-},H_{v}^+)$ which are not crossed by edges incident to $v$. The surface $\Sigma(G)$ is then obtained by the quotient of $\bigcup_{v} P_{v}$, where we for each edge in $G$, say between the vertices $v$ and $v'$, identify the boundary components of $P_{v}$ and $P_{v'}$ corresponding to the edge. In particular, two copies of $G$ are embedded in $\Sigma(G)$, and the image $\Sigma^+(G)$ of $\bigcup_{v} H_{v}^+$ in $\Sigma(G)$ deformation retracts onto $G$.

Since two different ways of choosing the topological spaces $H_{v}^{\pm}$ and the glueing are related by diffeomorphisms in a unique orientation-preserving isotopy class, $\Sigma(G)$ can be considered as a single object in $\Surf$ up to canonical isomorphisms.

If $G \in \mathbb{G}^{g,n}_{1}$, reversing the orientation of the boundary component of $\Sigma(G)$ containing the root gives a object in $\Surf_{1}$ which is uniquely defined up to canonical isomorphisms. From here on $G$ always stands for a graph in $\mathbb{G}_1^{g,n}$ and $\Sigma(G)$ the associated object in $\Surf_{1}$.
 
Let $\widetilde{\mathbb{M}}^G(\Sigma)$ be the set of morphisms in $\Surf_{1}$ from $\Sigma(G)$ to $\Sigma$. We observe that an element of $\widetilde{\mathbb{M}}^{G}(\Sigma)$ induces a pair of pants decomposition of $\Sigma$, up to ambient isotopy. Thus there is a well-defined subgroup $\Gamma^G_{\Sigma}$ associated to each $G$, which is the subgroup of $\Gamma_{\Sigma}$ preserving each pair of pants as well as each boundary component of a pair of pants in this decomposition. We denote the quotient by 
\[
\mathbb{M}^G(\Sigma) = \widetilde{\mathbb{M}}^{G}(\Sigma)/\Gamma^G_{\Sigma},
\]
and we will write 
\[
\widetilde{\mathbb{M}}(\Sigma) = \bigcup_{G \in \mathbb{G}^{g,n}_{1}}\widetilde{\mathbb{M}}^{G}(\Sigma),\qquad \mathbb{M}(\Sigma) = \bigcup_{G \in \mathbb{G}^{g,n}_{1}}\mathbb{M}^{G}(\Sigma).
\]

\vspace{0.4cm}

\noindent \textsc{GR as a sum over fatgraphs}

\vspace{0.2cm}

The spanning tree $\mathfrak{t}(G)$ of $G$ allows us to determine an ordering on the embedded pair of pants $(P_1,\ldots P_{2g-2+n})$ as well as an ordering of the set of boundary component of each $P_i$. We can therefore consider each $P_i$ as an object in $\Surf_{1}$, by declaring the orientation of its first boundary to disagree with the orientation of the surface. We put $\Sigma(G)_1 = \Sigma(G)$, and successively excise the pair of pants $P_i$, for $i = 1,2,\ldots$ to obtain a sequence of surfaces inductively defined by $\Sigma(G)_{i + 1} = \Sigma(G)_i - P_i$. When excising $P_i$, we can distinguish the following four mutually exclusive cases.

\begin{description}
\item[A --] $\partial P_i \subseteq \partial \Sigma(G)_i$.

\item[B --] Exactly two of the boundary components of $P_i$ are boundary components of $\Sigma(G)_i$.

\item[C --] Exactly one of the boundary components of $P_i$ is a boundary component of $\Sigma(G)_i$, and the two other boundary components of $P_i$ are disjointly embedded in $\Sigma(G)_i$.

\item[D --] Two of the boundary components of $P_i$ are glued together under the embedding into $\Sigma(G)_i$. In this case, by glueing these two boundary components we replace $P_i$ by a torus with one boundary component (which we keep denoting $P_i$ for convenience).

\end{description}

This defines for us a type map $X:\,\llbracket 1,2g-2+n \rrbracket \longrightarrow \{A,B^{b},C,D\}$. In the $B$ situation, the type also records the second boundary of the corresponding pair of pants with the letter $b$. We denote
\begin{equation}
\label{autfactor} |{\rm Aut}(G)| = 2^{X^{-1}(C)}.
\end{equation}

In a given target theory $E$, we introduce the composed glueing morphism
\[
\Theta^{G} = \Theta_{P_{2g - 2 + n}} \circ \cdots \circ \Theta_{P_1}:\,\prod_{i = 1}^{2g - 2 + n} E(P_i) \longrightarrow E(\Sigma(G)).
\] 
For given admissible initial data, we define $\omega_G \in E(\Sigma(G))$ by
\[ 
\omega^G  = \Theta^{G}\big((X_{P_i}(i))_{i = 1}^{2g - 2 + n}\big).
\]
For a $\tilde{\varphi} \in \widetilde{\mathbb{M}}(\Sigma)$ we have of course a morphism $E(\tilde{\varphi}):\, E(\Sigma(G)) \rightarrow E(\Sigma)$. Since the mapping classes in $\Gamma^{G}_{\Sigma}$ leave stable each boundary component of the $P_i$s, $\omega^{G}$ is $\Gamma^{G}_{\Sigma}$-invariant. Therefore $E(\tilde{\varphi})(\omega^G) \in E(\Sigma)$ only depends on the projection $\varphi \in \mathbb{M}(\Sigma)$ of $\tilde{\varphi}$ and can be denoted $E(\varphi)(\omega^{G})$.
\begin{proposition} 
\label{OmegaPP} For any surface $\Sigma$ in $\Surf_{1}$, the GR amplitudes satisfy
\[
\Omega_{\Sigma} = \sum_{G\in \mathbb{G}^{g,n}_1} \frac{\Omega_{\Sigma}^{G}}{|{\rm Aut}(G)|},\qquad \Omega_{\Sigma}^{G} = \sum_{\varphi \in \mathbb{M}_G(\Sigma)} E(\varphi)(\omega^G).
\]
\end{proposition} 
\noindent \textbf{Proof.} Unfolding the recursive definition \eqref{defGR} and thanks to the property of absolute convergence from Theorem~\ref{converge}, we can apply the Fubini theorem and obtain a sum over sequences $([P_1],\ldots,[P_{m}])$ where $[P_i] \in \mathcal{P}_{\Sigma^{(i)}}$, with $\Sigma^{(0)} = \Sigma$ and $\Sigma^{(i)} = \Sigma^{(i - 1)} - P_{i}$. At each step, the excision of a pair of pants produces the isotopy class of embeddings $[P_i] = [f_i:\,P_i \rightarrow \Sigma^{(i - 1)}]$ of a pair of pants with ordered boundary components. What is left, when no more excisions are possible (recall we stop before getting the empty surface), are connected components which are either pairs of pants (an $A$ factor) or tori with one boundary component (a $D$ factor). Tracking the ordering of boundary component of the $P_i$ in this nested structure reveals that this set of sequences is in bijection with $\mathbb{M}(\Sigma)$. For each $\varphi \in \mathbb{M}(\Sigma)$, we can then rearrange the nested application of glueing maps into the composite glueing map $E(\varphi) \circ \Theta^{G}$. It should be applied to a tuple formed by the factors contributing to each embedded piece $P_i$, either an $A$, $B$, $C$ or $D$, depending only on the abstract fatgraph $G$ as specified  by the type map. The automorphism factor \eqref{autfactor} takes into account the conventional power of $1/2$ attached to each $C$ in \eqref{GRdef}. \hfill $\Box$

\vspace{0.4cm}

\noindent \textsc{Fatgraph recursion}

\vspace{0.2cm}

We can also write a fatgraph recursion, but we stress that the glueing maps will be non local, \textit{i.e.} they \textit{a priori} not only depend on the excision of a trivalent vertex, but on the full structure of the fatgraph. This makes our recursion quite different from usual fatgraph recursions encountered in enumerative geometry, as we will comment on in Section~\ref{S4strict}.

Let us define $E_G =  E(\Sigma(G))$ and let
\begin{equation}
\label{OmG} \Omega^{G} = \Omega^{G}_{\Sigma(G)} = \sum_{\varphi \in \mathbb{M}^{G}(\Sigma(G))} E(\varphi)(\omega^G) \in E_{G}.
\end{equation}

Let $\mathfrak{t}_{1} \subset \mathfrak{t}(G)$ be the union of the univalent vertex $u$, the trivalent vertex $v$ incident to it, and the half edges attached to these two vertices. $G-\mathfrak{t}_{1}$ may or may not be a connected graph. If it is connected we put $\tilde{E}_{G - \mathfrak{t}_1} := E_{G - \mathfrak{t}_1}$. If it is not connected, $G - \mathfrak{t}_1$ will consist of an ordered (by cyclic order at $v$ starting with $u \rightarrow v$) pair of connected components $(G',G'')$, and we define $\tilde{E}_{G - \mathfrak{t}_1} := E_{G'} \times E_{G''}$. In both cases, we construct a \emph{fatgraph glueing map}
\[
\Theta^{G,\mathfrak{t}_{1}}:\,E_{\mathfrak{t}_1} \times \tilde{E}_{G - \mathfrak{t}_{1}} \longrightarrow E_{G}
\]
as follows. Let $\mathcal{P}(G,\mathfrak{t}_1)$  be the set of isotopy classes of all embeddings of the pair of pants $\Sigma(\mathfrak{t}_{1})$ into $\Sigma(G)$ preserving the orientation of the surfaces and of their boundary components. For each $[p] \in \mathcal{P}(G,\mathfrak{t}_{1})$ we choose an orientation-preserving identification $\varphi_p$ of the complement of $p$ with $\Sigma(G - \mathfrak{t}_{1})$. This way we obtain glueing map $\Theta_{p}:\,E_{\mathfrak{t}_{1}}\times \tilde{E}_{G - \mathfrak{t}_1} \ra  E_{G}$, which does not depend on the choice of $\phi_p$, and we define
\begin{equation}
\label{ThetaGT1} \Theta^{G,\mathfrak{t}_1} := \sum_{[p] \in \mathcal{P}(G,\mathfrak{t}_{1})} \Theta_{p}\,.
\end{equation}
The right-hand side is only defined when evaluated on certain elements of $E_{\mathfrak{t}_1} \times \tilde{E}_{\tilde{G}}$ making the series absolutely convergent.

\begin{proposition} For any $G \in \mathbb{G}^{g,n}_{1}$ with $2g - 2 + n \geq 2$, we have $\Omega^G = \Theta^{G,\mathfrak{t}_1}(\Omega^{\mathfrak{t}_1}, \Omega^{G - \mathfrak{t}_{1}})$.
\end{proposition}

\noindent \textbf{Proof.} As the series \eqref{OmG} is absolutely convergent, it can be repackaged as
\[
\Omega^{G} = \sum_{[p] \in \mathcal{P}(G,\mathfrak{t}_{1})} \sum_{\substack{\varphi \in \mathbb{M}^{G}(\Sigma(G)) \\  \varphi|_{\mathfrak{t}_1} = p}} E(\varphi)(\omega^G).
\]
By the definition of $\Theta^{G,\mathfrak{t}_{1}}$, $\Omega^{G - \mathfrak{t}_{1}}$ and $\Omega^G$, this sum may be rewritten as $\Omega^{G} = \Theta^{G,\mathfrak{t}_1}(\Omega^{\mathfrak{t}_{1}}, \Omega^{G - \mathfrak{t}_{1}})$. \hfill $\Box$


Recall that any $G \in \mathbb{G}_{1}^{g,n}$ determines a pair of pants decomposition, with ordered pair of pants $(P_1,\ldots,P_{2g - 2 + n})$ of $\Sigma(G)$, as well as a type map $X:\,\llbracket 1,2g - 2 + n \rrbracket \rightarrow \{A,B^{j},C,D\}$. In a strict target theory ${\rm e}$ we have a glueing morphism
\[
\theta^{G}:\,\prod_{i = 1}^{3g - 3 + n} {\rm e}(P_i) \longrightarrow {\rm e}(g,n) .
\]
Unfolding the recursive definition \eqref{omgn} as we did in Section~\ref{MainProp2} (without the complications due to the consideration of isotopy classes) yields

\begin{proposition}
\label{GRAPHTR} Let $(A,B,C,D)$ be an admissible initial data for a strict target theory ${\rm e}$. For any stable $(g,n)$, the strict GR amplitude takes the form
\[
w_{g,n} = \sum_{G \in \mathbb{G}_{1}^{g,n}} \theta^{G}\big((X_{P_i})_{i = 1}^{3g - 3 + n}\big).
\]
\end{proposition}

 \section{Relation to topological recursion}
 
\label{Section7}
 
This section gives examples of strict symmetric GR which demonstrate that the topological recursion (TR) and its many variants are particular cases of the strict symmetric GR.  The most obvious class of example comes from two-dimensional topological quantum field theories (TQFTs). We then explain separately the class of examples governing the original topological recursion of \cite{EOFg} which is based on the geometry of complex curves, and then the Kontsevich-Soibelman approach to topological recursion (KS-TR) \cite{KSTR,TRABCD} based on algebraic structures called ``quantum Airy structures''. Although KS-TR contains the other two examples as particular cases, it is instructive to see in each different language how the target theory should be formed.
 
The strict GR can be described in some sense as a (possibly non-symmetric) generalisation of the topological recursion (TR) of \cite{EOFg,KSTR} in which the glueing maps are allowed to ``increase complexity'', aligning with the possibly increasing complexity of the spaces ${\rm e}(g,n)$ when $2g - 2 + n$ increases. This should be compared with Section~\ref{2dTQFT}-\ref{Airystr}, where they do not depend on $g$. Adopting this perspective, we may say in light of Corollary~\ref{npin} that a natural transformation from a target theory to a strict target theory maps GR to TR. The richer structure which appears in the construction of GR due to infinite sums over homotopy class and mapping class group considerations is not seen at the level of the corresponding TR, as everything has been incorporated into finitely many glueing maps and the action of mapping class groups has been replaced by permutation groups.
 
Conversely, if we start from a strict target theory $E$ and initial data, we will describe in \S~\ref{FIBERING} a way to construct some target theory $E$ and initial data such that GR recomputes the outcome of strict GR. This is based on the idea of fibering $E$ over Teichm\"uller spaces, and on Mirzakhani's generalisation of the McShane identity.

We nevertheless stress that target theories are allowed to contain much more information than their strict counterpart. Indeed, they are functors from the category $\Surf_{1}$ which reflects $(2 + 1)$-dimensional topology (surfaces and their mapping classes) while $\overline{\Surf}_{1}$ only reflects the classification of $2$-dimensional compact orientable manifolds and therefore has a simple combinatorial structure. One can find many examples of $E(\Sigma)$ having a geometric meaning and carrying naturally non-trivial mapping class group actions. The gain between GR and strict GR will become clear in Chapter II (and III), where we construct and use the first non-trivial examples of target theories based on the topology (and the geometry) of the Teichm\"uller spaces, and provide examples of natural transformations from GR to strict GR (in application of Proposition~\ref{GRtos}) by integration over moduli spaces.

\subsection{Example from 2d TQFTs}
\label{2dTQFT}
We consider 2d TQFTs in the category ${\rm Vect}_{\mathbb{C}}$ of finite-dimensional complex vector spaces to illustrate the idea -- which is standard and elementary -- but one may try to adapt a similar construction for more fancy, higher-categorical notions of TQFTs.

A 2d TQFT is equivalent to the data of a unital Frobenius algebra $\mathcal{A}$, \textit{i.e.} an object in ${\rm Vect}_{\mathbb{C}}$ equipped with a pairing and a commutative associative product which is invariant for the pairing -- see \textit{e.g.} \cite{At2,Abrams}. Using the pairing, we have canonical identifications between $\mathcal{A}$ and $\mathcal{A}^*$, and therefore distinguished elements
\[
b^{\dagger} \in {\rm Hom}({\rm Sym}^2(\mathcal{A}^*),\mathbb{C}),\qquad \mu \in {\rm Sym}^3(\mathcal{A}^*),
\]
representing the pairing and the product. The TQFT amplitudes are the element 
\[
F_{g,n} \in {\rm Hom}(\mathcal{A}^{\otimes n},\mathbb{C})
\]
defined for any $(g,n)$ as follows. Take a pair of pants decomposition $(P_1,\ldots,P_{2g - 2 + n})$ of a connected surface of genus $g$ with $n$ boundaries, and form the tensor product
\[
\bigg(\bigotimes_{i = 1}^{2g - 2 + n} \mu\bigg) \in (\mathcal{A}^*)^{6g - 6 + 3n}.
\]
Then, each time $P_i$ has a common boundary component with $P_j$, we use the pairing to contract a copy of $\mathcal{A}$ in the $i$-th factor together with a copy of $\mathcal{A}$ in the $j$-th factor. This makes a total of $3g - 3 + n$ pairings, so we are left with an element of $(\mathcal{A}^*)^{\otimes n} = {\rm Hom}(\mathcal{A}^{\otimes n},\mathbb{C})$ which is by definition $F_{g,n}$. The axioms of a Frobenius algebra implies that the result is independent of the choice of a pair of pants decomposition, and thus only depends on $g$ and $n$, and is symmetric under permutation of the $n$ remaining factors of $\mathcal{A}$.

\vspace{0.2cm}

\noindent \textsc{Strict symmetric target theory.} ${\rm Vect}_{\mathbb{C}}$ is a rather simple full subcategory of our category $\mathcal{C}$. For each stable $(g,n)$, we choose ${\rm e}(g,n) = {\rm Hom}(\mathcal{A}^{\otimes n},\mathbb{C})$. The union map comes from the tensor product (in this case the functor ${\rm e}$ is monoidal) and the morphisms are composition of tensors -- using as many times as necessary the canonical identification $\mathcal{A} \cong \mathcal{A}^*$.

\vspace{0.2cm}

\noindent \textsc{Initial data.} We take $A = B = C = \mu$ representing the product. We also take $D =  \langle H,\cdot \rangle$ where $H := \sum_{i} e_i^2$ in terms of an orthogonal basis $(e_i)_{i}$. The four relations characterising symmetric initial data are satisfied because, in each of them, each of the three terms in the left-hand side are separately symmetric under the permutation $\sigma_{1,2}$ due to the properties of the product.

\vspace{0.2cm}

\noindent \textsc{The amplitudes.} With this initial data, it is easy to see that $w_{1,1} = F_{1,1} = D$, and to prove recursively that, if one writes the $w_{g,n}$ as sum over the set $\mathbb{G}_{g,n}$ of fatgraphs described in Section~\ref{MainProp2}, each term uniquely defines a topological class of a pair of pants decomposition, and the value assigned to each fatgraph is by construction the TQFT amplitude $F_{g,n}$ -- computed with the help of this pair of pants decomposition. So, for any $2g - 2 + n > 0$, the strict GR amplitudes are
\[
w_{g,n} = \#\mathbb{G}^{g,n}_1 \cdot F_{g,n} .
\]

\subsection{Example from quantum Airy structures}
\label{Airystr}
Another class of examples of strict symmetric GR is provided by the approach of Kontsevich-Soibelman approach to topological recursion \cite{KSTR,TRABCD}.

We first outline this theory. Let $V$ be a vector space over $\mathbb{C}$, which we here assume to be finite-dimensional for simplicity. One can easily handle infinite-dimensionality if equipped with an increasing filtration by finite-dimensional subspaces. Let $\mathcal{W}_{V}^{\hbar}$ be the Weyl algebra of $V$, \textit{i.e.} the unital algebra generated over $\mathbb{C}[[\hbar]]$ by $T^*V = V \oplus V^*$ with relations
\[
\forall (v,\lambda) \in V \times V^*,\qquad [v,\lambda] = \hbar\,\lambda(v) \in \mathbb{C}[[\hbar]] .
\]

\begin{definition}
\label{qA} A \emph{quantum Airy structure} is the data of a linear map $L:\,V \rightarrow \mathcal{W}_{V}^{\hbar}$ such that there exists a basis $(x_i)_{i \in I}$ of linear coordinates on $V$ in which $L$ takes the form
\[
L_i = \hbar \partial_{x_i} - \sum_{a,b \in I} \big(\tfrac{1}{2} A^i_{a,b}x_{a}x_{b} - \hbar B^i_{a,b}x_{a}\partial_{x_b} - \tfrac{\hbar^2}{2} C^i_{a,b}\partial_{x_a}\partial_{x_b}\big) - \hbar D^i\,,
\]
and satisfies the Lie algebra relations
\begin{equation}
\label{Liecomut}\forall (i,j) \in I^2,\qquad  [L_i,L_j] = \sum_{a \in I} \hbar\,f_{i,j}^a\,L_{a}
\end{equation}
for some scalars $A^i_{j,k} = A^i_{k,j}$, $B^i_{j,k}$, $C^i_{j,k} = C^i_{k,j}$, $D^i$ and $f^k_{i,j} = -f^k_{j,i}$.
\end{definition}
This notion is closely related to the quantisation of Lagrangians in $T^*V$ which are tangent to the zero section at $0$ and defined by quadratic equations. The constraints \eqref{Liecomut} impose relations on the coefficients of $L$.

\begin{lemma} \cite{TRABCD}
\label{Relel}Equation \eqref{Liecomut} is equivalent to the system of equations indexed by $i,j,k,\ell \in I$
\begin{itemize}
\item[$\bullet$] $A^i_{j,k} = A^j_{i,k}$;
\item[$\bullet$] $f^k_{i,j} = B^i_{j,k} - B^j_{i,k}$;
\item[$\bullet$] $\sum_{a \in I} B^i_{j,a}A^a_{k,\ell} + B^i_{k,a}A^j_{a,\ell} + B^i_{\ell,a}A^j_{a,k} = (i \leftrightarrow j)$;
\item[$\bullet$] $\sum_{a \in I} B^i_{j,a}B^a_{k,\ell} + B^i_{k,a}B^j_{a,\ell} + C^i_{\ell,a}A^j_{a,k} = (i \leftrightarrow j)$;
\item[$\bullet$] $\sum_{a \in I} B^i_{j,a}C^a_{k,\ell} + C^i_{k,a}B^j_{a,\ell} + C^i_{\ell,a}B^j_{a,k} = (i \leftrightarrow j)$;
\item[$\bullet$] $\sum_{a \in I} B^i_{j,a}D^a + \tfrac{1}{2} \sum_{a,b \in I} C^i_{a,b}A^j_{a,b} = (i \leftrightarrow j)$.
\end{itemize}
\end{lemma}

The third, fourth and fifth relation have the same index structure, and take the form of H-I-X relations. This was explained from a Lie algebraic perspective in \cite{TRABCD}. The structure of the axioms for the initial data of a symmetric GR was modelled on these relations, and give them a geometric content as coming from different (homotopy class of) pair of pants decompositions of a sphere with four boundaries.

We can also assemble the coefficients of a quantum Airy structure into tensors
\[
A \in {\rm Hom}(V^{\otimes 3},\mathbb{C}),\qquad B \in {\rm Hom}(V \otimes V,V),\qquad C \in {\rm Hom}(V,V \otimes V),\qquad D \in {\rm Hom}(V,\mathbb{C}),
\]
and the relations in Lemma~\ref{Relel} are then tensorial relations, where the sum over intermediate indices are replaced by composition of linear maps.

To any quantum Airy structure, one may associates amplitudes $F_{g,n} \in {\rm Hom}(V^{\otimes n},\mathbb{C}$) in the following way.
\begin{theorem} 
\label{THF} \cite{KSTR,TRABCD} Let $L$ be a quantum Airy structure. There exists a unique $F \in \hbar^{-1}({\rm Sym}\,V^*)[[\hbar]]$ without constant term, which we can decompose into
\[
F = \sum_{\substack{g \geq 0 \\ n \geq 1}} \frac{\hbar^{g - 1}}{n!} F_{g,n},\qquad F_{g,n} \in {\rm Sym}^n V^*\,,
\]
and such that $F_{0,1} = 0$, $F_{0,2} = 0$ and
\[
\forall i \in I,\qquad L_i \cdot \exp(F) = 0.
\]
In fact, $F_{g,n}$ is uniquely determined by the initial data $F_{0,3} = A$ and $F_{1,1} = D$, and a recursion on $2g - 2 + n > 0$. If we denote $(e_i)_{i \in I}$ the basis of $V$ mentioned in Definition~\ref{qA}, and
\[
F_{g,n}[i_1,\ldots,i_n] := F_{g,n}(e_{i_1} \otimes \cdots \otimes e_{i_n}),
\]
this recursion takes the form
\begin{equation}
\label{KSTRRR}
\begin{split} F_{g,n}[i_1,\ldots,i_n] & = \sum_{m = 2}^n \sum_{a \in I} B^{i_1}_{i_m,a} F_{g,n - 1}[a,i_{\llbracket 2,n \rrbracket \setminus\{m\}}]  \\
& \quad  + \sum_{a,b \in I} \tfrac{1}{2}\,C^{i_1}_{a,b}\bigg(F_{g - 1,n + 1}[a,b,i_{\llbracket 2,n\rrbracket}] +  \sum_{\substack{J \cup J' = \{i_2,\ldots,i_n\} \\ h + h' = g}} F_{g,1 + \#J}[a,J]F_{g',1 + \#J'}[b,J']\bigg)\,.
\end{split}
\end{equation}
\end{theorem}  
Here $\exp(F)$ is sometimes called the partition function or wave function, and this theorem characterises it as the solution of differential constraints forming a Lie algebra. The only point which may not be obvious is the existence of $F$ -- which amounts to prove symmetry of the $F_{g,n}$ obtained by the recursion \eqref{KSTRRR}. This theorem is proved by a general argument in \cite{KSTR}, and by direct computation from Lemma~\ref{Relel} in \cite{TRABCD}.

\vspace{0.2cm}

\noindent \textsc{Target theory.} For any stable $(g,n)$, we put ${\rm e}(g,n) = {\rm Hom}(V^{\otimes n},\mathbb{C})$. The union map is provided by the tensor product, and the glueing map by the composition of linear maps.

\vspace{0.2cm}

\noindent \textsc{The amplitudes.} The constraints on the tensors $(A,B,C,D)$ described in Lemma~\ref{Relel} are exactly the ones characterising symmetric initial data. This is in fact the reason why we have postulated these relations in Definition~\ref{inininis}, as well as their non-strict counterpart in Definition~\ref{DEFGRSSSS}. The proof of Theorem~\ref{SGRSYMTH} (of which Proposition~\ref{thttsh} was a direct consequence) follows the same scheme as the proof of symmetry of $F_{g,n}$ in \cite{TRABCD}, except that it incorporates homotopy class considerations (which are absent in the strict version). Comparison of \eqref{KSTRRR} with Definition~\ref{SDEFGR} shows directly that $F_{g,n} = w_{g,n}$.

\begin{remark}
The third constraint in Lemma~\ref{Relel} motivates the adoption of the BB-CA relation instead of the relation (\ref{BBCAprime}) in the definition of symmetric initial data for GR. Indeed, only the BB-CA relation implies this third constraint when we apply a natural transformation to a strict target theory.
\end{remark}

\subsection{Relation to Eynard--Orantin formalism}
\label{EOfmum}
The original topological recursion of Eynard--Orantin \cite{EORev} and all its variants formulated in terms of spectral curves determine quantum Airy structures \cite{KSTR,TRABCD} and therefore provide examples of strict symmetric GR. For later use in Section~\ref{GUNIGUNFIOGSO}, we briefly describe this relation. Spectral curves $(S,x,y,\omega_{0,2})$ with simple ramification points provide examples of Airy structures on the vector space $V = \zeta \mathcal{A}[[\zeta^2]]$, where $\mathcal{A} = \oplus_{r \in \mathfrak{r}} \mathbb{C}.e_{r}$ and $\mathfrak{r} \subset S$ is the set of ramification points of $x$. We define a basis of $V$
\[
\xi^*_{r,d}(\zeta) = \frac{\zeta^{2d + 1}}{(2d + 1)!}\,e_{r},\qquad (d,r) \in \mathbb{N} \times \mathfrak{r}.
\] 
We have $T^*V = \mathcal{A}[\zeta^{-2},\zeta^2]]\dd \zeta$ and we make it a symplectic vector space with the skew-symmetric $2$-form
\[
\langle f_1,f_2 \rangle = \Res_{\zeta = 0} \bigg(f_1 \bigg| \int^{\zeta} f_2\bigg)_{\mathcal{A}},
\]
where $(\cdot | \cdot)_{\mathcal{A}}$ is the symmetric  pairing making $(e_r)_{r \in \mathfrak{r}}$ orthonormal. We already have an isotropic subspace $\dd V \subset T^*V$, and the data of $\omega_{0,2}$ specifies another isotropic subspace $V'$ such that $\dd V \oplus V' = T^*V$, as follows. We first define the basis of $1$-forms on $S$ indexed by $d \geq 0$ and $r \in \mathfrak{r}$
\begin{equation}
\label{xidkbasis} \dd\hat{\xi}_{r,d}(z_0) = \Res_{z = r} \frac{(2d + 1)!\,\dd \zeta(z)}{\zeta(z)^{2d + 2}} \int_{r}^{z} \omega_{0,2}(\cdot,z_0),
\end{equation}
where we use the local coordinate $\zeta(z) = \sqrt{2(x(z) - x(r))}$ when $z$ is in the neighborhood of $r$. The Laurent expansion of $\dd\hat{\xi}_{r,d}$ near $z \rightarrow r'$ defines an element $\xi_{r,d}^{(r')} \in \mathbb{C}[\zeta^{-1},\zeta^1]\dd \zeta$. We then set
\[
\xi_{r,d}(\zeta) = \sum_{r' \in \mathfrak{r}} \frac{\xi_{r,d}^{(r')}(\zeta) - \xi_{r,d}^{(r')}(-\zeta)}{2}\,e_{r'},\qquad V' := {\rm span}_{\mathbb{C}}\big(\xi_{r,d}\quad : \quad (d,r) \in \mathbb{N} \times \mathfrak{r}\big),
\]
which indeed gives an isotropic subspace supplement to $\dd V \subset T^*V$. The amplitudes $F_{g,n} \in \mathcal{A}^{\otimes n}$ of Theorem~\ref{THF} give the decomposition of the Eynard--Orantin multidifferentials $\omega_{g,n} \in H^0(S^n,K_{S}(*\mathfrak{r})^{\boxtimes n})^{\mathfrak{S}_{n}}$ on the basis \eqref{xidkbasis}, that is for $2g - 2 + n > 0$
\begin{equation}
\label{omigyngufybgub}
\omega_{g,n}(z_1,\ldots,z_n) = \sum_{\substack{d_1,\ldots,d_n \geq 0 \\ r_1,\ldots,r_n}} F_{g,n}(\xi^*_{r_1,d_1} \otimes \cdots  \otimes \xi^*_{r_n,d_n}) \otimes_{i = 1}^n \xi_{r_i,d_i}(z_i),
\end{equation}
and the sum truncates to $d_1 + \cdots + d_n \leq 3g - 3 + n$.

\subsection{Fibering over Teichm\"uller space}
\label{sfiber}
\label{FIBERING}

Imagine that $\tilde{E}$ is a pre-target theory, \textit{i.e.} a target theory without the data of length functions. The typical case we have in mind are functors coming from topological field theories or conformal field theories. We can still construct a target theory by considering functions from Teichm\"uller space with values in $\tilde{E}$,
\[
E(\Sigma) = \mathscr{F}(\mathcal{T}_{\Sigma},\tilde{E}(\Sigma)).
\]
We define the union and glueing morphisms for $E$ by combining those of $\tilde{E}$ and those of Sections~\ref{TargeTeich}. We equip $E(\Sigma)$ with the seminorms indexed by $i \in \tilde{\mathscr{I}}^{\Sigma}$, $\alpha \in \tilde{\mathscr{A}}^{(i)}_{\Sigma}$, $\varepsilon \in (0,1)$ and $\sigma \in \mathcal{T}_{\Sigma}^{\epsilon}$ defined by
\[
|f|_{(i,\epsilon),(\alpha,\sigma)} = \sup_{\sigma' \in F_{\epsilon;\sigma}} |f(\sigma')|_{i,\alpha}.
\]
We equip $E$ with the collection of length functions induced by hyperbolic lengths as in \eqref{lengthC0}
\[
\ell_{(i,\epsilon),(\alpha,\sigma)}(\gamma) = \ell_{\sigma}(\gamma),
\]
which in fact does not depend on $i,\epsilon,\alpha$.

If $(\tilde{A},\tilde{B},\tilde{C},\tilde{D})$ are initial data for $\tilde{E}$, we can turn it canonically to admissible initial data $(A,B,C,D)$ for $E$ with the formulae
\begin{equation}
\label{inifiber} 
\begin{split}
 A_{P}(L_1,L_2,L_3) & = \tilde{A}_{P} , \\
B_{P}^b(L_1,L_b,\ell) & = B^{{\rm M}}(L_1,L_2,\ell)\,\tilde{B}_{P} ,\\ 
C_{P}(L_1,\ell,\ell') & = C^{{\rm M}}(L_1,\ell,\ell')\,\tilde{C}_{P}, \\
D_{T}(\sigma) & = \tilde{D}_{T}.
\end{split}
\end{equation}
 If the mapping class groups act trivially on $\tilde{E}$, using the Mirzakhani-McShane identity, we deduce that the GR amplitudes for $E$ attached to \eqref{inifiber} will be constant functions on Teichm\"uller space. This constant value in $\tilde{E}(\Sigma)$ can be considered as a definition of the GR amplitude for $\tilde{E}$ with initial data $(\tilde{A},\tilde{B},\tilde{C},\tilde{D})$. An example will appear in Section~\ref{Chehwam}.

This trick circumvents the potential absence of length functions in pre-target theories. If the mapping class group act non-trivially on $\tilde{E}(\Sigma)$, this argument does not apply and the GR amplitudes attached to \eqref{inifiber} could \textit{a priori} be interesting functions on the Teichm\"uller space.

 \subsection{(Lack of) symmetry}
 
 If $\Omega_{\Sigma}$ is obtained by symmetric GR  (Section~\ref{SGRSym}), then it is invariant under braidings of all boundary components of $\Sigma$ and the corresponding functions $V\Omega_{g,n}(L_1,\ldots,L_n)$ are symmetric in their $n$ variables. Otherwise, the GR amplitudes $\Omega_{\Sigma}$ are not a priori invariant under braidings of $\partial_-\Sigma$ with other boundary components and $V\Omega_{g,n}(L_1,\ldots,L_n)$ are only invariant under permutation of $L_2,\ldots,L_n$. Sufficient conditions to ensure that $V\Omega_{g,n}$ is symmetric in all variables already appear in a slightly different setting, in Section~\ref{SSSS}. Here they amount to requiring that $A$ is symmetric in its three variables, and the following four conditions to hold for any $L_1,L_2,L_3,L_4 > 0$
\begin{equation}
\label{theuytbsugfngoig}\begin{split}
\int_{\mathbb{R}_{+}}\big(B(L_1,L_2,\ell)A(\ell,L_3,L_4) + B(L_1,L_3,\ell)A(L_2,\ell,L_4) + B(L_1,L_4,\ell)A(L_2,\ell,L_3)\big) \ell\, \dd \ell & = (L_1 \leftrightarrow L_2),   \\
\int_{\mathbb{R}_{+}} \big(B(L_1,L_2,\ell)B(\ell,L_3,L_4) + B(L_1,L_3,\ell)B(L_2,\ell,L_4) + C(L_1,L_4,\ell)C(L_2,\ell,L_3)\big) \ell\, \dd \ell & = (L_1 \leftrightarrow L_2),  \\
\int_{\mathbb{R}_{+}} \big(B(L_1,L_2,\ell)C(\ell,L_3,L_4) + C(L_1,L_3,\ell)B(L_2,\ell,L_4) + C(L_1,L_4,\ell)B(L_2,\ell,L_3)\big) \ell\, \dd \ell  & = (L_1 \leftrightarrow L_2),  \\
\int_{\mathbb{R}_{+}} B(L_1,L_2,\ell)\cdot V\!D(\ell)\ell\, \dd \ell  + \frac{1}{2} \int_{\mathbb{R}_{+}^2} C(L_1,\ell,\ell')A(L_2,\ell,\ell') \ell \ell'\,\dd\ell\,\dd\ell'& = (L_1 \leftrightarrow L_2) .
\end{split}
\end{equation}
It is possible -- see \cite{WKarticle} -- to have geometrically meaningful GR amplitudes $\Omega_{\Sigma}$ which are not invariant under all braidings, but for which the integrals $V\Omega_{g,n}$ are symmetric functions of all the length variables.

\subsection{Coupling to 2d TQFTs}

Let $(\mathcal{A},\cdot,\langle \, , \, \rangle)$ be a Frobenius algebra, together with an (arbitrary) choice of hermitian norm $|\cdot|_{\mathcal{A}}$. The previous construction can be adapted to 
\begin{equation}
\label{E2d} E(\Sigma) = \mathscr{F}\big(\mathcal{T}_{\Sigma},\mathcal{A}^{\otimes \pi_0(\partial \Sigma)}\big),
\end{equation}
where the glueing morphisms now incorporate the pairing in $\mathcal{A}$, and the seminorms of these spaces are provided by the supremum of $|\cdot|_{\mathcal{A}}$. The integration of GR amplitudes $\Omega_{\Sigma}$ then produce functions $V\Omega_{g,n}$ from $\mathbb{R}_{+}^n$ to $\mathcal{A}^{\otimes n}$.

Another way to describe \eqref{E2d} is to make the tensor product of the strict target theory based on the 2d TQFT of $\mathcal{A}$, with the aforementioned target theories of functions over Teichm\"uller spaces. This tensor product is well-defined because $\mathcal{A}$ is a finite-dimensional vector space. When $\mathcal{A}$ is semi-simple, the choice of a canonical basis on $\mathcal{A}$ gives a natural choice of hermitian norm on $\mathcal{A}$. This is a basic example of the fibering procedure described later in Section~\ref{sfiber}.

\subsection{Various ways to compute with Corollary~\ref{costatunfg}}
\label{GUNIGUNFIOGSO}

As it is the most useful for applications, let us comment further on the case of Corollary~\ref{costatunfg} when $A \in \mathbb{R}[L_1^2,L_2^2,L_3^2]$, $V\!D \in \mathbb{R}[L_1^2]$ and the operators
\begin{equation*}
\begin{split}
\hat{B}[\phi_1](L_1,L_2) & = \int_{\mathbb{R}_{+}} B(L_1,L_2,\ell)\phi_1(\ell)\ell\,\dd\ell,  \\
\hat{C}[\phi_2](L_1) & = \int_{\mathbb{R}_{+}^2} C(L_1,\ell,\ell')\phi_2(\ell,\ell')\ell\,\ell'\,\dd\ell\dd\ell'
\end{split}
\end{equation*}
send even polynomials to even polynomials (from one variable to two variables for $\hat{B}$, and from two variables  to one variable variable for $\hat{C}$). We then say the initial data is of polynomial type. This occurs for instance for Mirzakhani initial data \eqref{iniMirza} -- see \cite[Section 5]{Mirzaint} or \cite{TRlecture} -- and for the Kontsevich initial data \eqref{Kontinfigufn}. Being of polynomial type is preserved by the twisting operation \eqref{twsitini}. The nature of this operation becomes more transparent if we decompose the initial data on a basis of even polynomials. Namely, if we introduce the scalars $\mathsf{A}^{i}_{j,k}, \mathsf{B}^{i}_{j,k}, \mathsf{C}^{i}_{j,k}$ and $\mathsf{D}^{i}$ by the formulas
\begin{equation*}
\begin{split}
A(L_1,L_2,L_3) & =  \sum_{i,j,k \geq 0} \mathsf{A}^{i}_{j,k}\,L_{1}^{2i} L_2^{2j} L_3^{2k}, \\
\hat{B}[\ell^{2k}] & =  \sum_{i,j \geq 0} \mathsf{B}^{i}_{j,k}\,L_1^{2i}L_2^{2j} ,\\
\hat{C}[\ell^{2j} (\ell')^{2k}] & =  \sum_{i \geq 0} \mathsf{C}^{i}_{j,k}\,L_1^{2i} , \\
V\!D(L_1) & =  \sum_{i \geq 0} \mathsf{D}^{i}\,L_1^{2i},
\end{split}
\end{equation*}
we then have a similar decomposition  
\begin{equation*}
\label{thetwistng} \begin{split} 
\phantom{ss} \mathsf{A}[f] & = \mathsf{A}^{i}_{j,k}, \\ 
\phantom{ss} \mathsf{B}[f] & = \mathsf{B}^{i}_{j,k} + \sum_{a \geq 0} \mathsf{A}^{i}_{j,a}u_{a,k}, \\ 
\phantom{ss} \mathsf{C}[f]_{j,k}^{i} & =  \mathsf{C}^i_{j,k} + \sum_{a \geq 0} \mathsf{B}^i_{a,k} u_{a,k} + \mathsf{B}^i_{j,a} u_{a,j} + \sum_{a,b \geq 0} \mathsf{A}^{i}_{a,b}u_{j,a}u_{k,a}, \\
 \phantom{ss} \mathsf{D}[f]^{i} & =  \mathsf{D}^i + \frac{1}{2} \sum_{a,b \geq 0} \mathsf{A}^{i}_{a,b} u_{a,b},
 \end{split}
 \end{equation*}
where $u$ is computed in terms of the odd moments of $f$
\begin{equation}
\label{uyeah}
u_{i,j} = m_{2(i + j) + 1}[f],\qquad m_{k}[f] = \int_{\mathbb{R}_{+}} \ell^{k}\,f(\ell) \dd \ell.
\end{equation}

\vspace{0.2cm}

\noindent \textsc{Relation to quantum Airy structures --} Recall the terminology and notations of Section~\ref{Airystr}. Let us assume that the coefficients $(\mathsf{A},\mathsf{B},\mathsf{C},\mathsf{D})$ coming from some GR initial data of polynomial type in the target theory ${\rm Mes}(\mathcal{T}_{\Sigma},\mathcal{A}^{\otimes \pi_0(\partial \Sigma)})$, describe the quantum Airy structure coming from a spectral curve $(S,x,y,\omega_{0,2})$ with Jacobi ring $\mathcal{A}$. Since we have
\begin{equation}
\label{dnddddiugnfgungu}\frac{(2d + 1)!}{\zeta^{2d + 2}} = \int_{\mathbb{R}_{+}} \ell^{2d}\,e^{-\zeta \ell}  \ell\,\dd \ell,
\end{equation}
applying the Laplace transform to \eqref{omigyngufybgub} gives the following relation between Eynard--Orantin multidifferentials and the TR amplitudes $V\Omega_{g,n}$
\begin{equation}
\label{theLpalceun}\begin{split}
& \quad \omega_{g,n}(z_1,\ldots,z_n) \\
& = \sum_{r_1,\ldots,r_n \in \mathfrak{r}} \Res_{z_1' = r_1} \cdots \Res_{z_n' = r_n} \bigg[\prod_{i = 1}^n  \dd \zeta(z_i') \int_{r_i}^{z_i'} \omega_{0,2}(\cdot,z_i)\bigg] \int_{\mathbb{R}_{+}^n} \big(V\Omega_{g,n}(L) \,\big|\, e_{r_1}\otimes \cdots \otimes e_{r_n}\big)_{\mathcal{A}} \prod_{i = 1}^n L_i\,e^{-L_i\zeta(z_i')} \dd L_i
\end{split}
\end{equation}
where we used the pairing in $\mathcal{A}$ for each of the $n$ tensor factors.

We recognise by comparing \eqref{thetwistng} to \cite[Section 4.1]{TRABCD} or \cite[Section 1.5]{TRlecture} that the twisting by $f$ corresponds to a certain change of polarisation, therefore to a shift of $\omega_{0,2}$; from the precise relation between the $u$s and the choice of $\omega_{0,2}$ in \cite[Section 9]{TRABCD}, we find from \eqref{uyeah} that the twisting amounts to adding to $\omega_{0,2}$ the symmetric holomorphic bidifferential
\[
\bigg(\int_{\mathbb{R}_{+}} f(\ell) e^{-\ell(\zeta_1 + \zeta_2)} \ell\,\dd \ell\bigg)\dd \zeta_1\dd \zeta_2 \in {\rm Sym}^2(\mathcal{A}[[\zeta]])
\]
defined locally near $\mathfrak{r}^2 \subset S^2$.

\vspace{0.2cm}
 
\noindent \textsc{Relation with Givental group action --} The operation of twisting is in fact known in another context. Sums over stable graphs describe the action of the $R$ element of the symplectic loop group on correlation functions of cohomological field theories introduced by Givental \cite{FSZ,Givental}. There, one considers a Frobenius algebra $\mathcal{A}$ and a formal series $R \in {\rm End}(\mathcal{A})[[t]]$ such that $R(t) = {\rm id}_{\mathcal{A}} + O(t)$ and $R(t)R^{{\rm T}}(-t) = {\rm Id}$, and defines\footnote{The Givental group action is usually formulated in the basis such that the index $d$ is directly coupled to the insertion of $\psi^{d}$ without prefactor. As one can check by comparing with Theorem~\ref{IntPsi}, this basis is $\tfrac{\ell^{2d}}{2^{d}d!}$ in the $\ell$-variable, hence after Laplace transform $\frac{(2d + 1)!!\dd \zeta}{\zeta^{2d + 2}}$ in the $\zeta$-variable. This explains the $(2d + 1)!!$ factors appearing in \eqref{IRD}.} $u_{a,b} \in {\rm Sym}^2(\mathcal{A})$ by 
\begin{equation} 
\label{IRD} \bigg(\frac{{\rm Id} - R(t_1) \otimes R(t_2)}{t_1 + t_2}\bigg)[\eta] = \sum_{d_1,d_2 \geq 0} (2d_1 + 1)!!(2d_2 + 1)!!\,u_{d_1,d_2}\,t^{d_1}_1t^{d_2}_2,
\end{equation}
where $\eta \in \mathcal{A}^{\otimes 2}$ describes the pairing in $\mathcal{A}$. The $u$s coming from an $R$ matrix are not always of the form \eqref{uyeah}. However, if we only consider $\mathcal{A} = \mathbb{C}$, we can make the following observation.

\begin{lemma}
\label{Ltssshet}Assume there exists a measurable test function $f\,\colon\,\mathbb{R}_{+} \rightarrow \mathbb{C}$ and $R \in \mathbb{C}[[t]]$ such that $u_{d_1,d_2}$ satisfy both \eqref{IRD} and \eqref{uyeah} and is not identically $0$. Then, there exists a unique $H \in \mathbb{R}_{+}$ such that, up to adding to $f$ a test function whose odd moments vanish, we have
\[
f(\ell) = -\eta\,\Theta(H^{1/2} - \ell),\qquad R(t) = e^{Ht/2},
\]
\end{lemma}
\noindent \textbf{Proof.}  We compute from \eqref{uyeah}
\[
\sum_{d_1,d_2 \geq 0} (2d_1 + 1)!!(2d_2 + 1)!!\,u_{d_1,d_2}\,t_1^{d_1}t_2^{d_2} = \int_{\mathbb{R}_{+}} f(\ell)\,e^{\ell (t_1 + t_2)/2}\,\ell\,\dd \ell.
\]
In particular, this generating series depends only on $t_1 + t_2$. Taking into account $R(0) = 1$, this is only possible when $R(t_1)R(t_2) = R(t_1 + t_2)$, that is $R(t) = e^{Ht/2}$ for some $H \in \mathbb{C}$. When this is the case, we compute with \eqref{IRD}
\[
(2d_1 + 1)!!(2d_2 + 1)!!\,u_{d_1,d_2} = -\eta\,\frac{H^{d_1 + d_2 + 1}}{2^{d_1 + d_2 + 1}(d_1 + d_2 + 1)}\,\frac{1}{d_1!d_2!}.
\]
Thus
\begin{equation}
\label{udd2} u_{d_1,d_2} = -\frac{\eta}{2(d_1 + d_2 + 1)}\,\frac{H^{d_1 + d_2 + 1}}{(2d_1 + 1)!(2d_2 + 1)!}.
\end{equation}
This is satisfied if we choose $H \geq 0$ and $f(\ell) = - \eta\,\Theta(H^{1/2} - \ell)$ and the case $H = 0$ gives vanishing $u$. If $\hat{f}$ is another test function that yields \eqref{udd2}, then
\[
\forall k \geq 0,\qquad \int_{\mathbb{R}_{+}} \big(f(\ell) - \hat{f}(\ell)\big)\,\ell^{2k + 1}\,\dd \ell = 0.
\]
\hfill $\Box$

\subsection{Chern character of the bundle of conformal blocks and length statistics}

\label{Chehwam}

Let $\mathcal{Z}$ is a modular functor. We denote $\Lambda$ its associated label set, $\dagger$ the involution on $\Lambda$, $\mathcal{A}$ its associated Frobenius algebra with its natural basis $(e_{\lambda})_{\lambda \in \Lambda}$ and product $\mu$ encoding the Verlinde rules. Using the pairing to identify $\mathcal{A} \cong \mathcal{A}^*$, we shall consider $\mu$ as an element of $\mathcal{A}^{\otimes 3}$.

For any $g,n \geq 0$ such that $2g - 2 + n > 0$ and $\vec{\lambda} \in \Lambda^{n}$, we have constructed in \cite{ABO1} a bundle $\mathcal{Z}_{g,n}(\vec{\lambda}) \rightarrow \overline{\mathfrak{M}}_{g,n}$ and shows that its Chern character defines a semi-simple cohomological field theory $w_{g,n}\,\colon\,\mathcal{A}^{\otimes n} \rightarrow H^{\bullet}(\overline{\mathfrak{M}}_{g,n})$ via the formula
\begin{equation}
\label{ifugfffffnfdugng} w(e_{\lambda_1} \otimes \ldots \otimes e_{\lambda_n}) = {\rm Ch}(\mathcal{Z}_{g,n}(\vec{\lambda})).
\end{equation}
Furthermore, this cohomological theory is identified as the result of the Givental group action for the $R$-operator\footnote{This formula differs by a sign of $r_{\lambda}$ compared to \cite{ABO1}. The reason is that in \cite{ABO1} the contribution of $\tilde{\mathcal{L}}_{p}$ to the Chern character was mistaken for $\psi_p$ instead of $-\psi_p$, which should be corrected with this sign.}\begin{equation}
\label{thetRT}
R(t) = e^{t(c/24 - L_0)},
\end{equation}
where $L_0 \in {\rm End}(\mathcal{A})$ is the operator whose eigenvectors are $e_{\lambda}$ with eigenvalue given by the conformal dimensions $r_{\lambda}$, and $c$ is the central charge. We refer to \cite{ABO1} for details. 

Here we find a rather surprising interpretation of the intersection indices of this cohomological field theories with $\psi$-classes in terms of integrals (with respect to the Weil--Petersson measure) of statistics of hyperbolic lengths of multicurves by comparing with Corollary~\ref{costatunfg}.

\begin{proposition}
\label{chengunu}Assume that for any $\lambda \in \Lambda$ we have $c/24 - r_{\lambda} > 0$. Consider the following initial data for the target theory $E(\Sigma) = {\rm Mes}(\mathcal{T}_{\Sigma},\mathcal{A}^{\otimes \pi_0(\Sigma)})$:
\begin{equation}
\label{songifgng}
\begin{split}
A^{\mathcal{Z}}(L_1,L_2,L_3) & = \mu, \\
B^{\mathcal{Z}}(L_1,L_2,\ell) & = \bigg(\frac{2{\rm i}\pi}{\sqrt{c/12}}\bigg)^2\,B^{{\rm M}}(L_1,L_2,\ell)\,\mu, \\
C^{\mathcal{Z}}(L_1,\ell,\ell') & = \bigg(\frac{2{\rm i}\pi}{\sqrt{c/12}}\bigg)^4\,C^{{\rm M}}(L_1,\ell,\ell')\,\mu,  \\
D^{\mathcal{Z}}(\sigma) & =  \bigg(\frac{2{\rm i}\pi}{\sqrt{c/12}}\bigg)^{2} \sum_{\lambda \in \Lambda} {\rm rk}(\mathcal{Z}_{1,1}(\lambda))\,e_{\lambda},
\end{split}
\end{equation}
and set
\begin{equation}
\label{fchoicefh} f(\ell) := - \sum_{\lambda \in \Lambda} \Theta(H_{\lambda}- \ell)\,e_{\lambda} \otimes e_{\lambda^{\dagger}},\qquad {\rm where}\,\,\, H_{\lambda} := \sqrt{2(c/24 - r_{\lambda})} .
\end{equation}
We denote $\Omega^{\mathcal{Z}}[f]$ the GR amplitude for the initial data \eqref{songifgng} twisted by $f$. We have for $2g - 2 + n > 0$ and $\vec{\lambda} \in \Lambda^n$
\begin{equation}
\label{yepur} (-1)^{n} \int_{\overline{\mathfrak{M}}_{g,n}} {\rm Ch}[\mathcal{Z}_{g,n}(\vec{\lambda})] \exp\bigg(\sum_{i = 1}^n \frac{L_i^2}{2}\,\psi_i\bigg)
 = \Big(V\Omega_{g,n}^{\mathcal{Z}}[f]\big(\tfrac{\sqrt{c/12}}{2{\rm i}\pi} L\big)\,\Big|\,e_{\lambda_1} \otimes \cdots \otimes e_{\lambda_n}\Big)_{\mathcal{A}},
 \end{equation}
where in the left-hand side, it is understood that the polynomial function $V\Omega^{\mathcal{Z}}[f](L)$ is specialised to $L' =\frac{\sqrt{c/12}}{2{\rm i}\pi}\,L \in \mathbb{C}^n$.
\end{proposition}

When the positivity conditions are not satisfied, we can still conclude that the right-hand side of \eqref{yepur} is the analytic continuation of a function $V\Omega_{g,n}^{\mathcal{Z}}[f]$ in the variables $(H_{\lambda})_{\lambda \in \Lambda}$, which coincide in the region $H_{\lambda} > 0$ with the integral of statistics of hyperbolic lengths. Notice in any case that \eqref{yepur} involves the specialisation of Weil--Petersson volumes to purely imaginary lengths. This is sometimes related to surfaces with conic singularities instead of boundaries, see e.g. \cite{NorburyDoconic}. Besides, the operator $L_0 - \frac{c}{24}$ gives the hamiltonian in conformal field theory, while $\sqrt{L_0 - \frac{(c - 1)}{24}}$ -- which would appear in \eqref{fchoicefh} if one had tensored $\mathcal{Z}_{g,n}(\vec{\lambda})$ with a squareroot of the Hodge bundle -- is the momentum. It is therefore tempting to ask for a physical explanation of Proposition~\ref{chengunu} in the conformal field theory corresponding to $\mathcal{Z}$ on surfaces of genus $g$ with $n$ conic singularities.

\vspace{0.2cm}

\noindent \textbf{Proof of Proposition~\ref{chengunu}.} As reviewed in \cite{ABO1}, the Frobenius algebra $\mathcal{A}$ of a modular functor admits a basis $(\underline{\epsilon}_{\lambda})_{\lambda \in \Lambda}$ which is orthogonal for the pairing and such that $\underline{\epsilon}_{\lambda} \times \underline{\epsilon}_{\mu} = \delta_{\lambda,\mu} \underline{\epsilon}_{\lambda} $ for the product. The $S$-matrix of $\mathcal{Z}$ gives the change of basis
\[ 
\underline{\epsilon}_{\lambda} = \sum_{\mu \in \Lambda} \underline{S}_{\lambda,\mu}\,e_{\mu},
\] 
and the Verlinde formula computes the ranks
\begin{equation} 
\label{Verlindesss} {\rm rk}[\mathcal{Z}_{g,n}(\vec{\lambda})] = \sum_{\mu \in \Lambda} \frac{\prod_{i = 1}^n \underline{S}_{\lambda_i,\mu}^{-1}}{(\underline{S}_{1,\mu}^{-1})^{2g - 2 + n}}.
\end{equation}

We consider the Eynard--Orantin multidifferentials $\omega_{g,n}^0$ for the local spectral curve $S = \mathbb{C}$ with
\begin{equation}
\label{thspcssssruve} x(\zeta) = \zeta^2/2,\qquad y(\zeta) = \sum_{\lambda \in \Lambda} -\underline{S}_{1,\lambda}^{-1}\,\frac{{\rm sinh}\big(\zeta \sqrt{c/12}\big)}{\sqrt{c/12}},\qquad \omega_{0,2}(\zeta_1,\zeta_2) = \sum_{\lambda \in \Lambda} \frac{\dd \zeta_1 \otimes \dd \zeta_2}{(\zeta_1 - \zeta_2)^2}\,\underline{\epsilon}_{\lambda} \otimes \underline{\epsilon}_{\lambda}.
\end{equation}
By comparison with Theorem~\ref{iuniuufgbdfginugiun}, we find
\[
\frac{\omega_{g,n}^0(\zeta_1,\ldots,\zeta_n)}{\dd \zeta_1 \otimes \cdots \otimes \dd \zeta_n} = \sum_{\lambda \in \Lambda} \bigg(-\frac{(2{\rm i}\pi)^3}{(c/12)^{3/2}}\,\underline{S}_{1,\lambda}^{-1}\bigg)^{2 - 2g - n}\,\bigg(\frac{2{\rm i}\pi}{\sqrt{c/12}}\bigg)^{n} \int_{\mathbb{R}_{+}^n} V\Omega_{g,n}^{{\rm M}}\big(\tfrac{\sqrt{c/12}}{2{\rm i}\pi}\,L\big)\,\underline{\epsilon}_{\lambda}^{\otimes n} \prod_{i = 1}^n \,L_i\,e^{-\zeta_i L_i}\,\dd L_i
\] 
Hence, taking scalar products and using \eqref{Verlindesss}, we find
\begin{equation}
\label{jiniuynKKKK}
\begin{split}
& \quad \bigg(\frac{\omega_{g,n}^0(\zeta_1,\ldots,\zeta_n)}{\dd \zeta_1 \otimes \cdots \otimes \dd\zeta_n}\,\Big|\,e_{\lambda_1} \otimes \cdots \otimes e_{\lambda_n}\bigg)_{\mathcal{A}} \\
& = (-1)^{n} \bigg(\frac{2{\rm i}\pi}{\sqrt{c/12}}\bigg)^{6 - 6g - 2n}\,{\rm rk}[\mathcal{Z}_{g,n}(\vec{\lambda})]\,\int_{\mathbb{R}_{+}^n} V\Omega_{g,n}^{{\rm M}}\big(\tfrac{\sqrt{c/12}}{2{\rm i}\pi}\,L\big) \prod_{i = 1}^n \,L_i\,e^{-\zeta_i L_i}\,\dd L_i.
\end{split}
\end{equation}
Using Mirzakhani--McShane identities (Theorem~\ref{thMirza}), Verlinde formula and tracking the powers of $\tfrac{2{\rm i}\pi}{\sqrt{c/12}}$ we also find that the integral of the GR amplitudes for the initial data \eqref{songifgng} are given by
\[
\big(\Omega_{g,n}^{\mathcal{Z}}(\sigma)\,\big|\,e_{\lambda_1} \otimes \cdots \otimes e_{\lambda_n}\big)_{\mathcal{A}} = \bigg(\frac{2{\rm i}\pi}{\sqrt{c/12}}\bigg)^{6g - 6 + 2n}\,{\rm rk}[\mathcal{Z}_{g,n}(\vec{\lambda})]\,\Omega_{\Sigma}^{{\rm M}}(\sigma),
\]
and therefore reproduce up to prefactors the right-hand side of \eqref{jiniuynKKKK} after integration over $\mathcal{M}_{g,n}(L')$ and specialisation of the result (which is a polynomial in $L'$) to the values $L' = \frac{\sqrt{c/12}}{2{\rm i}\pi} L$.

The $R$-operator \eqref{thetRT} is of the form considered in Lemma~\ref{Ltssshet}, so it realises simultaneously \eqref{IRD} and \eqref{uyeah} for $f$ given by \eqref{fchoicefh}. Theorem~\cite[Theorem 4.5]{ABO1} identifies the Laplace transform of the right-hand side of \eqref{yepur} with the Eynard--Orantin multidifferentials for the spectral curve obtained from \eqref{thspcssssruve} by twisting via this $R$-operator. We then use the correspondence \eqref{theLpalceun} and the fact that the Laplace transform \eqref{dnddddiugnfgungu} between odd polynomials in $\zeta^{-1}$ and even polynomials in $L$ is an isomorphism to conclude. \hfill $\Box$

\medskip

\section{A symmetric version of the geometric recursion}
\label{SGRSym}
\label{S5}
In the recursive definition \eqref{defGR}, the single boundary component in $\partial_-\Sigma$ plays a special role, and this propagates at each step of the recursion. This is the reason why we have used the category $\Surf_{1}$ instead of $\Surf$. We now describe a version of GR where all boundary components play the same role. In fact, we give sufficient conditions on the initial data for this symmetric GR to be well-defined.

\begin{definition}
Let $\Surf_{{\rm s}}$ be the full subcategory of $\mathcal{B}$ whose objects $\Sigma$ are such that $\partial_-\Sigma = \emptyset$.
\end{definition}

In this section, it is assumed that all surfaces are objects of $\Surf_{{\rm s}}$, unless stated otherwise. If $\Sigma$ is a connected object in $\Surf_{{\rm s}}$ and $b_1$ a choice of boundary component, we can form the object $\Sigma_{b_1}$ of $\Surf_1$ by reversing the orientation of $b_1$. The sets of homotopy classes of embedded pairs of pants $\mathcal{P}_{\Sigma_{b_1}}$ and $\mathcal{P}_{\Sigma_{b_1}}^{d}$ for $d \in \{\emptyset\} \cup \pi_0(\partial_{+}\Sigma_{b_1})$ described in Section~\ref{glueinpant} the type of excisions from the object $\Sigma_{b_1}$ in $\Surf_{1}$ are here denoted $\mathcal{P}_{\Sigma,b_1}$ and $\mathcal{P}_{\Sigma,b_1}^{d}$. If $[P] \in \mathcal{P}_{\Sigma,b_1}$,  the surface $\Sigma - P$ is a priori an object in $\Surf_{1}$. We automatically make it an object of $\Surf_{{\rm s}}$ by forcing the orientation of each boundary components to agree with the orientation of the surface.

The notion of symmetric target theory $E:\,\Surf_{{\rm s}} \rightarrow \mathcal{C}$ is defined as in Section~\ref{STarget}, by replacing $\Surf_1$ with $\Surf_{{\rm s}}$ and requiring functorial glueing morphisms denoted $\Theta_{P}$ and indexed by $b_1 \in \pi_0(\partial \Sigma)$ and $[P] \in \mathcal{P}_{\Sigma,b_1}$.

\subsection{Initial data}
\label{SGRSymIni}

\begin{definition}
\label{DEFGRSSSS} Initial data for a given symmetric target theory $E$ are functorial assignments
\begin{itemize} 
\item[$\bullet$] of $A_{P} \in E(P)$ for any pair of pants $P$.
\item[$\bullet$] of $B_{P}^{b_1,b_2} \in E(P)$ for any pair of pants $P$ for which an ordered pair $(b_1,b_2)$ of distinct boundary components has been selected.
\item[$\bullet$] of $C_{P}^{b_1} \in E(P)$ for any pair of pants $P$ for which some $b_1 \in \pi_0(\partial P)$ has been selected.
\item[$\bullet$] of $D_{T}\in E(T)$ for any torus $T$ with one boundary component.
\end{itemize} 
An initial data is called admissible if it satisfies the decay properties (as in Definition~\ref{DEFDECAY}) and the four symmetry properties stated below.
\end{definition} 

In the following relations, $X$ will denote any genus $0$ surface with $4$ ordered boundary components $(b_i)_{i = 1}^4$, and \foreignlanguage{russian}{Yu} any surface with genus $1$  with $2$ ordered boundary components $(b_1,b_2)$. If $[P] \in \mathcal{P}_{X,b_1}$, we recall that $\gamma_P$ is the unique boundary component of $P$ which is interior to $X$. If $(\beta,\beta')$ is an ordered pair of boundary components of $X$ or \foreignlanguage{russian}{Yu}, we denote by $\sigma_{\beta,\beta'}$ any mapping class that exchanges $\beta$ and $\beta'$. It is called a braiding of $\beta$ and $\beta'$. 

\vspace{0.2cm}

\noindent\textsc{BA relation.}
\begin{equation}
\label{BA} 
(\sigma_{b_1,b_2}-{\rm Id})\bigg[\sum_{[P] \in \mathcal{P}_{X,b_1}^{b_2}} ({\rm Id} + \sigma_{b_2,b_3} + \sigma_{b_2,b_4})\Theta_{P}(B_{P}^{b_1,b_2},A_{X - P})\bigg] = 0\,.
\end{equation}

\noindent\textsc{BB-CA relation.}

\begin{equation}
\label{BBAC} 
(\sigma_{b_1,b_2}-{\rm Id})\bigg[\sum_{[P] \in \mathcal{P}_{X,b_1}^{b_2}}\!\!\!\Theta_{P}(B^{b_1,b_2}_{P},B_{X - P}^{\gamma_P,b_3}) + \!\!\!\sum_{[P] \in \mathcal{P}_{X,b_1}^{b_3} }\Theta_{P}(B^{b_1,b_3}_{P},B^{b_2,\gamma_P}_{X - P}) + \!\!\! \sum_{[P] \in \mathcal{P}_{X,b_1}^{b_4}}\!\!\!\Theta_{P}(C^{b_1}_{P},A_{X - P})\bigg] =0\,. 
\end{equation}

\noindent \textsc{BC relation.}
\begin{equation}
\label{BC} 
(\sigma_{b_1,b_2}-{\rm Id})\bigg[\sum_{[P] \in \mathcal{P}_{X,b_1}^{b_2}}\Theta_{P}(B^{b_1,b_2}_{P},C^{\gamma_P}_{X - P}) + \sum_{[P] \in \mathcal{P}_{X,b_1}^{b_3}} \Theta_{P}(C^{b_1}_{P},B^{b_2,\gamma_P}_{X - P}) + \sum_{[P] \in \mathcal{P}_{X,b_1}^{b_4}} \Theta_{P}(C^{b_1}_{P},B^{b_2,\gamma_P}_{X - P}) \bigg] =0.
\end{equation}

\noindent \textsc{D relation.}
\begin{equation}
\label{Drel} 
(\sigma_{b_1,b_2}-{\rm Id})\bigg[\sum_{[P] \in \mathcal{P}_{\text{\foreignlanguage{russian}{Yu}},b_1}^{b_2}} \Theta_{P}(B^{b_1,b_2}_{P},D_{\text{\foreignlanguage{russian}{Yu}} - P}) + \tfrac{1}{2} \sum_{[P] \in \mathcal{P}_{\text{\foreignlanguage{russian}{Yu}},b_1}^{\emptyset}} \Theta_{P}(C^{b_1}_{P},A_{\text{\foreignlanguage{russian}{Yu}} - P})\bigg] = 0.
\end{equation}
\vspace{0.2cm}

We note that, due to the properties of the initial data, the validity of these relations for one choice of braidings $\sigma_{b_1,b_2}, \sigma_{b_2,b_3}$ and $\sigma_{b_3,b_4}$ imply their validities for all choice of braidings.

\vspace{0.7cm}

\begin{figure}[ht!]
\begin{center}
\includegraphics[width=0.95\textwidth]{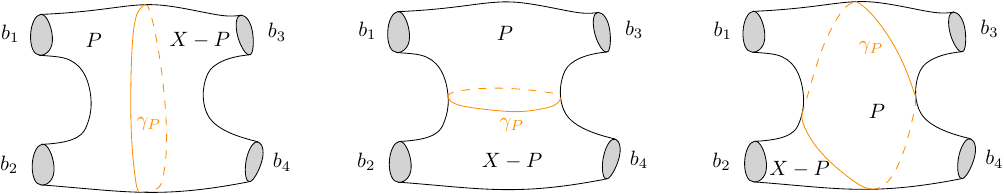}
\caption{\label{SREL1FIG} The three types of terms in the BB-CA relation. The same three types of terms appear in the BA and BC relations. These relations take an H-I-X form, and provide a geometric incarnation for the (purely algebraic) H-I-X relations found in quantum Airy structures and the topological recursion, see Section~\ref{Airystr}.}
\end{center}
\end{figure}

\vspace{0.5cm}

\begin{figure}[ht!]
\begin{center}
\includegraphics[width=0.55\textwidth]{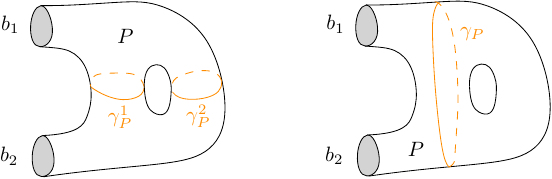}
\caption{\label{SREL2FIG} The two types of terms in the D-relation.}
\end{center}
\end{figure}

\subsection{Definition}

Let $(A,B,C,D)$ be admissible initial data for a symmetric target theory $E:\,\Surf_{{\rm s}} \rightarrow \mathcal{C}$.

\begin{definition}
We define the symmetric GR amplitudes as follows.

\begin{itemize}
\item[$\bullet$] We put $\Omega_{\emptyset} := 1 \in E(\emptyset)$.
\item[$\bullet$] If $P$ is a pair of pants in $\Surf_{{\rm s}}$, we put $\Omega_{P} = A_{P}$.
\item[$\bullet$] If $T$ is a torus with one boundary in $\Surf_{{\rm s}}$, we put $\Omega_{T} = D_{T}$.
\item[$\bullet$] For disconnected surfaces in $\Surf_{{\rm S}}$, we declare $\Omega_{\Sigma} = \bigsqcup_{a \in \pi_0(\Sigma)} \Omega_{\Sigma(a)}$.
\item[$\bullet$] If $\Sigma$ is a connected surface in $\Surf_{{\rm s}}$ with $\chi_{\Sigma} \leq -2$, let us make a choice of a boundary component $b_1$. We seek to inductively define as in \eqref{defGR}
\begin{equation}
\label{defGRSym} \Omega_{\Sigma} = \sum_{b \in \pi_0(\partial \Sigma) \setminus \{b_1\}} \sum_{[P] \in \mathcal{P}_{\Sigma,b_1}^{b}} \Theta_{P}(B_{P}^{b_1,b},\Omega_{\Sigma - P}) + \tfrac{1}{2} \sum_{[P] \in \mathcal{P}_{\Sigma,b_1}^{\emptyset}} \Theta_{P}(C_{P}^{b_1},\Omega_{\Sigma - P})
\end{equation}
as an element of $E(\Sigma)$.
\end{itemize}
\end{definition}

We now have to check that, step by step, this assignment does not depend on the choice of $b_1$.

\begin{theorem}
\label{SGRSYMTH} The assignment $\Sigma \mapsto \Omega_{\Sigma}$ is well-defined for any surface in $\Surf_{{\rm s}}$. More precisely, the series \eqref{defGRSym} converges absolutely for any of the semi-norms $|\cdot|_{i,\alpha}$ to a unique limit. This limit is an element of ${}^{\flat} E(\Sigma)$ which is functorial. In particular $\Omega_{\Sigma}$ is invariant under braidings of all boundary components.
\end{theorem}

\noindent \textbf{Proof.} As the convergence and functoriality under mapping classes that respect $b_1$ can be adressed as in Theorem~\ref{welldefGR}, we focus on the justification -- at each step of the induction -- of the invariance under the braiding of $b_1$ with another component.
We denote $\Omega_{\Sigma,b_1}$ the right-hand side of \eqref{defGRSym}, to stress its potential dependence in the choice of the boundary component. As the result for connected surfaces automatically implies the result for all surfaces,  we assume in the rest of the proof that $\Sigma$ is connected. The result obviously holds when $\chi_{\Sigma} = -1$.

\vspace{0.2cm}

\noindent\textsc{Step 1.}  We first study the case of genus $0$ with $4$ boundaries, that is $\Sigma = X$. The functoriality of the glueing morphisms imply that $\sigma_{b_1,b_2}\Omega_{X,b_1} = \Omega_{X,b_2}$. The \textbf{BA} relation therefore means that $\Omega_{X,b_1} = \Omega_{X,b_2}$. So, $\Omega_{X} := \Omega_{X,b_1}$ is a well-defined, functorial assignment in $E(X)$.

\vspace{0.2cm} 

\noindent \textsc{Step 2.} We prove the result likewise for a surface \foreignlanguage{russian}{Yu} with genus $1$ with $2$ boundary components, which we arbitrarily order $(b_1,b_2)$. Applying \eqref{defGRSym} with chosen boundary component $b_1$ yields
\[
\Omega_{\textrm{\foreignlanguage{russian}{Yu}},b_1} := \sum_{[P] \in \mathcal{P}_{\textrm{\foreignlanguage{russian}{Yu}},b_1}^{b_2}} \Theta_{P}(B^{b_1,b_2}_{P},D_{\textrm{\foreignlanguage{russian}{Yu}} - P}) + \tfrac{1}{2} \sum_{c \in \mathcal{P}_{\textrm{\foreignlanguage{russian}{Yu}},b_1}^{\emptyset}} \Theta_{P}(C^{b_1}_{P},A_{\textrm{\foreignlanguage{russian}{Yu}} - P})
\]
and this is invariant under braidings $\sigma_{b_1,b_2}$ according to relation \textsc{D}.

\vspace{0.2cm}
 
\noindent\textsc{Step 3.} So far we have proved the result up to $\chi_{\Sigma} \geq -2$. Assume it holds for surfaces with Euler characteristic larger than or equal to a $\chi_0 \leq -2$, and consider $\Sigma$ of Euler characteristic $\chi_0 - 1$.
 
If $\Sigma$ has only $1$ boundary component, then $\Omega_{\Sigma}$, defined by \eqref{defGRSym}, does not depend on any choice, and functoriality from the induction hypothesis, of the initial data and glueing maps imply that $\Omega_{\Sigma}$ is functorial for such surfaces. If $\Sigma$ has $n \geq 2$ boundary components, we order them arbitrarily $(b_1,\ldots,b_n)$ and we shall prove invariance of $\Omega_{\Sigma,b_1}$ under the braiding $\sigma_{b_1,b_2}$. Let us compute $\Omega_{\Sigma,b_1}$ from \eqref{defGRSym}, \textit{i.e.} excise pair of pants $P$ in all possible classes $[P] \in \mathcal{P}_{\Sigma,b_1}$.  We want to replace the contribution of some of the connected components of $\Sigma - P$ with the GR formula unambiguously defined from the previous induction steps. We distinguish several cases in doing so.
\begin{itemize}
\item[$\bullet$] When $[P] \in \mathcal{P}_{\Sigma,b_1}^{\emptyset}$,  $(\gamma_P^1,\gamma_P^2)$ are the connected components of the multicurve $\gamma_{P}$ along which we excised -- they are ordered by our definition of excision in \S~\ref{glueinpant}. We denote by $\Sigma_{P}$ the connected component of $\Sigma - P$, which contains $b_2$ and $\Sigma_P'$ the other connected component of $\Sigma - P$. If $\Sigma_{P}$ is a pair of pants bounded by $\gamma_P^k$ for some $k \in \{1,2\}$, $b_2$ and another $b_i$ for $i \geq 3$ we have that $\Omega_{\Sigma_{P}} = A_{\Sigma_P}$, and leave $\Omega_{\Sigma_{P}'}$ as it is. We denote by $\mathcal{P}_{\Sigma,b_1}^{\emptyset;b_i, k}$ the set of $[P]$s leading to this situation. We observe that the contributions for $k = 1$ and $2$ from $\mathcal{P}_{\Sigma,b_1}^{\emptyset;b_i,k}$ are equal, since $C^{b_1}_{P}$ is invariant under braidings of its two last boundaries. The subset $\mathcal{P}_{\Sigma,b_1}^{\emptyset;0}\subset \mathcal{P}_{\Sigma,b_1}^{\emptyset}$ will consist of those $[P]$s for which $\Sigma_{P}$ is not a pair of pants. For these we have that $\Omega_{\Sigma_{P}}$ is given by the GR formula \eqref{defGRSym} with $b_2$ the chosen boundary, \textit{i.e.} excise $Q$ from $\Sigma - P$ in all possible ways specified by $[Q] \in \mathcal{P}_{\Sigma - P,b_2}$. If $[Q] \in \mathcal{P}_{\Sigma - P,b_2}^{\emptyset}$ then it gives the contribution $C^{b_2}_{Q}$ and if $[Q] \in \mathcal{P}_{\Sigma - P,b_2}^{b_j}$, for $b_j$ a boundary component of $\Sigma - P$ distinct from $b_2$, we get $B^{b_2,b_j}_{Q}$.

\item[$\bullet$] When $c \in \mathcal{P}_{\Sigma,b_1}^{b_2}$, we have that $\Omega_{\Sigma - P}$ is given by the GR formula \eqref{defGRSym} with $\gamma_{P}$ as the chosen boundary component. We are therefore excising a second time with all possible $[Q] \in \mathcal{P}_{\Sigma - P,\gamma_{P}}$. If $[Q] \in \mathcal{P}_{\Sigma - P,\gamma_{P}}^{\emptyset}$, then we get a $C^{\gamma_P}_{Q}$, while if $[Q] \in \mathcal{P}_{\Sigma - P,\gamma_{P}}^{b_i}$ for some $i \geq 3$ we get a $B^{\gamma_P,b_i}_{Q}$.

\item[$\bullet$] When $c \in \mathcal{P}_{\Sigma,b_1}^{b_i}$ for $i \geq 3$, we have that $\Omega_{\Sigma - P}$ is given by the GR formula \eqref{defGRSym} with $b_2$ the chosen boundary component, \textit{i.e.} perform a second excision specified by $[Q] \in \mathcal{P}_{\Sigma - P,b_2}$.  If $[Q] \in \mathcal{P}_{\Sigma - P,b_2}^{\emptyset}$, we denote $(\gamma_{Q}^1,\gamma_{Q}^2)$ the ordered connected components of the multicurve $\gamma_{Q}$, and we get a $C^{b_2}_{Q}$. For $[Q] \in \mathcal{P}_{\Sigma - P,b_2}^{\gamma_{P}}$, we get $B^{b_2,\gamma_P}_{Q}$. In the remaining cases, $[Q] \in \mathcal{P}_{\Sigma - P,b_2}^{b_j}$ for some $j \geq 3$ with $j \neq i$, we get a $B^{b_2,b_j}_{Q}$.

\end{itemize}
Keeping in mind the constraint $\chi_{\Sigma} < -2$, this list exhausts all possible situations, in particular this process will not pull out contributions of a $D$. We also stress that $\Sigma'$ may be empty in some situations, in which case we use $\Omega_{\emptyset} = 1 \in \mathbb{K}$ and the union axiom of \S~\ref{STarget} for the emptyset. The result is
\begin{equation}
\label{RecSym}
\begin{split}
 & \quad \Omega_{\Sigma,b_1}  \\
& = \sum_{i = 3}^{n} \sum_{c \in \mathcal{P}_{\Sigma,b_1}^{\emptyset;b_i,1}} \Theta_{P}(C^{b_1}_{P},A_{S_P},\Omega_{S'_P})  \\
& \quad + \!\!\!\!\!\!\! \sum_{[P] \in \mathcal{P}_{\Sigma,b_1}^{\emptyset;0}} \!\sum_{[Q] \in \mathcal{P}_{\Sigma - P,b_2}^{\emptyset}} \!\!\!\!\!\!\!\!\! \tfrac{1}{4}\Theta_{P}\big(C^{b_1}_{P},\Theta_{Q}(C^{b_2}_{Q},\Omega_{\Sigma - P - Q})\big) + 
\sum_{i = 3}^{n} \sum_{[P] \in \mathcal{P}_{\Sigma,b_1}^{\emptyset;0}} \! \sum_{[Q] \in \mathcal{P}_{\Sigma - P,b_2}^{b_i}} \!\!\!\!\!\!\!\!\! \tfrac{1}{2}\Theta_{P}\big(C^{b_1}_{P},\Theta_{Q}(B^{b_2,b_i}_{Q},\Omega_{\Sigma - P - Q})\big)
\\
& \quad + \!\!\!\!\!\!\!\sum_{[P] \in \mathcal{P}_{\Sigma,b_1}^{\emptyset;0}} \bigg\{\sum_{[Q] \in \mathcal{P}_{\Sigma - P,b_2}^{\gamma_P^{1}}} \!\!\!\!\!\!\! \tfrac{1}{2}\Theta_{P}\big(C^{b_1}_{P},\Theta_{Q}(B^{b_2,\gamma_{P}^1}_{Q},\Omega_{\Sigma - P - Q})\big) +  \sum_{[Q] \in \mathcal{P}_{\Sigma - P,b_2}^{\gamma_P^{2}}} \!\!\!\!\!\!\! \tfrac{1}{2}\Theta_{P}\big(C^{b_1}_{P},\Theta_{Q}(B^{b_2,\gamma_{P}^2}_{Q},\Omega_{\Sigma - P - Q})\big)\bigg\} \\
& \quad + \!\!\!\!\!\!\! \sum_{[P] \in \mathcal{P}_{\Sigma,b_1}^{b_2}} \bigg\{\sum_{[Q] \in \mathcal{P}_{\Sigma - P,\gamma_{P}}^{\emptyset}}\!\!\!\!\!\!\!\!\! \tfrac{1}{2}\Theta_{P}\big(B^{b_1,b_2}_{P},\Theta_{Q}(C^{\gamma_{P}}_{Q},\Omega_{\Sigma - P - Q})\big) + \sum_{i = 3}^n \sum_{[Q] \in \mathcal{P}_{\Sigma - P,\gamma_{P}}^{b_i}} \!\!\!\!\!\!\! \Theta_{P}\big(B^{b_1,b_2}_{P},\Theta_{Q}(B^{\gamma_{P},b_i}_{Q},\Omega_{\Sigma - P - Q})\big)\bigg\}  \\  
& \quad + \!\!\sum_{i = 3}^n \sum_{[P] \in \mathcal{P}_{\Sigma,b_1}^{b_i}} \bigg\{\sum_{[Q] \in \mathcal{P}_{\Sigma - P,b_2}^{\emptyset}} \!\!\!\!\!\!\! \tfrac{1}{2} \Theta_{P}\big(B^{b_1,b_i}_{P},\Theta_{Q}(C^{b_2}_{Q},\Omega_{\Sigma - P - Q})\big) + \!\!\!\! \sum_{[Q] \in \mathcal{P}_{\Sigma - P,b_2}^{\gamma_P}} \!\!\!\!\!\!\! \Theta_{P}\big(B^{b_1,b_i}_{P},\Theta_{Q}(B^{b_2,\gamma_{P}}_{Q},\Omega_{\Sigma - P - Q})\big)  \\
& \quad \qquad\qquad + \sum_{\substack{j \geq 3 \\ j \neq i}} \sum_{[Q] \in \mathcal{P}_{\Sigma - P,b_2}^{b_j}} \!\!\!\!\!\!\! \Theta_{P}\big(B^{b_1,b_i}_{P},\Theta_{Q}(B^{b_2,b_j}_{Q},\Omega_{\Sigma - P - Q})\big)\bigg\}.
\end{split}
\end{equation}
We have ten multisums on the right, which we will denote $(\mathscr{S}_i)_{i = 1}^{10}$ respecting their order of appearance.

First we observe that $\mathscr{S}_2$ is a sum over the set of $([P],[Q]) \in \mathcal{P}_{\Sigma,b_1}^{\emptyset} \times \mathcal{P}_{\Sigma,b_2}^{\emptyset}$, such that there exist non-intersecting representatives $\gamma_P$ and $\gamma_{Q}$ associated to $P$ and $Q$ respectively and cutting along $\gamma_P \cup \gamma_{Q}$ will not produce a component which is a pair of pants containing as boundaries $b_2$ and some other $b_i$ for some $i \geq 3$. This sum is manifestly invariant under $\sigma_{b_1,b_2}$. The same kind of argument also applies to $\mathscr{S}_{10}$. 

Inspection of $\mathscr{S}_3+ \mathscr{S}_{8}$ reveals that the sum of these two terms must be invariant under $\sigma_{b_1,b_2}$, since we can use composition of glueings on one multicurve at a time.

\newpage

\begin{figure}[ht!]
\begin{center}
\includegraphics[width=0.95\textwidth]{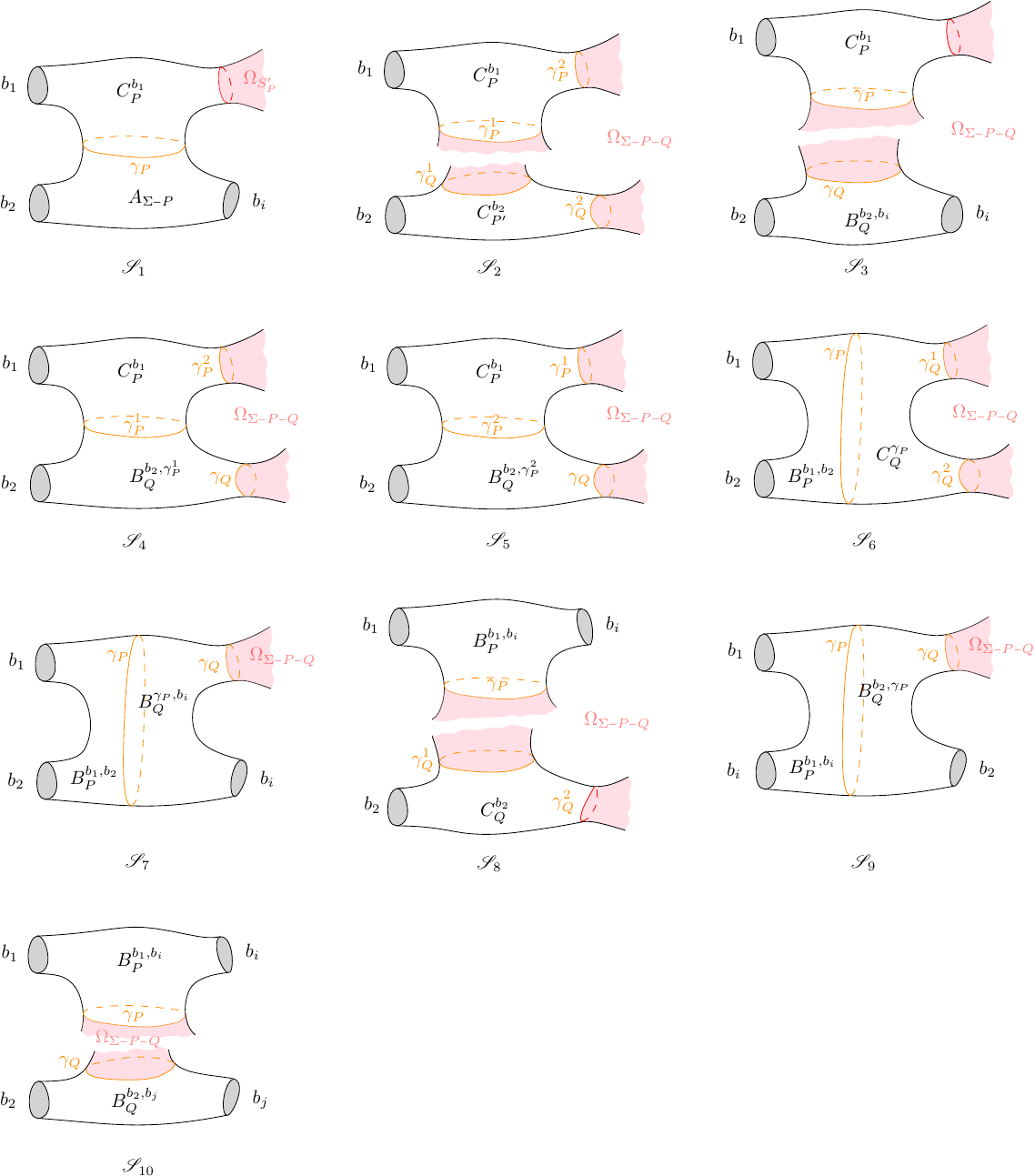}
\caption{\label{AllSFig} Geometric origin of the ten terms (respecting their order) in \eqref{RecSym}.}
\end{center}
\end{figure}

We now consider $\mathscr{S}_9+\mathscr{S}_7+\mathscr{S}_1$. First we invoke composition of glueing for the summands in $\mathscr{S}_9$ and $\mathscr{S}_7$, so as to reduce each summand to one glueing along one multicurve. If we now compare the summands of $\mathscr{S}_9$, $\mathscr{S}_{7}$ and $\mathscr{S}_{1}$ with the three terms in the \textsc{BB-CA} relation, we find they match. We now reorganise the three multisums of $\mathscr{S}_{9}$, $\mathscr{S}_{7}$ and $\mathscr{S}_{1}$, into a sum over all possible isotopy classes of embeddings of a surface $X$ of genus $0$ and $4$ boundary components, where three of its boundary components go to $b_1,b_2$ and $b_i$, for some fixed $i \geq 3$, followed by a sum as in the \textsc{BB-CA} relation locally in each of these embedded surfaces $X$. For each embedding of $X$ we can now invoke the \textsc{BB-CA} relation to conclude that all together $\mathscr{S}_9+ \mathscr{S}_7+ \mathscr{S}_1$ is invariant under $\sigma_{b_1,b_2}$.

Finally, we consider $\mathscr{S}_4+\mathscr{S}_5+\mathscr{S}_6$. We reorganise these three sums to first be a sum over all isotopy classes of embeddings of a surface $X$ of genus zero and four boundary components, where two of its boundary components goes to $b_1,b_2$ and the other two are mapped to the interior of $\Sigma$, followed by a sum as in the \textsc{BC} relation in each of these embedded surfaces $X$. Examining the summands of $\mathscr{S}_4,\mathscr{S}_5$ and $\mathscr{S}_6$ and using associativity of glueings, we observe that for each embedding of $X$, we can now invoke the \textsc{BC} relation to conclude all together that  $\mathscr{S}_4+\mathscr{S}_5+\mathscr{S}_6$ is invariant under $\sigma_{b_1,b_2}$.

\begin{remark}\label{remark-symmetry}
It is easy to see that this proof is still valid if the BB-CA relation is replaced by the relation 
\begin{equation} 
\label{BBCAprime} 
(\sigma_{b_1,b_2}-{\rm Id})\left(\sum_{[P] \in \mathcal{P}_{\Sigma,b_1}^{b_2}}\Theta_{P}(B^{b_1,b_2}_{P},B^{b_3,\gamma_P}_{\Sigma - P})  + \!\! \sum_{[P] \in \mathcal{P}_{\Sigma,b_1}^{b_3}}\Theta_{P}(B^{b_1,b_3}_{P},B^{b_2,\gamma_P}_{\Sigma - P}) + \!\!   \sum_{[P] \in \mathcal{P}_{\Sigma,b_1}^{b_4}}\Theta_{P}(C^{b_1}_{P},A_{\Sigma - P})\right) =0,
\end{equation}
obtained by replacing $B^{\gamma_P,b_3}_{\Sigma - P}$ by $B^{b_3,\gamma_P}_{\Sigma - P}$ in the BB-CA relation. This correspond to choosing $b_3$ instead of $\gamma_P$ as the reference boundary in the step 3 of the proof when $[P] \in \mathcal{P}_{\Sigma,b_1}^{b_2}$. By the induction hypothesis, this choice does not affect the result.
\end{remark}

\subsection{Main properties}
\label{S33}
There is a forgetful functor $\Surf_{1} \rightarrow \Surf_{{\rm s}}$ which reverses the orientation of the $-$ boundary component and is compatible with union and glueing morphisms. If $E_{{\rm s}}:\,\Surf_{{\rm s}}  \rightarrow \mathcal{C}$ is a symmetric target theory, the composition of this forgetful functor with $E_{{\rm s}}$ gives a target theory $E$.  If $(A,B,C,D)$ is an $E_{{\rm s}}$-valued admissible initial data and $\Omega_{\Sigma}^{{\rm s}}$ the outcome of symmetric GR, then $(A,B,C,D)$ can be considered as an $E$-valued convergent initial data in the obvious way. And, for this convergent initial data, GR generates $\Omega_{\Sigma} = \Omega_{\Sigma}^{{\rm s}}$, as the defining formula \eqref{defGR} is identical to \eqref{defGRSym}. This little argument shows that symmetric GR is a particular case of GR where the initial data satisfies extra relations. Therefore, all the properties listed in \S~\ref{MainProp1}-\ref{MainProp2} also hold for symmetric GR.

\subsection{Strict symmetric GR}
\label{SSSS}
One can also formulate a symmetric version of the strict GR. The category $\overline{\Surf}_{{\rm s}}$ is constructed in the same way as $\overline{\Surf}_{1}$, except that we now define automorphisms in this category by choosing braiding representatives for each permutation of $\llbracket 1,n \rrbracket$, hence replace everywhere the groups $\mathfrak{S}_{\llbracket 2,n\rrbracket}$ with $\mathfrak{S}_{\llbracket 1,n\rrbracket} \cong \mathfrak{S}_{n}$. Strict symmetric target theories in this context are functors ${\rm e}:\,\overline{\Surf}_{{\rm s}} \longrightarrow \mathcal{C}$ satisfying axioms parallel to \S~\ref{S41}.

\begin{definition}\label{inininis} Initial data for strict symmetric GR is a quadruple $(A,B,C,D)$ where
\[
A \in {\rm e}(0,3)^{\mathfrak{S}_{3}},\qquad B \in {\rm e}(0,3),\qquad C \in {\rm e}(0,3)^{\mathfrak{S}_{\{2,3\}}},\qquad D \in {\rm e}(1,1).
\]
We also require that the initial data satisfies the four relations below. The first three ones involve glueing maps in $(0,4)$, the last one glueing maps in $(2,1)$, and $\sigma_{i,j}$ stands for the transposition of $i$ and $j$. 
\vspace{0.2cm}

\noindent \textsc{BA relation}. 
\[
(\sigma_{1,2} - {\rm Id})\big(\theta_{(0,3)_{2}}(B,A) + \theta_{(0,3)_{3}}(B,A) + \theta_{(0,3)_{4}}(B,A)\big) = 0.
\]

\vspace{0.2cm}

\noindent  \textsc{BB-AC relation}.
\[
(\sigma_{1,2} - {\rm Id})\big(\theta_{(0,3)_2}(B,B) +\theta_{(0,3)_3}(B,B) +\theta_{(0,3)_4}(C,A)\big) = 0.
\] 

\vspace{0.2cm}

\noindent  \textsc{BC relation}.
\[
(\sigma_{1,2} - {\rm Id})\big(\theta_{(0,3)_2}(B,C) + \theta_{(0,3)_3}(C,B) + \theta_{(0,3)_{4}}(C,B)\big) = 0.
\]

\vspace{0.2cm}

\noindent \textsc{D relation}. 
\[
(\sigma_{1,2} - {\rm Id})\big(\theta_{(1,1)}(B,D) + \tfrac{1}{2}\theta_{(0,3)}(C,A)\big) = 0.
\]
\end{definition}

Given a strict symmetric target theory and initial data, we define $w_{s} \in {\rm e}(s)$ for any object $s$ in $\overline{\Surf}_{{\rm s}}$ by the same formulas as in Definition~\ref{SDEFGR}.

\begin{proposition}
\label{thttsh}For any stable $(g,n)$, $w_{g,n}$ is a well-defined element of ${\rm e}(g,n)^{\mathfrak{S}_{n}}$. In other words, $(g,n) \mapsto w_{g,n}$ is a functorial assignment from $\overline{\Surf}_{{\rm s}}$.
\end{proposition}
\textbf{Proof.} The proof uses the same recollection of terms as done in the proof of functoriality of the symmetric GR in Theorem~\ref{SGRSYMTH}, with only (heavily simplifying) difference that the sums over homotopy class is replaced by one of our $(0,3)$-glueing maps, as is already apparent in comparing the 4 relations here to the 4 relations in \S~\ref{SGRSymIni}. 
\hfill $\Box$

\vspace{0.2cm}

As explained for GR/symmetric GR in \S~\ref{S33}, strict symmetric GR inherits  in the obvious way the properties of strict GR described in Section~\ref{S42}, upon replacing everywhere $\mathfrak{S}_{\llbracket 2,n\rrbracket }$ with $\mathfrak{S}_{n}$ allowing the permutation of all boundary components in $\Sigma_{g,n}$.

\providecommand{\bysame}{\leavevmode\hbox to3em{\hrulefill}\thinspace}
\providecommand{\MR}{\relax\ifhmode\unskip\space\fi MR }
\providecommand{\MRhref}[2]{%
  \href{http://www.ams.org/mathscinet-getitem?mr=#1}{#2}
}
\providecommand{\href}[2]{#2}

\end{document}